\documentclass[12pt,a4paper]{article}
\usepackage[12pt]{extsizes} 
\usepackage{amsfonts}
\usepackage{amssymb}
\usepackage{amsthm}
\usepackage{amsmath}
\usepackage{graphicx}
\usepackage{mathrsfs}
\usepackage{graphics}
\usepackage{abstract}
\usepackage{color}
\usepackage{cite}
\usepackage{physics}
\usepackage[shortlabels]{enumitem}
\usepackage{empheq}
\usepackage{appendix}
\usepackage{mathrsfs}
\usepackage{bm}
\usepackage{mathtools}
\mathtoolsset{showonlyrefs}
\usepackage[colorlinks=true,citecolor=red]{hyperref}
\hypersetup{linkcolor=blue}
\graphicspath{ {./result/} }
\usepackage{titling}
\usepackage{authblk}
\usepackage{fancyhdr}
\usepackage{lineno}
\usepackage{titlesec}
\titleformat{\section}[block]{\centering\Large\bfseries}{\thesection}{1em}{}
\titleformat{\subsection}[block]{\centering\large\bfseries}{\thesubsection}{1em}{}
\titleformat{\subsubsection}[block]{\centering\normalsize\bfseries}{\thesubsubsection}{1em}{}

\usepackage{ulem}

\renewcommand{\div}{\operatorname{div}}
\renewcommand{\ip}[2]{\left \langle #1, #2 \right \rangle}
\newcommand{\coloneqq}{\overset{\mathrm{def}}{=}}

\newtheorem{theorem}{\textbf{Theorem}}[section]
\newtheorem{lemma}{\textbf{Lemma}}[section]
\newtheorem{proposition}{\textbf{Proposition}}[section]
\newtheorem{corollary}{\textbf{Corollary}}[section]
\newtheorem{remark}{\textbf{Remark}}[section]
\newtheorem{definition}{\textbf{Definition}}[section]

\allowdisplaybreaks[4]

\newcommand{\psol}{(\mathbf{u},q)}

\newcommand{\app}{\mathrm{a}}

\newcommand{\dif}{\mathrm{d}}
\newcommand{\ini}{\mathrm{in}}
\newcommand{\R}{\mathrm{R}}

\numberwithin{equation}{section}

\renewcommand{\hat}[1]{\widehat{#1}}
\newcommand{\supp}{\mathop{\mathrm{supp}}}

\usepackage{tikz}
\usetikzlibrary{calc, arrows, intersections}
\title{\bf{Inertial instability of Couette flow with Coriolis force}}

\author{
 Yanlong Fan$^{a,}$\!\!
 \thanks{email:\protect\url{fanyanlong1@stu.scu.edu.cn}} \quad\quad 
 Daozhi Han$^{b,}$\!\!
\thanks{email:\protect\url{daozhiha@buffalo.edu} }
\quad\quad 
Quan Wang$^{c,}$\!\!
\thanks{email:\protect\url{xihujunzi@scu.edu.cn} } 
\\
\footnotesize $^{a}$ School of mathematics, Southwest Jiaotong University, Chengdu, Sichuan, 611756, China \\ 
\footnotesize $^{b}$ Department of Mathematics, The State University of New York at Buffalo, Buffalo, NY 14260 USA \\ 
\footnotesize $^{c}$ College of Mathematics, Sichuan University, Chengdu, Sichuan, 610065,  China
\\  
} 

\date{\today}

\begin{document}
\maketitle
\begin{abstract} 
    We analyze the nonlinear inertial instability of Couette flow under Coriolis forcing in  \(\mathbb{R}^{3}\). For the Coriolis coefficient \(f \in (0,1)\), we show that the non-normal operator associated with the linearized system admits only continuous spectrum. Hence, there are no exponentially growing eigenfunctions for the linearized system. Instead, we construct unstable solutions in the form of pseudo-eigenfunctions that exhibit non-ideal spectral properties. By identifying the growth bound of the linear semigroup with the mode rate \(\Lambda=\sqrt{f(1-f)}\) and by a bootstrap argument that resolves the challenges posed by the non-ideal spectral behavior of pseudo-eigenfunctions, we establish the velocity instability of Couette flow in the Hadamard sense for every \(f\in(0,1)\).
    
\end{abstract}


\newpage
\tableofcontents

\newpage
\section{Introduction}  
\subsection{Presentation of the problem}\label{subsec_problem}
Inertial instability of shear flows is a crucial concept in fluid dynamics, manifesting in diverse natural and engineered fluid systems. It occurs when the equilibrium between the pressure gradient and the Coriolis force in a rotating fluid is disrupted. 
The inertial instability of shear flows drives the redistribution of momentum, heat, and mass within the fluid. It can give rise to complex flow patterns and vortices, affecting the overall energy balance and dynamic processes of the system. Understanding this instability is essential for accurately predicting and explaining various phenomena in fields ranging from climate science \cite{Thompson2021,Rapp2018,Grisouard2020} to astrophysics \cite{Park2020,Park2021}.

In this article, we consider inertial instability governed by the Navier--Stokes equations with Coriolis force in \(\mathbb{R} ^3 \)
\begin{align}
	\begin{cases}
		\frac{\partial{\mathbf{v}}}{\partial{t}}+ \mathbf{v} \cdot \nabla \mathbf{v} + \mathbf{f} \times \mathbf{v} = \nu \Delta \mathbf{v} - \nabla p  , \\
		\div   \mathbf{v}= 0,
	\end{cases}
\end{align}
where \( \mathbf{v} = (v_1, v_2, v_3) \) is the velocity, \( p \) is the pressure, \( \mathbf{f} = (0, 0, f )\) is the Coriolis vector and \( f>0 \) is called the Coriolis coefficient. A steady shear flow takes the form of  
\begin{equation}\label{steady_state}
	 \overline{\mathbf{v}}=\big(U(y),0,0\big),\quad  \overline{p} =\overline{p}(y),
\end{equation}
where \( U \) and \( \overline{p} \) satisfy the geostrophic balance
\begin{align}
	f  U(y) = - \frac{\partial{\overline{p}}}{\partial{y}}.
\end{align}
The inertial instability of shear flows is determined by the Rayleigh discriminant \cite{Park2020}, also known as the Bradshaw--Richardson number \cite{Bradshaw1969,Huang2018}, defined as
\[ \mathrm{R}(y;f,U(y)) \coloneqq f (f-U'(y)). \]
When $\mathrm{R}(y;f,U(y)) <0$ at some \( y_0 \), the flow becomes inertially unstable. In this article, we consider  the Couette flow
\begin{equation}\label{cff}
    U(y) = y.
\end{equation}
Thus, the Rayleigh discriminant becomes a constant
$
    \mathrm{R}(f) = \mathrm{R}(y;f,y)=f (f- 1).
$

Couette flow is a phenomenon where viscous fluids move between parallel plates or concentric cylinders, driven by viscous shear forces as one surface moves tangentially. The Couette flow along with other monotonic shear flows are common in geophysical fluids and have been comprehensively studied from both the physical \cite{Avila2023,Orlandi2015,Tillmark1992,Grossmann2016,Hiwatashi2007} and the mathematical perspectives \cite{Chen2020,Bedrossian2015,Chen2022,Masmoudi2022,Masmoudi2024,Ma2010,Castro2023,Cui2025,Huang2024b}. As the Couette flow is widespread within geophysical fluids, the primary factor leading to the instability of this flow is the Coriolis force. This type of instability is geophysically 
referred to as inertial instability \cite{DJ2019}.

Introducing the perturbations around the Couette flow \eqref{cff} 
\begin{equation}
		\mathbf{u} = (u_1, u_2 ,u_3) = \mathbf{v} - \overline{\mathbf{v}}, \quad 
		q = p - \overline{p}, 
\end{equation}
we derive the perturbed equations
\begin{equation}\label{nonlinear_perturbed_equation}
	\begin{cases}
		\frac{\partial{\mathbf{u}}}{\partial{t}} + \mathbf{u} \cdot  \nabla \mathbf{u} + 
        \begin{pmatrix}
            u_2 \\0 \\0
        \end{pmatrix} 
        + y \partial_x  \mathbf{u}
		+ f 
         \begin{pmatrix}
            -u_2 \\ u_1 \\0
        \end{pmatrix} 
        = \nu \Delta \mathbf{u}-  \nabla q  , \\ 
		\div   \mathbf{u}= 0,
	\end{cases}
\end{equation}
where the pressure \( q \) is given by \(  q \coloneqq q^{L} + q^{NL}  \) and 
\begin{equation}
    \begin{aligned}
    q^{L} &\coloneqq - \Delta ^{-1} \div \qty( \mathbf{u} \cdot \nabla \overline{\mathbf{v}} +\overline{\mathbf{v}}\cdot \nabla \mathbf{u}
		+ \mathbf{f} \times \mathbf{u}) \\&= - \Delta ^{-1} \qty( 2\partial_x u_2  + f ( \partial_y u_1 - \partial_x u_2 )), \\ 
        q^{NL} &\coloneqq - \Delta ^{-1} \div \qty(\mathbf{u} \cdot \nabla \mathbf{u}).
    \end{aligned}
\end{equation}
The system is supplemented with the initial data
\begin{equation}\label{initial_conditions_for_perturbed_equation}
	\left. \mathbf{u}(t,\mathbf{x}) \right|_{t = 0} = \mathbf{u}_{\ini}(\mathbf{x}) \coloneqq (u_{1,\ini}(\mathbf{x}),u_{2,\ini}(\mathbf{x}),u_{3,\ini}(\mathbf{x}))
\end{equation}
and the far-field condition at infinity
\begin{equation}\label{boundary_conditions_for_perturbed_equation}
	\lim_{ \abs{\mathbf{x}}  \to + \infty} \abs{\mathbf{u}(t,\mathbf{x})}  = 0, \quad \forall t\geq 0.
\end{equation}

Specifically, we aim to establish the nonlinear inertial instability in the following sense: 
 \begin{definition}[Nonlinear instability in the sense of Hadamard]\label{def_linear_instability}
    The steady state solution \( \mathbf{\overline{v}} \) is nonlinearly unstable  if there are constants $\sigma$ and \( C \) such that for every $\delta$ arbitrarily small there exists a solution ${\bf u}$ of \eqref{nonlinear_perturbed_equation} satisfying
    \begin{align*}
    &||{\bf u}(0,\mathbf{x}) ||_{H^{1}}  
    \leq\delta,\quad  ||\mathbf{u}(T^{\delta},\mathbf{x}) ||_{L^{2}}  \geq\sigma,
    \end{align*}
    where $T^{\delta}\leq C\ln \delta^{-1}  + C$ is an escape time.
    \end{definition}

    \autoref{def_linear_instability} implies the failure of Hadamard well-posedness: it violates the uniform continuous dependence of the solution on the initial data, in the sense that the escape time $T^{\delta}$ depends on $\delta$ and diverges as $\delta\to0$. We refer to \cite{hwang2003,mao2024,Jiang2014,Nguyen2024} among others for the analysis of Rayleigh--Taylor instability in the sense of Hadamard.

\subsection{Roadmap and main results}

The linear part of system \eqref{nonlinear_perturbed_equation} is identified as
\begin{equation}\label{linearized_perturbed_equation}
	\begin{cases}
		\frac{\partial{\mathbf{u}}}{\partial{t}}  + 
        \begin{pmatrix}
            u_2 \\0 \\0
        \end{pmatrix} 
        + y \partial_x  \mathbf{u}
		+ f 
         \begin{pmatrix}
            -u_2 \\ u_1 \\0
        \end{pmatrix} 
        = \nu \Delta \mathbf{u} + \nabla \Delta ^{-1} \qty( 2\partial_x u_2  + f ( \partial_y u_1 - \partial_x u_2 ))  , \\ 
		\div   \mathbf{u}= 0,
	\end{cases}
\end{equation}
with the initial-boundary conditions \eqref{initial_conditions_for_perturbed_equation}-\eqref{boundary_conditions_for_perturbed_equation}.  A key step in the analysis is to obtain the spectral properties of the linearized operator $\mathscr{L}$
\begin{equation}
    \label{Linear_eq_abstract_form} 
    \partial_t \mathbf{u} = \mathscr{L} \mathbf{u} , 
\end{equation}
with the domain of definition 
\begin{equation}\label{domain_of_L}
	\begin{aligned}
		\mathcal{D} (\mathscr{L}) & \coloneqq  \qty{ \mathbf{u} =(u_1,u_2,u_3) \in[L^2(\mathbb{R}^3)]^3 \ \middle|  \
		 \nu \Delta\mathbf{u} , y \partial_x \mathbf{u} \in[L^2(\mathbb{R}^3)]^3 }  .                                \\
	\end{aligned}
\end{equation}
We split the operator \( \mathscr{L} \) as follows 
\begin{equation}\label{definition_of_L}
   \begin{aligned}
	 \mathscr{L} &\coloneqq \mathscr{L}_0  + \mathscr{L}_1 + \mathscr{L}_{2}, \quad  \mathscr{L}_{2}  \coloneqq - \mathop{\mathrm{diag}} (y \partial_x, y \partial_x, y \partial_x), \\ 
    \mathscr{L}_0 & \coloneqq 
    \begin{pmatrix}
        \nu \Delta  & f - 1 & 0 \\ 
       -f + f \Delta^{-1}  \partial_{y}^2 & \nu \Delta & 0\\ 
       f \Delta^{-1}  \partial_{yz} & 0 &\nu \Delta  \\
    \end{pmatrix} , \\
    \mathscr{L}_1 & \coloneqq  \Delta ^{-1}
    \begin{pmatrix}
       f  \partial_{xy}   & \qty(2 - f )\partial_x^2  & 0 \\
        0  & \qty(2 - f )\partial_{xy}  & 0 \\
         0 & \qty(2 - f )\partial_{xz}  & 0 \\
    \end{pmatrix}.  
\end{aligned}
\end{equation}
That is
\begin{itemize}
    \item \( \mathscr{L}_0 \) denotes the component whose form remains invariant when \( \mathbf{u} \) is independent of \( x \);
    \item \( \mathscr{L}_1 \) is the bounded part of the operator involving partial derivatives with respect to \( x \);
    \item \( \mathscr{L}_2 \) is the bad part with \( \partial_x \) that renders the operator non-normal.
\end{itemize}
Note that \(\mathscr{L}_0\) is a normal operator for which the analysis is standard, cf. {\rm\cite[{\rm Section 5.3}]{Satoh2014}}.

\begin{remark}
     As observed in {\rm\cite{Bedrossian2020,Zeng2022}} and references therein, the unknowns \( w_{i} = \Delta u_{i}\), \(  i = 1,2,3 \) are more convenient to use than \( \mathbf{u} \) in the study of stability threshold, see also   {\rm\cite{Li2023}}. In this article, however, the main focus is on nonlinear instability of \(\mathbf{u}\).
\end{remark}

For nonlinear instability analysis, one seeks the linearly growing modes and a bound on the growth rate of the linear equation, cf. \cite{Jiang2014,Wang2012}. However, the growth rate of the linear operator does not need to be optimal as long as it is less than twice the growth rate of the linearly growing modes. This is sufficient to overcome the effect of the nonlinear terms, thereby proving nonlinear instability in the Hadamard sense. Crude bounds may still miss even this factor-two threshold.

In what follows, we define the exponential growth rate 
\begin{equation}
    \Lambda \coloneqq \sqrt{f( 1 - f)}.
\end{equation}
The quantity \( \Lambda \) is the growth rate of the linearly growing modes constructed below, and is essentially sharp as a uniform bound for the linear semigroup, up to an arbitrarily small excess. Precisely,
\begin{theorem}\label{thm:simp_growth}
Let $0<f<1$ and $\nu>0$. Then $\omega_{0}(\mathscr{L})=\Lambda$. 
Moreover, for every $\varepsilon>0$ there exists a constant $C_{\varepsilon}>0$ depending on $\varepsilon$ such that
\begin{equation}\label{eq:simp_semigroup_bnd}
\norm{e^{t\mathscr{L}}}_{L^{2}\to L^{2}}\le C_{\varepsilon}e^{(\Lambda+\varepsilon)t}
\qq{for all} t\ge0.
\end{equation}
Furthermore, for every $\varepsilon>0$,
\begin{equation}\label{eq:simp_est}
\norm{(\lambda-\mathscr{L})^{-1}}_{L^2 \to L^2}\le C\varepsilon^{-1}\bigl(1+\varepsilon^{-1}\bigr)e^{C\varepsilon^{-1}}
\qq{for any $\lambda$ with}  \operatorname{Re}\lambda\ge\Lambda+\varepsilon.
\end{equation}
\end{theorem}

The main difficulties in establishing \autoref{thm:simp_growth} are as follows. For bounded profiles the Kato--Pazy theorem does give analyticity when the perturbation is \(A\)-bounded with relative bound \(0\) \cite[Chap.~3, Thm.~2.1]{Pazy1992}, and then spectral and growth bounds coincide as in Friedlander--Pavlovi\'c--Shvydkoy \cite{Friedlander2006}. This hypothesis fails for the unbounded \(y\partial_{x}\), and indeed \(\mathscr{L}_{*}=\nu\Delta-y\partial_{x}\) generates only a \(C_{0}\)-group, since the drift group \(e^{t y\partial_{x}}\) is unitary and non-analytic. Thus, the resolvent of \(\mathscr{L}_{*}\) does not satisfy the sectorial estimate, so \(\mathscr{L}\) is not analytic. 

Crude bounds on the semigroup growth, such as those obtained from the numerical range of \(\mathscr{L}\), fall short of the sharp rate \(\Lambda\); even the relaxed requirement of a growth bound strictly less than \(2\Lambda\), which would already suffice for the nonlinear argument, remains out of reach for some \( f\in (0,1) \). 

To address this obstacle, the spectral properties of \(\mathscr{L}\) are approached via a Fourier analogue of Kirchgässner's reduction~\cite{Kirchgassner1982,Kirchgassner1988}. In the classical setting one treats a distinguished spatial coordinate as an evolutionary variable, thereby converting an elliptic problem into an ODE in that ``time''. As a toy model, consider the scalar operator \(\mathscr{L}_{\star}=\nu\Delta-y\partial_{x}\). Its resolvent problem \((\lambda-\mathscr{L}_{\star})u=F\) becomes, after Fourier transformation,
\begin{equation}
\xi_{1}\partial_{\xi_{2}}\widehat{u}=(\lambda+\nu|\bm{\xi}|^{2})\widehat{u}-\widehat{F},
\end{equation}
an inhomogeneous ODE in the frequency \(\xi_{2}\) whenever \(\xi_{1}\neq0\). For the full operator \(\mathscr{L}\), the Couette transport likewise becomes \(\xi_{1}\partial_{\xi_{2}}\), so that for fixed \((\xi_{1},\xi_{3})\) with \(\xi_{1}\neq0\) both the eigenvalue problem and the resolvent equation reduce to ODEs in \(\xi_{2}\). This ODE viewpoint yields the absence of point spectrum \autoref{lemma_no_point_spectrum} and underpins the resolvent bound in \autoref{thm:simp_growth}.

For the construction of linearly unstable modes, one notes that the linear system \eqref{linearized_perturbed_equation} neither possesses a natural variational structure nor admits an integrable separation-of-variables solution in \( L^2 \), cf. \cite{Jiang2014,Fan2025,Jiang2014a,Guo2011,Nguyen2024} for applications of the variational method. Nonetheless, we observe that the zonal-independent variant of the equations \eqref{linearized_perturbed_equation} admits separation-of-variables solutions, though these solutions are not integrable in \( \mathbb{R}^3 \) and thus not eigenfunctions of the linear operator. Indeed, as shown in \autoref{lemma_no_point_spectrum}, the linear operator does not have any point spectrum--this implies that it is impossible to find an exponentially growing eigenfunction for the linearized system. To tackle this problem, we take a perturbative approach and look for the linear solution near the exponentially growing solution of the zonal-independent version of the equations \eqref{linearized_perturbed_equation}. Precisely, we have the following.
\begin{theorem}\label{thm:Linear_Instability}
	Suppose \( 0 < f < 1 \), then for any given time \( T > 0 \) and \( \epsilon > 0 \), equations \eqref{linearized_perturbed_equation} with initial-boundary conditions \eqref{initial_conditions_for_perturbed_equation}-\eqref{boundary_conditions_for_perturbed_equation} has a real-valued solution \( \mathbf{u}^{T,\epsilon}(t) \in C([0,T];[H^k(\mathbb{R}^3)]^3) \) with \(\div \mathbf{u}^{T,\epsilon}=0\) and  \( \nabla q^{T,\epsilon} (t) \in C([0,T];[H^k(\mathbb{R}^3)]^3) \), \( k \geq 0 \) satisfying
	\begin{enumerate}[{\rm(i)}]
		\item For any \( t\in [0,T] \) 
              \begin{equation}
                \begin{aligned}
                    &  (e^{t \Lambda } - \epsilon) \norm{ \mathbf{u}^{T,\epsilon}(0)}_{L^2}   \leq   \norm{\mathbf{u}^{T,\epsilon}(t)}_{L^2}\leq (e^{t \Lambda  } + \epsilon) \norm{\mathbf{u}^{T,\epsilon}(0)}_{L^2}; \\
\end{aligned}
            \end{equation}   
\item   
For any positive integer \( k \), there exists a constant \( C_k   > 0  \) independent of \( T \) and \( \epsilon \) such that
\begin{equation}
    \norm{\mathbf{u}^{T,\epsilon}(0)}_{H^k} \leq C_k \norm{\mathbf{u}^{T,\epsilon}(0)}_{L^2}.   
\end{equation}
	\end{enumerate}
\end{theorem}

 The nonlinear instability of the perturbed system \eqref{nonlinear_perturbed_equation} is established in the following theorem.
   \begin{theorem}\label{thm:nonlinear_instability_Hadamard}
    The Couette flow is unstable in the sense of Hadamard. That is, for any \(0 < f < 1\), \( \nu > 0 \), there exist constants \( \overline{\delta}_0 \leq 1 \), \( \varepsilon_0 > 0 \) such that for any \( \delta \in (0,\overline{\delta}_0) \) there exists an initial condition \( \mathbf{u}_{\ini}^{\delta} \in [\mathcal{S}(\mathbb{R}^3)]^3\) with \(\div \mathbf{u}_{\ini}^{\delta}=0\) and \(\norm{\mathbf{u}_{\ini}^{\delta}}_{H^{1}}=\delta\), and a unique strong solution \( \mathbf{u}^{\delta} \in C([0,T^{\max});[H^1(\mathbb{R}^3)]^3)\cap L^2(0,T^{\max};[H^2(\mathbb{R}^3)]^3) \) with \(\div \mathbf{u}^{\delta}=0\) of \eqref{nonlinear_perturbed_equation}-\eqref{boundary_conditions_for_perturbed_equation} emanating from \(\mathbf{u}_{\ini}^{\delta}\) with an associated pressure \( q^{\delta} \in C([0, T^{\max});H^{1}(\mathbb{R}^3)) \), such that
    \begin{equation}
        \begin{aligned}
            \norm{\mathbf{u}^{\delta}(T^{\delta})}_{L^2} > \varepsilon_0
        \end{aligned}
    \end{equation}
    for some escape time \( T^{\delta} =\frac{1}{\Lambda} \ln \frac{2\varepsilon_0}{L\delta} < T^{\max}\), where \( L>0 \) is the constant and \( T^{\max} \) denotes the maximal time of existence of \( \mathbf{u}^{\delta} \). 

\end{theorem} 

Note that the linearly unstable solution is given by a pseudo-eigenfunction, which does not behave like an exponential function for large time. Indeed, as time increases, the error between the pseudo-eigenfunction and an exponential function accumulates. Only for a fixed time \( T \) can one find an unstable solution that is uniformly close to an exponential on \( [0,T] \). To establish the nonlinear instability result, one needs a precise estimate of the dependence of $T$ on the size of the initial condition, and a careful construction of the initial condition. In particular the admissible initial datum $\mathbf{u}_{\ini}^{\delta}$ in \autoref{thm:nonlinear_instability_Hadamard} is allowed to depend on $\delta$, exactly as in \autoref{def_linear_instability}, via the choice $T\geq T^{\delta}\sim\Lambda^{-1}\ln(\delta^{-1})$; this is unavoidable when no genuine eigenfunction exists and is consistent with the Hadamard formulation.

 \subsection{Literature review}
Without the Coriolis effect \( f = 0 \), the linearized system around Couette flow is known to be spectrally stable, in the sense that there are no unstable eigenmodes. However, it has been observed that the linear system has large pseudo-spectra and may lead to significant transient growth \cite{Trefethen1993}. In this case, the study of nonlinear stability, namely the transition threshold, of 3D Couette flow is carried out in some recent works \cite{Bedrossian2017,Bedrossian2020,Bedrossian2022a, Wei2021,Chen2024}. For nonlinear instability, Li et al. \cite{Li2023} adopted a dynamical approach to study the instability of steady-state profiles near the Couette flow, and showed that the vorticity field is unstable in the $L^2$ norm. Recently, boundary-driven instability \cite{Bian2025} and viscosity-driven instability \cite{Li2025} of shear flows have also been investigated. 

The nonlinear stability of Couette flow with Coriolis force is established in \cite{Zelati2025,Huang2024} when \( \R(f) > 0 \). In particular,  Guo et al. \cite{Guo2022} constructed a class of axisymmetric global solutions near the stationary state to the 3d Euler equations with uniform rigid body rotation \(f = 1\) (\( \R(f)= 0 \)). To the best of our knowledge, when \( \R(f) < 0 \), the inertial instability of the Couette flow in $\mathbb{R}^3$ has not been resolved in the literature, which is the major undertaking of this article. Recent results regarding the effect of the Coriolis force on the stability and instability of Couette flow are summarized in the table below.
\begin{table}[htbp]
  \centering
  \caption{The effect of Coriolis coefficient \( f \) on (in)stability of Couette flow}
  \label{tab:linear_effects}
  \begin{tabular}{c|c|c|c}
    \hline
    References    & domain           & range of \( f \) &  small perturbation is                 \\
    \hline
       \cite{Bedrossian2017,Wei2021}  & \( \mathbb{T} \times \mathbb{R} \times \mathbb{T} \)     & \( f = 0 \)                & nonlinear stable                      \\
   \cite{Huang2024}            & \( \mathbb{T} \times \mathbb{R} \times \mathbb{T} \)     & \( ( - \infty, - 2] \cup [2, + \infty)    \)                & nonlinear stable                      \\
      \cite{Huang2024a}            & \( \mathbb{T} \times \mathbb{R} \times \mathbb{T} \)     & \( f = 1 \)                & nonlinear stable                      \\ 
   \cite{Zelati2025}        & \( \mathbb{T} \times \mathbb{R} \times \mathbb{T} \)        & \( ( - \infty,0) \cup (1, + \infty)\)            & nonlinear stable       \\
    Current paper    & \( \mathbb{R}^3 \)        & \( (0,1) \)   & nonlinear unstable             \\ 
    \hline
  \end{tabular}
\end{table}

From the linear analysis, it is expected that the Couette flow is nonlinearly unstable when \(0 < f < 1\), while remaining stable outside this interval. Combining the linearly growing modes of \autoref{thm:Linear_Instability} with the sharp growth bound \(\omega_{0}(\mathscr{L})=\Lambda\) of \autoref{thm:simp_growth}, we obtain nonlinear Hadamard instability for every \(f\in(0,1)\). This means that deriving results regarding the nonlinear stability threshold is infeasible. It also implies that the stabilizing effect, such as inviscid damping and enhanced dissipation, fails to suppress the linear instability effect generated by the Coriolis force.


\subsection{Organization of the paper}
The article is structured as follows. \autoref{sec_Preliminary} recalls some preliminary concepts. \autoref{sec_linear_analysis} establishes the semigroup and spectral properties of \(\mathscr{L}\) and proves the resolvent estimate of \autoref{thm:simp_growth}. \autoref{sec_linear_instability} constructs the linearly unstable modes of \autoref{thm:Linear_Instability}. \autoref{sec_nonlinear_instability} proves the nonlinear instability \autoref{thm:nonlinear_instability_Hadamard}. The appendix provides some analysis tools.

\section{Preliminary}\label{sec_Preliminary}
Throughout, \( \ip{\cdot}{\cdot} \) is the \( L^{2} \) inner product. \( \cdot^{\intercal} \) denotes the transpose of a matrix or a vector. The space \(\mathcal{H}=[L^{2}(\mathbb{R}^{3})]^{3}\). We denote the divergence-free subspaces by
\begin{equation}
    \begin{aligned}
        L^{2}_{\sigma}(\mathbb{R}^{3}) &=\{\mathbf u\in[L^{2}(\mathbb{R}^{3})]^{3}:\div\mathbf u=0\},\\
        H^{k}_{\sigma}(\mathbb{R}^{3}) &=\{\mathbf u\in[H^{k}(\mathbb{R}^{3})]^{3}:\div\mathbf u=0\},\\
        \mathcal S_{\sigma}(\mathbb{R}^{3}) &= [\mathcal S(\mathbb{R}^{3})]^{3}\cap L^{2}_{\sigma}(\mathbb{R}^{3}),
    \end{aligned}
\end{equation}
understood in the sense of distributions.

Suppose \( \mathscr{T} :\mathcal{H}_1  \to \mathcal{H}_2 \) is an operator from the Hilbert space \( \mathcal{H}_1 \) to another Hilbert space \( \mathcal{H}_2 \) in which the linear subspace \( \mathcal{D}(\mathscr{T}) \) is the domain of \( \mathscr{T} \). We denote the range of \( \mathscr{T} \) by
\begin{equation*}
	\mathcal{R}(\mathscr{T}) \coloneqq \mathscr{T}(\mathcal{D}(\mathscr{T})) = \qty{\mathscr{T}(x) \ | \ x\in \mathcal{D}(\mathscr{T})  }.
\end{equation*}
We write the kernel of \( \mathscr{T} \) as
\begin{equation}
	\mathcal{N}(\mathscr{T}) \coloneqq \qty{x \in \mathcal{D}(\mathscr{T}) \ | \ \mathscr{T}(x) = 0  }.
\end{equation}
 The graph of  \( \mathscr{T} \)  is the set
\[
	\mathcal{G}(\mathscr{T} ) =  \qty{(x, \mathscr{T} x)\ | \ x \in \mathcal{D}(\mathscr{T} )}.
\]
The norm \( \norm{x}_{\mathscr{T}} = \norm{x}_{\mathcal{H}_1} + \norm{\mathscr{T}x}_{\mathcal{H}_2} \) on \( \mathcal{D}(\mathscr{T}) \) is called the graph norm of the operator \( \mathscr{T} \).
\begin{definition}
	An operator \( \mathscr{T} \) is called closed if its graph \( \mathcal{G}(\mathscr{T} ) \) is a closed subset of the Hilbert space \( \mathcal{H}_1 \times \mathcal{H}_2 \), and \( \mathscr{T} \) is called closable (or pre-closed) if there exists a closed linear operator \( \mathscr{\overline{T}} \) from \( \mathcal{H}_1 \) to \( \mathcal{H}_2 \) such that \( \mathscr{T} \subseteq \mathscr{\overline{T}} \). The operator \( \mathscr{\overline{T}} \) is called the closure of the closable operator \( \mathscr{T}\).
\end{definition}
 
Another useful notion is that of a core of an operator which allows us to prove statements of a closed operator on its core rather than the domain.
\begin{definition}
	A linear subspace \( \mathcal{D} \) of \( \mathcal{D}(\mathscr{T}) \) is called a core for \( \mathscr{T} \) if \( \mathcal{D} \) is dense in \( (\mathcal{D}(\mathscr{T}), \norm{\cdot}_{\mathscr{T}}   ) \), that is, for each \( x\in  \mathcal{D}(\mathscr{T}) \), there exists a sequence \( (x_{n})_{n\in \mathbb{N}} \) of vectors \( x_n \in \mathcal{D} \) such that \( x = \lim_{n \to \infty} x_n\) in \( \mathcal{H}_1 \) and \( \mathscr{T}x = \lim_{n \to \infty} \mathscr{T}x_n\) in \( \mathcal{H}_2 \).
\end{definition}

Let us define \( \rho(\mathscr{T}) \) and \( \sigma(\mathscr{T}) \) to be the resolvent set and the spectrum of a closed operator \( \mathscr{T} \) on a Banach space \( \mathcal{X} \), respectively. The spectral bound is defined by
$
	\alpha(\mathscr{T}) = \sup_{\lambda \in \sigma(\mathscr{T})} \Re \lambda.
$
\begin{theorem}[\cite{Trefethen2005} p.150]\label{theorem_spectral_instability}
	Let \( \mathscr{T} \) be a closed linear operator on a Banach space \( \mathcal{X} \) generating a \( C_0 \)-semigroup. Then
	\begin{equation}
		\norm{e^{t \mathscr{T}}} \geq e^{t \alpha(\mathscr{T})}  ,\quad \forall t \geq 0,
	\end{equation}
	where we use \( e^{t \mathscr{T}} \) to denote the \( C_0 \)-semigroup generated by \( \mathscr{T} \).
\end{theorem}

The following concept of pseudo-spectra is from the monograph by Lloyd N. Trefethen and Mark Embree
\cite{Trefethen2005}.
\begin{definition}[\cite{Trefethen2005} p.31]\label{pspectral_definition}
	Let \( \mathscr{T} \) be a closed operator on a Banach space \( \mathcal{X} \) and \( \varepsilon > 0\) be arbitrary. The \( \varepsilon \)-pseudo-spectrum \( \sigma_{\varepsilon}(\mathscr{T}) \) of \( \mathscr{T} \) is the set of \( \zeta\in \mathbb{C} \) defined equivalently by any of the conditions
	\begin{enumerate}[(i)]
		\item \( \norm{(\zeta - \mathscr{T})^{-1} } > \varepsilon^{-1}    \);
		\item \( \zeta \in \sigma(\mathscr{T} + \mathscr{E}) \) for some bounded operator \( \mathscr{E} \) with \( \norm{\mathscr{E}} < \varepsilon  \);
		\item \( \zeta\in \sigma(\mathscr{T}) \) or \( \norm{(\zeta - \mathscr{T}) u} < \varepsilon \) for some \( u \in \mathcal{D}(\mathscr{T}) \) with \( \norm{u} = 1   \). Then \( \zeta \) is an \( \varepsilon \)-pseudo-eigenvalue of \( \mathscr{T} \) and \( u \) is the corresponding \( \varepsilon \)-pseudo-eigenfunction.
	\end{enumerate}
\end{definition}
The \( \varepsilon \)-spectral bound is defined as
$
	\alpha_{\varepsilon}(\mathscr{T}) = \sup_{\lambda \in \sigma_{\varepsilon}(\mathscr{T})} \Re \lambda.
$
\begin{theorem}[\cite{Trefethen2005} p.31]\label{theorem_pspectral_properties}
	Given a closed operator \( \mathscr{T} \) on a Banach space \( \mathcal{X} \), the pseudo-spectra \( \{ \sigma_{\varepsilon}(\mathscr{T})\}_{\varepsilon> 0} \) have the following properties. They can be defined equivalently by any of the conditions (i)-(iii) in \autoref{pspectral_definition}. Each \( \sigma_{\varepsilon}(\mathscr{T}) \) is a nonempty open subset of \( \mathbb{C} \), and any bounded connected component of \( \sigma_{\varepsilon}(\mathscr{T}) \) has a nonempty intersection with \( \sigma(\mathscr{T}) \). The pseudo-spectra are strictly nested supersets of the spectrum: \( \cap_{\varepsilon > 0 }\sigma_{\varepsilon}(\mathscr{T}) =\sigma(\mathscr{T}) \), and conversely, for any \( \delta > 0 \), \( \sigma_{\varepsilon + \delta}(\mathscr{T}) \supseteq \sigma_{\varepsilon }(\mathscr{T}) + B(\delta)\); where \( B(\delta) \) is the open disk of radius \( \delta \).
\end{theorem}

Finally, we recall the Gearhart--Pr\"uss theorem, which characterizes the growth bound of a bounded \( C_{0} \)-semigroup by means of the resolvent; we refer to {\rm\cite[Theorem V.1.11]{Engel2008}} for the precise statement.
\begin{theorem}[Gearhart--Pr\"uss]\label{thm:gearhart_pruss}
	Let \( \mathscr{T} \) generate a bounded \( C_{0} \)-semigroup \( e^{t\mathscr{T}} \) on a Hilbert space \( \mathcal{X} \), with growth bound
	\begin{equation}\label{eq:gap_growth_bound}
		\omega_{0}(\mathscr{T}) \coloneqq \lim_{t\to\infty}\frac1t\log\norm{e^{t\mathscr{T}}}_{\mathcal{X}\to\mathcal{X}} ;
	\end{equation}
	then
	\begin{equation}\label{eq:gap_identity}
		\omega_{0}(\mathscr{T})
		=\inf\Bigl\{\omega\in\mathbb R:\sup_{\Re\lambda>\omega}\norm{(\lambda-\mathscr T)^{-1}}<\infty\Bigr\}.
	\end{equation}
\end{theorem}

\section{The linear analysis}\label{sec_linear_analysis}

In this section, we analyze the linearized operator \( \mathscr{L} \) defined in \eqref{Linear_eq_abstract_form}--\eqref{definition_of_L}.

\subsection{The semigroup and spectral analysis}
In this subsection, we prove that the operator \( \mathscr{L} \) generates a \( C_0 \)-semigroup on \( [L^2(\mathbb{R}^3)]^3 \), and that its spectrum consists entirely of continuous spectrum. For the sake of convenience, we define the key part of \( \mathscr{L} \) by
\begin{equation}
    \mathscr{L}_{\star} \coloneqq \nu \Delta - y \partial_x,
\end{equation}
the domain of which is defined by 
\begin{equation}\label{domain_of_L_star}
    \mathcal{D}(\mathscr{L}_{\star})  \coloneqq   \qty{u\in L^2(\mathbb{R}^3)\ | \ \Delta u, y\partial_x u \in L^2(\mathbb{R}^3)}  . 
\end{equation}
One may recall from \eqref{definition_of_L} that we can rewrite the definition of \( \mathscr{L} \) by
\begin{equation}\label{decomposition_of_L}
    \mathscr{L} = \mathscr{L}_{\star }I_3 + \mathscr{L}_{\bullet } , 
\end{equation} 
where \( I_3 \) is the identity matrix of order 3 and 
\begin{equation}
    \mathscr{L}_{\bullet} \coloneqq  \begin{pmatrix}
       f \Delta^{-1}  \partial_{xy}   & f - 1 + \qty(2 - f ) \Delta^{-1}\partial_x^2  & 0 \\
       - f + f \Delta^{-1}  \partial_{y}^2  & \qty(2 - f ) \Delta^{-1}\partial_{xy}  & 0 \\
         f \Delta^{-1}  \partial_{yz}  & \qty(2 - f ) \Delta^{-1}\partial_{xz}  & 0 \\
    \end{pmatrix} , 
\end{equation}
is a bounded operator on \( [L^2(\mathbb{R}^3)]^3 \).

\begin{lemma}\label{lemma:domain_of_L}
    The operator \( \mathscr{L} \) is closed with the domain \( \mathcal{D}(\mathscr{L}) \) given by \eqref{domain_of_L}. 
\end{lemma}
\begin{proof}
    First, note that \( \mathscr{L}_{\star} \) with domain \( C_0^{\infty}(\mathbb{R}^3) \) is closable (cf. \cite{Schmuedgen2012}). For any \( u \in C_0^{ \infty}(\mathbb{R}^3) \)  we have, using integration by parts, that 
\begin{equation}\label{sec_1:eq_1}
    \begin{aligned}
	\norm{\mathscr{L}_{\star} u}_{L^2}^2 & =  \norm{\nu \Delta u}_{L^2}^2 + \norm{y \partial_x u}_{L^2}^2 - 2   \ip{\nu \Delta u}{y \partial_x u} \\ 
    & = \norm{\nu \Delta u}_{L^2}^2 + \norm{y \partial_x u}_{L^2}^2 + 2 \nu \int_{\mathbb{R}^3}^{} y \partial_x \nabla u \cdot\nabla u
    +\partial_x u \partial_y u
    \dd{\mathbf{x}}  \\ 
    & = \norm{\nu \Delta u}_{L^2}^2 + \norm{y \partial_x u}_{L^2}^2 + 2 \nu \int_{\mathbb{R}^3}^{} \partial_x u \partial_y u
    \dd{\mathbf{x}} . \\ 
\end{aligned}
\end{equation}
Now, in view of Plancherel's formula, one notices that 
\begin{equation}
    \begin{aligned}
	 &\abs{\int_{\mathbb{R}^3}^{} \partial_x u \partial_y u
    \dd{\mathbf{x}}}  = \abs{\int_{\mathbb{R}^3}^{}  u \partial_x\partial_y u
    \dd{\mathbf{x}}}
    \\&\leq \norm{u}_{L^2} \norm{\partial_x\partial_y u}_{L^2} 
    =\norm{u}_{L^2} \norm{\xi_1 \xi_2 \hat{u}}_{L^2}  \\ 
    &\leq  \norm{u}_{L^2} \norm{\frac{\xi_1^2 + \xi_2^2 + \xi_3^2}{2}  \hat{u}}_{L^2} =
    \frac{1}{2} \norm{u}_{L^2} \norm{\Delta u}_{L^2}.
\end{aligned}
\end{equation}
Plugging the above estimate into \eqref{sec_1:eq_1}, we obtain 
\begin{equation}
    \begin{aligned}
	\norm{\mathscr{L}_{\star} u}_{L^2}^2 \geq {} &\norm{\nu \Delta u}_{L^2}^2 + \norm{y \partial_x u}_{L^2}^2 -\nu\norm{u}_{L^2} \norm{\Delta u}_{L^2} \\
     \geq {} &\norm{\nu \Delta u}_{L^2}^2 + \norm{y \partial_x u}_{L^2}^2 - \frac{1}{2} \qty(\norm{u}_{L^2}^2 + \nu^2\norm{\Delta u}_{L^2}^2) ,\\
\end{aligned}
\end{equation}
which further leads to 
\begin{equation}
    \begin{aligned}
	\norm{\mathscr{L}_{\star} u}_{L^2}^2 + \norm{u}_{L^2}^2\geq C \qty(\norm{\nu \Delta u}_{L^2}^2 + \norm{y \partial_x u}_{L^2}^2 + \norm{u}_{L^2}^2),\\
\end{aligned}
\end{equation}
for some constant \( C > 0\). And the inverse inequality 
\begin{equation}
    \begin{aligned}
	\norm{\mathscr{L}_{\star} u}_{L^2}^2 + \norm{u}_{L^2}^2 \leq  C \qty(\norm{\nu \Delta u}_{L^2}^2 + \norm{y \partial_x u}_{L^2}^2 + \norm{u}_{L^2}^2),\\
\end{aligned}
\end{equation}
is obvious by virtue of the definition of \( \mathscr{L}_{\star} \). That is, we found the equivalent norm for the graph norm of \( \mathscr{L}_{\star} \). Hence 
 \( \mathscr{L}_{\star} \) is closed with domain defined in \eqref{domain_of_L_star}. Since \( \mathscr{L}_{\bullet} \) is bounded on \( [L^2(\mathbb{R}^3)]^3 \), the operator \( \mathscr{L} = \mathscr{L}_{\star}I_3 + \mathscr{L}_{\bullet} \) is closed on \( \mathcal{D}(\mathscr{L})=\mathcal{D}(\mathscr{L}_{\star})^3 \) as in \eqref{domain_of_L}. 
\end{proof}

\begin{lemma}
	The operator \( \mathscr{L} \) generates a \( C_0 \)-semigroup on \( [L^2(\mathbb{R}^3)]^3 \).
\end{lemma}
\begin{proof}
	We prove the lemma following a perturbation argument. Recalling \eqref{decomposition_of_L}, we only need to show that \( \mathscr{L}_{\star } \) generates a \( C_0 \)-semigroup, since \( \mathscr{L}_{\bullet} \) is a bounded operator on \( [L^2(\mathbb{R}^3)]^3 \). 
  According to the Hille--Yosida Theorem \cite{Pazy1992}, the operator \( \mathscr{L}_{\star} \) generates a \( C_0 \)-semigroup on \( L^2(\mathbb{R}^3) \) if
	\begin{enumerate}[(i)]
		\item \( \mathscr{L}_{\star} \) is closed and \( \mathcal{D} (\mathscr{L}_{\star})\) is dense in \( L^2(\mathbb{R}^3) \);
		\item The resolvent set \( \rho(\mathscr{L}_{\star}) \) contains \( (0, + \infty) \), and for every \( \lambda > 0\) the following resolvent estimate holds
		      \begin{equation}\label{resolvent_estimate}
			      \norm{(\lambda I -\mathscr{L}_{\star})^{-1} } \leq  \frac{1}{\lambda }.
		      \end{equation}

	\end{enumerate}
	The domain is obviously dense in \( L^2(\mathbb{R}^3) \) since \( C_0^{ \infty}(\mathbb{R}^3) \subset \mathcal{D}(\mathscr{L}_{\star}) \). And the operator \( \mathscr{L}_{\star} \) is closed by \autoref{lemma:domain_of_L}. Therefore, \( C_0^{ \infty}(\mathbb{R}^3) \) is a core of \( \mathscr{L}_{\star}  \).

	Now, we derive the resolvent estimate \eqref{resolvent_estimate}. Consider for any \( F \in C_0^{ \infty}(\mathbb{R}^3) \subset  \mathcal{D}(\mathscr{L}_{\star}) \), \( \lambda > 0\),
	\begin{equation}
		\begin{aligned}
			\norm{(\lambda I-\mathscr{L}_{\star})F }_{L^2}^2  =  {} & \lambda^2 \norm{F}_{L^2}^2 + \norm{\mathscr{L}_{\star} F}_{L^2}^2 - 2\lambda \ip{ F}{\mathscr{L}_{\star} F}                \\
			\geq                                         & \lambda^2 \norm{F}_{L^2}^2 - 2\lambda \int_{\mathbb{R}^3}^{} ( \nu \Delta F- y \partial_x F)F  \dd{\mathbf{x}} \\
			\geq                                         & \lambda^2 \norm{F}_{L^2}^2 + 2\lambda \int_{\mathbb{R}^3}^{} \nu \abs{\nabla F}^2  \dd{\mathbf{x}}             \\
		\end{aligned}
	\end{equation}
	where the last inequality follows from integrating by parts. Hence, owing to the fact that \( C_0^{ \infty}(\mathbb{R}^3) \) is a core of \( \mathscr{L}_{\star} \) we obtain
	\begin{equation}\label{resolvent_estimate_1}
		\frac{1}{\lambda} \norm{(\lambda I -\mathscr{L}_{\star})F }_{L^2} \geq \norm{F}_{L^2}, \quad \forall  F \in \mathcal{D}(\mathscr{L}_{\star}).
	\end{equation}
	This shows that \( \lambda I -\mathscr{L}_{\star} \) is injective and \( (\lambda I-\mathscr{L}_{\star})^{-1} :\mathcal{R}(\lambda I -\mathscr{L}_{\star}) \to \mathcal{D}(\mathscr{L}_{\star}) \) is bounded by \( \frac{1}{\lambda} \). Also, \( \mathcal{R}(\lambda I -\mathscr{L}_{\star}) \) is closed in \( L^2(\mathbb{R}^3) \) \cite[Proposition 2.1 (iii)]{Schmuedgen2012}.

	Finally, we show that the range \( \mathcal{R}(\lambda -\mathscr{L}_{\star}) \) is also dense in \( L^2(\mathbb{R}^3) \).  Due to \cite[Corollary 2.2]{Schmuedgen2012}, we deduce that \( L^2(\mathbb{R}^3) = \mathcal{N}(\lambda -\mathscr{L}_{\star}^{*}) \oplus  \mathcal{R} (\lambda -\mathscr{L}_{\star})\), which suggests that we only have to verify \( \mathcal{N}(\lambda -\mathscr{L}_{\star}^{*}) = \{0\}\).  It is clear that, for \( F\in \mathcal{D}(\mathscr{L}_{\star}) \) and \( G\in C_0^{ \infty} (\mathbb{R}^3)\)
	\begin{equation}
		\ip{ \mathscr{L}_{\star} F}{G} = \ip{F}{( \nu \Delta + y \partial_x)G}.
	\end{equation}
	Hence the restriction of \( \mathscr{L}_{\star}^{*} \) on \( C_0^{ \infty}(\mathbb{R}^3) \) is \( \nu \Delta + y \partial_x \). Since \( \mathscr{L}_{\star}^{*} \) is closed and \( C_0^{ \infty}(\mathbb{R}^3) \subset \mathcal{D}( \mathscr{L}_{\star}^{*} ) \) is a core of \( \mathscr{L}_{\star}^{*} \), then \( \mathscr{L}_{\star}^{*} \) is the closure of \( \nu \Delta + y \partial_x \). Now, similar to the derivation of \eqref{resolvent_estimate_1}, one obtains that for \( \lambda > 0\) and \( G \in \mathcal{D}(\mathscr{L}_{\star}^{*}) \)
	\begin{equation}
		\frac{1}{\lambda} \norm{(\lambda I-\mathscr{L}_{\star}^{*})G }_{L^2} \geq \norm{G}_{L^2}.
	\end{equation}
	This shows \( \mathcal{N}(\lambda I-\mathscr{L}_{\star}^{*}) = \{0 \}\). The proof is complete.

\end{proof} 

\begin{lemma}\label{lemma_no_point_spectrum}
	The spectrum of the operator \( \mathscr{L} :\mathcal{D}(\mathscr{L}) \to \mathcal{H} \) consists entirely of continuous spectrum.
\end{lemma}
\begin{proof}
We prove the lemma by showing that neither the operator \( \mathscr{L} \) nor the adjoint \( \mathscr{L}^{*} \) has point spectrum. Indeed, if \( \lambda_r \) is in the residual spectrum of \( \mathscr{L} \), then \( \mathcal{R}(\lambda_r - \mathscr{L})^{\perp} \neq \{ \bm{0}\}  \). Since \( \mathcal{R}(\lambda_r - \mathscr{L})^{\perp} = \mathcal{N}(\overline{\lambda_r} - \mathscr{L}^{*}) \) (see \cite[Proposition 1.6 (ii), p. 9]{Schmuedgen2012}), then \( \overline{\lambda_r} \) is a point spectrum for \( \mathscr{L}^{*} \). So if both \( \mathscr{L} \) and \( \mathscr{L}^{*} \) have no point spectrum, we can conclude that the spectrum of \( \mathscr{L} \) consists entirely of continuous spectrum. 

    Suppose \( (\lambda, \mathbf{u}) \in \mathbb{C}\times \mathcal{D}(\mathscr{L})\) is an eigen-pair such that
	\begin{equation}\label{eigen_prob_0}
		\mathscr{L}  \mathbf{u} (\mathbf{x}) = \lambda  \mathbf{u} (\mathbf{x}).
	\end{equation}
    Then 
    \begin{equation}\label{eigen_prob_1-1}
        \mathscr{\hat{L}} \hat{ \mathbf{u} }(\bm{\xi}) = \lambda \hat{ \mathbf{u} }(\bm{\xi}),
    \end{equation}
	holds for a.e. \( \bm{\xi} \in \mathbb{R}^3\).  Here the operator \( \mathscr{\widehat{L}} \) is defined as
    \begin{align}
        \mathscr{\widehat{L}} \coloneqq
     \begin{pmatrix}
        - \nu \abs{\bm{\xi}}^2 + \xi_1 \partial_{\xi_2} + f \frac{\xi_1 \xi_2}{ \abs{\bm{\xi}}^2 }& f - 1 + (2 - f) \frac{\xi_1^2}{ \abs{\bm{\xi}}^2 } & 0\\ 
        - f + f \frac{\xi_2^2}{ \abs{\bm{\xi}}^2 }& - \nu \abs{\bm{\xi}}^2 + \xi_1 \partial_{\xi_2} + (2 - f) \frac{\xi_1 \xi_2}{ \abs{\bm{\xi}}^2 }& 0\\ 
        f \frac{\xi_2 \xi_3}{ \abs{\bm{\xi}}^2 }& (2 - f) \frac{\xi_1 \xi_3}{ \abs{\bm{\xi}}^2 }& - \nu \abs{\bm{\xi}}^2 + \xi_1 \partial_{\xi_2}
     \end{pmatrix} ,
    \end{align}
   with the domain 
    \begin{equation}
        \mathcal{D}(\mathscr{\hat{L}}) \coloneqq \qty{\hat{ \mathbf{u} } \in [L^2(\mathbb{R}^3)]^3\ \middle|\ \abs{\bm{\xi}}^2 \hat{ \mathbf{u} },\   \xi_1 \partial_{\xi_2}\hat{ \mathbf{u} } \in [L^2(\mathbb{R}^3)]^3 }.
    \end{equation}
    
Next, we show by contradiction that there exists no non-zero \(\mathbf{\widehat{u}}\) such that \eqref{eigen_prob_1-1} holds.
We rewrite the eigenvalue problem as the $\xi_1$ and $\xi_3$ parameterized ODE system of the independent variable $\xi_2$:
	\begin{equation}\label{linear_analysis_4}
		\xi_1 \partial_{\xi_2} \mathbf{\widehat{u}}
		=
		\begin{pmatrix}
        \lambda + \nu \abs{\bm{\xi}}^2 - f \frac{\xi_1 \xi_2}{ \abs{\bm{\xi}}^2 }& -f + 1 - (2- f) \frac{\xi_1^2}{ \abs{\bm{\xi}}^2 } & 0\\ 
         f - f \frac{\xi_2^2}{ \abs{\bm{\xi}}^2 }& \lambda + \nu \abs{\bm{\xi}}^2 - (2 - f) \frac{\xi_1 \xi_2}{ \abs{\bm{\xi}}^2 }& 0\\ 
        -f \frac{\xi_2 \xi_3}{ \abs{\bm{\xi}}^2 }& - (2 -f) \frac{\xi_1 \xi_3}{ \abs{\bm{\xi}}^2 }& \lambda + \nu \abs{\bm{\xi}}^2 
     \end{pmatrix} 
		\mathbf{\widehat{u}}.
	\end{equation} 

 Since $\mathbf{u}\not\equiv 0$, there exists $\bm{\xi}^\ast$ with $\xi_1^\ast \neq 0$ such that $\hat{ \mathbf{u} }(\bm{\xi}^\ast)\neq 0$. Without loss of generality, we assume $\xi_1^\ast>0$ throughout the proof. By the uniqueness of solutions to the ODE system \eqref{linear_analysis_4}, one has \( \hat{ \mathbf{u}  }(\xi_2;\xi_1^\ast, \xi_3^\ast) \neq \bm{0} \) for all \( \xi_2 \in \mathbb{R} \).
Since \(  \hat{ \mathbf{u} } \in \mathcal{D} (\mathscr{\widehat{L}}) \), then \(  \abs{\bm{\xi}}^2 \hat{ \mathbf{u} } \in  [L^2(\mathbb{R}^3)]^3\) and  \( \xi_1^\ast \partial_{\xi_2}\hat{ \mathbf{u} } \in [L^2(\mathbb{R}^3)]^3 \), hence \( \hat{ \mathbf{u}  }(\xi_2;\xi_1^\ast,\xi_3^\ast) \in[ C(\mathbb{R}) ]^3\) is bounded by Fubini's theorem and the Sobolev embedding. 

 On the other hand,  one obtains by multiplying \eqref{linear_analysis_4} by \( \overline{\hat{\mathbf{u}}}(\xi_2;\xi_1^\ast,\xi_3^\ast)  \) that
	\begin{equation}
		\begin{aligned}
            \frac{\xi_1^\ast}{2} \partial_{\xi_2} \abs{\hat{\mathbf{u}}} ^2 & =  (\Re \lambda +\nu \abs{\bm{\xi}}^2) \abs{\hat{\mathbf{u}}} ^2 + \Re \qty( \mathscr{M} \mathbf{\widehat{u}}   \cdot \mathbf{\overline{\widehat{u}}}),       
		\end{aligned}
	\end{equation}
    where $|\bm{\xi}|^2=|\xi_1^\ast|^2+|\xi_2|^2+|\xi_3^\ast|^2$, and the operator \( \mathscr{M} = \mathscr{M}(\xi_2;\xi_1^\ast,\xi_3^\ast) \) is identified as follows
    \begin{equation}
        \mathscr{M}(\xi_2;\xi_1^\ast,\xi_3^\ast) = \begin{pmatrix}
         - f \frac{\xi_1^\ast \xi_2}{ \abs{\bm{\xi}}^2 }& -f + 1 - (2 -f) \frac{{\xi_1^\ast}^2}{ \abs{\bm{\xi}}^2 } & 0\\ 
         f - f \frac{\xi_2^2}{ \abs{\bm{\xi}}^2 } & - (2 - f) \frac{\xi_1^\ast \xi_2}{ \abs{\bm{\xi}}^2 }& 0\\ 
        -f \frac{\xi_2 \xi_3^\ast}{ \abs{\bm{\xi}}^2 }& - (2 -f) \frac{\xi_1^\ast \xi_3^\ast}{ \abs{\bm{\xi}}^2 }& 0
     \end{pmatrix}.
    \end{equation}
    It is clear that  \( \abs{\mathscr{M} \hat{\mathbf{u}}}  \leq C(f) \abs{\hat{\mathbf{u}}}  \) for a constant \( C(f) > 0 \) depending on \( f \). It follows that
	\begin{equation}
		\begin{aligned}
			\frac{\xi_1^\ast}{2} \partial_{\xi_2}  \abs{\hat{\mathbf{u}}} ^2 & \geq   \qty(\Re \lambda +\nu \abs{\bm{\xi}}^2) \abs{\hat{\mathbf{u}}} ^2  - C(f) \abs{\hat{\mathbf{u}}}^2.
		\end{aligned}
	\end{equation}
	{
    Then, for large enough \( \abs{\xi_2} \) satisfying  
    \[
        \abs{\xi_2} \geq \xi_2^{*}(f,\nu,\lambda) \coloneqq \sqrt{\frac{1}{\nu} \qty(\frac{C(f)}{2} -\Re \lambda )} 
    \] 
    we have 
	\begin{equation}
		\begin{aligned}
			\xi_1^\ast\partial_{\xi_2}  \abs{\hat{\mathbf{u}}} ^2 & \geq (\Re \lambda +\nu \abs{\bm{\xi}}^2)  \abs{\hat{\mathbf{u}}} ^2.
		\end{aligned}
	\end{equation} 
One concludes that 
\begin{align*}
    |\widehat{\mathbf{u}}(  \xi_2;\xi_1^{*},\xi_3^{*})|^2 \geq |\widehat{\mathbf{u}}(\xi_2^{*};\xi_1^{*},\xi_3^{*})|^2 \cdot \exp\left( \frac{1}{\xi_1^{*}} \left[ \bigl(\Re\lambda + \nu \bigl((\xi_1^{*})^{2} + (\xi_3^{*})^{2}\bigr)\bigr)(\xi_2 - \xi_2^*) + \frac{\nu}{3} \left( \xi_2^3 - (\xi_2^*)^3 \right) \right] \right), 
\end{align*}
which is a contradiction. Hence $\hat{\mathbf{u}}\equiv 0$. Therefore, \( \mathscr{L} \) does not possess point spectrum.   
}

 The adjoint \( \mathscr{L}^{*} \) is given by
\begin{equation} 
   \begin{aligned}
	 \mathscr{L}^{*} &\coloneqq \mathscr{L}_0^{*}  + \mathscr{L}_1^{*} - \mathscr{L}_{2}, \\ \mathscr{L}_{2} & \coloneqq - \mathop{\mathrm{diag}} (y \partial_x, y \partial_x, y \partial_x) \\ 
    \mathscr{L}_0^{*} & \coloneqq 
    \begin{pmatrix}
        \nu \Delta &-f + f \Delta^{-1}  \partial_{y}^2& f \Delta^{-1}  \partial_{yz} \\ 
          f - 1 &  \nu \Delta & 0 \\ 
           0  & 0 &\nu \Delta  \\
    \end{pmatrix} , \\
    \mathscr{L}_1^{*} & \coloneqq  \Delta ^{-1}
    \begin{pmatrix}
       f  \partial_{xy}   & 0  & 0 \\
       \qty(2 - f )\partial_x^2  & \qty(2 - f )\partial_{xy} & \qty(2 - f )\partial_{xz} \\ 
       0 & 0  & 0 \\
    \end{pmatrix} ,  
\end{aligned}
\end{equation}
   with the domain \( \mathcal{D}(\mathscr{L}) = \mathcal{D}(\mathscr{L}^{*}) \). The same argument shows that \( \mathscr{L}^{*} \) also has no point spectrum. Consequently, the spectrum of both \( \mathscr{L} \) and \( \mathscr{L}^{*} \) consists entirely of continuous spectrum. This completes the proof.
\end{proof}

\begin{remark}\label{rem_affine_viscous}
The argument above uses in an essential way that the background flow is the affine Couette profile $U(y)=y$. After full Fourier transformation the transport term $y\partial_{x}$ becomes $\xi_{1}\partial_{\xi_{2}}$, and for fixed $(\xi_{1},\xi_{3})$ with $\xi_{1}\neq0$ the eigenvalue problem reduces to a first-order ODE system in $\xi_{2}$ alone, cf.~\eqref{linear_analysis_4}. For a general shear $U(y)$ the term $U(y)\partial_{x}$ becomes a convolution in $\xi_{2}$ after Fourier transformation, i.e. $\widehat{U\partial_{x}\mathbf{u}}=\xi_{1} \widehat{U}*\partial_{\xi_{2}}$ in the $\xi_{2}$-variable, and the ODE reduction is lost. Moreover viscosity is essential: the contradiction follows from the coercive term $\nu|\bm{\xi}|^{2}$, which forces super-exponential growth $\exp(\nu\xi_{2}^{3}/3\xi_{1}^{*})$ as $|\xi_{2}|\to\infty$ unless the slice vanishes identically; the estimate collapses when $\nu=0$.
\end{remark}

\subsection{Resolvent estimate}\label{sec:resolvent_simple} 

Throughout the section, we use $C$ to denote a generic positive constant depending only on $f$, whose value may change from line to line; for a real number $x$ we also write
\begin{equation}\label{eq:simp_pospart}
[x]_{+}\coloneqq\max\{x,0\}
\end{equation}
for its positive part. Every matrix norm in this section, written $\norm{\cdot}$, is the operator norm induced by the Euclidean norm on $\mathbb C^{n}$, and $\norm{\cdot}_{F}$ denotes the Frobenius norm.
  
Define the growth bound of an operator \( \mathscr{L} \) by
\begin{equation}\label{eq:simp_omega0}
\omega_{0}(\mathscr L)\coloneqq\lim_{t\to\infty}\frac1t\log\norm{e^{t\mathscr L}}_{L^{2}\to L^{2}}.
\end{equation}
In this subsection, we prove \autoref{thm:simp_growth}.

\medskip
\noindent\textbf{Setup.}
Under the Fourier transform, the decomposition \eqref{decomposition_of_L}
\begin{equation} 
    \mathscr{L} = \mathscr{L}_{\star }I_3 + \mathscr{L}_{\bullet } , 
\end{equation}
is given by 
\begin{equation}\label{eq:simp_A2}
\begin{aligned}
     \mathscr{\widehat{L}}_{\star } &\coloneqq -\nu|\bm{\xi}|^{2}+\xi_{1}\partial_{\xi_{2}} ,
     &
     \mathbf{s} &\coloneqq\frac{\bm{\xi}}{|\bm{\xi}|} ,\\
      \mathscr{\widehat{L}}_{\bullet } &\coloneqq
      \begin{pmatrix}
        A(\mathbf{s}) &   \bm{0}   \\ 
        \begin{matrix} fs_2s_3 & (2-f)s_1s_3 \end{matrix} & 0
      \end{pmatrix} ,
     &
      A(\mathbf{s}) & \coloneqq
      \begin{pmatrix}
        fs_1s_2 & f-1+(2-f)s_1^2\\
        -f+fs_2^2 & (2-f)s_1s_2
      \end{pmatrix} ,
\end{aligned}
\end{equation} 
where $A(\mathbf{s})$ denotes the top-left $2\times2$ block of $\mathscr{\widehat{L}}_{\bullet}$, acting on $(\widehat{u}_{1},\widehat{u}_{2})$.

Denote the horizontal quantities by 
\begin{equation}
   \mathbf{\widehat{u}} = (\mathbf{\widehat u}_{h},\widehat{u}_3 ),\quad
    \mathbf{\widehat{F}} = (\mathbf{\widehat F}_{h},\widehat{F}_3 ),\quad
   \mathbf{\widehat u}_{h}\coloneqq(\widehat u_{1},\widehat u_{2})^{\intercal} ,\quad 
\mathbf{\widehat F}_{h}\coloneqq(\widehat F_{1},\widehat F_{2})^{\intercal}
\end{equation}
The resolvent problem $(\lambda-\mathscr L)\mathbf{\widehat u}=\mathbf{\widehat F}$ decomposes into the coupled system
\begin{equation}\label{eq:simp_passive}
\begin{aligned}
(\lambda+\nu|\bm{\xi}|^{2})\mathbf{\widehat u}_{h}
-\xi_{1}\partial_{\xi_{2}}\mathbf{\widehat u}_{h} 
&=A(\mathbf{s}) \mathbf{\widehat u}_{h}+ \mathbf{\widehat F}_{h},
\\
(\lambda+\nu|\bm{\xi}|^{2})\widehat{u}_{3}-\xi_{1}\partial_{\xi_{2}}\widehat{u}_{3}
&=fs_{2}s_{3}\widehat{u}_{1}+(2-f)s_{1}s_{3}\widehat{u}_{2}+\widehat{F}_{3}.
\end{aligned}
\end{equation}
The block \( A(\mathbf{s}) \) in \eqref{eq:simp_A2} has trace $2s_{1}s_{2}$ and determinant
\begin{equation}\label{eq:simp_det}
\det A 
=fs_{1}^{2}-f(1-f)s_{3}^{2},
\end{equation}  
hence, with
\begin{equation}\label{eq:simp_Ddef}
D(\mathbf{s})\coloneqq f(1-f)s_{3}^{2}-s_{1}^{2}\bigl(f-s_{2}^{2}\bigr),
\end{equation}
the eigenvalues of $A$ are
\[
    \lambda_{\pm}=s_{1}s_{2}\pm\sqrt{D} \qand \operatorname{Re}\lambda_{+}=s_{1}s_{2}+\sqrt{\max\{D,0\}}.
\] 
Substituting $s_{2}^{2}=1-s_{1}^{2}-s_{3}^{2}$ further produces an alternative form of \eqref{eq:simp_Ddef}
\begin{equation}\label{eq:simp_Did}
D=\Lambda^{2}s_{3}^{2}+(1-f)s_{1}^{2}-s_{1}^{2}\bigl(s_{1}^{2}+s_{3}^{2}\bigr).
\end{equation}
Set the traceless matrix
\[
    N\coloneqq A-s_{1}s_{2}I_{2}.
\]
The matrix $N$ is uniformly bounded in the operator norm:
\begin{equation}\label{eq:simp_Nbdd}
\norm{N}\le\norm{N}_{F}\le C.
\end{equation} 
The following algebraic identity is used repeatedly below:
\begin{equation}\label{eq:simp_a12a21}
A_{12}A_{21}=D-N_{11}^{2},
\qquad
N_{11}=A_{11}-s_{1}s_{2}=-(1-f)s_{1}s_{2}.
\end{equation}  

In preparation for the proof of \autoref{thm:simp_growth}, we first give three lemmas.

\begin{lemma} \label{lem:simp_excess}
Set $\beta(\mathbf{s})\coloneqq[\operatorname{Re}\lambda_{+}(\mathbf{s})-\Lambda]_{+}$. Then, we have
\begin{equation}\label{eq:simp_excesspt}
\beta(\mathbf{s})\le C s_{1}^{2}
\qq{for all} \mathbf{s}\in\mathbb S^{2} ,
\end{equation}
and, along every characteristic $\bm{\xi}(\xi_{2})=(\xi_{1},\xi_{2},\xi_{3})$ with $\xi_{1}\neq0$
\begin{equation}\label{eq:simp_excessint}
\sup_{s<t}\ \int_{s}^{t}\beta\bigl(\mathbf{s}(\xi_{2}')\bigr) \dd{\xi_{2}'}
 \le C|\xi_{1}|.
\end{equation}
\end{lemma}
\begin{proof}
We split on the sign of $D$.

If $D<0$, then $\beta=[s_{1}s_{2}-\Lambda]_{+}$ gives
\[
s_{1}s_{2}-\Lambda\le|s_{1}|-\Lambda\le\frac{s_{1}^{2}}{2\Lambda}-\frac{\Lambda}{2}\le\frac{s_{1}^{2}}{2\Lambda},
\]
hence $\beta\le Cs_{1}^{2}$.

If $D\ge0$, introduce the polar variable
\[
v\coloneqq\sqrt{1-s_{2}^{2}}=\sqrt{s_{1}^{2}+s_{3}^{2}},
\]
so that the discriminant splits as
\[
D=\Lambda^{2}v^{2}-s_{1}^{2}\bigl(f(2-f)-s_{2}^{2}\bigr).
\]
For $v\ge v_{0}\coloneqq\frac{\Lambda}{2\bigl(1+\sqrt{1-f}+\Lambda\bigr)}$, the concavity inequality $\sqrt{X-Y}\le\sqrt X-\frac Y{2\sqrt X}$, $X\ge Y$, applied with $X=\Lambda^{2}v^{2}$ and $Y=s_{1}^{2}\bigl(f(2-f)-s_{2}^{2}\bigr)$, gives
\begin{align*}
\operatorname{Re}\lambda_{+}-\Lambda
&\le s_{1}s_{2}+\Lambda\bigl(\sqrt{1-s_{2}^{2}}-1\bigr)
+\frac{s_{1}^{2}\bigl|f(2-f)-s_{2}^{2}\bigr|}{2\Lambda v}\\
&\le \Bigl(\frac{s_{1}^{2}}{2\Lambda}+\frac{\Lambda}{2} s_{2}^{2}\Bigr)
+\Lambda\bigl(\sqrt{1-s_{2}^{2}}-1\bigr)
+Cs_{1}^{2}\\
&\le \Bigl(\frac{s_{1}^{2}}{2\Lambda}+\frac{\Lambda}{2} s_{2}^{2}\Bigr)
-\frac{\Lambda}{2} s_{2}^{2}
+Cs_{1}^{2} = \frac{s_{1}^{2}}{2\Lambda}+Cs_{1}^{2}.
\end{align*}
Hence $\beta\le Cs_{1}^{2}$.

For $v<v_{0}$, parametrise $s_{1}=tv$ with $|t|\le1$; then $s_{3}^{2}=v^{2}(1-t^{2})$, $|s_{2}|=\sqrt{1-v^{2}}$, and \eqref{eq:simp_Did} gives
\[
D=v^{2}\Bigl[\Lambda^{2}(1-t^{2})+t^{2}(1-f-v^{2})\Bigr]
\le v^{2}\bigl(\Lambda+\sqrt{1-f}\bigr)^{2},
\qquad
\sqrt{D}\le v\bigl(\Lambda+\sqrt{1-f}\bigr).
\]
$|s_{1}s_{2}|\le|s_{1}|\le v$. Collecting,
\[
\operatorname{Re}\lambda_{+}-\Lambda
\le|s_{1}s_{2}|+\sqrt{\max\{D,0\}}
\le v\bigl(1+\sqrt{1-f}+\Lambda\bigr)-\Lambda<0,
\]
for $v<v_{0}=\frac{\Lambda}{2(1+\sqrt{1-f}+\Lambda)}$. Hence $\beta$ vanishes identically when $v<v_{0}$. This proves \eqref{eq:simp_excesspt}.

Finally, we insert
\[
s_{1}^{2}=\frac{\xi_{1}^{2}}{|\bm{\xi}(\xi_{2}')|^{2}}
=\frac{\xi_{1}^{2}}{r^{2}+\xi_{2}'^{2}},
\qquad
\int_{-\infty}^{+\infty}\frac{\dd{\xi_{2}'}}{r^{2}+\xi_{2}'^{2}}=\frac{\pi}{r}\le\frac{\pi}{|\xi_{1}|}
\]
into the integral in \eqref{eq:simp_excessint} and conclude
\[
\sup_{s<t}\ \int_{s}^{t}\beta\bigl(\mathbf{s}(\xi_{2}')\bigr) \dd{\xi_{2}'}
\le C\int_{-\infty}^{+\infty}\frac{\xi_{1}^{2}}{r^{2}+\xi_{2}'^{2}} \dd{\xi_{2}'}
=\pi C\frac{\xi_{1}^{2}}{|\xi_{1}|}
=\pi C|\xi_{1}|.
\]
This proves \eqref{eq:simp_excessint}.
\end{proof}

\begin{lemma}\label{lem:simp_frame}
There exists $C>0$, depending only on $f$, such that, on every characteristic $\bm{\xi}(\xi_{2})=(\xi_{1},\xi_{2},\xi_{3})$ with $\xi_{1}\ne0$, there is a unitary matrix $\mathcal U(\xi_{2})\in\mathrm U(2)$, i.e.\ $\mathcal U(\xi_{2})\in\mathbb C^{2\times2}$ with $\mathcal U^{*}(\xi_{2})\mathcal U(\xi_{2})=I_{2}$, continuous on $\mathbb R$ and $C^{1}$ at every $\xi_{2}$ with $D(\mathbf{s}(\xi_{2}))\ne0$, for which:
\begin{enumerate}[label=\textup{(\roman*)}]
\item {\rm(Triangularization.)}
\begin{equation}\label{eq:simp_tri}
\mathcal U^{*}(\xi_{2}) A(\mathbf{s}(\xi_{2})) \mathcal U(\xi_{2})
=\begin{pmatrix}\lambda_{+}(\mathbf{s}(\xi_{2}))&\chi(\xi_{2})\\0&\lambda_{-}(\mathbf{s}(\xi_{2}))\end{pmatrix},
\qquad
|\chi(\xi_{2})|\le C.
\end{equation}
\item {\rm(Derivative estimate.)} We have the following estimate on $\mathbb R$,
\begin{equation}\label{eq:simp_adia}
\norm{\mathcal U^{*}\partial_{\xi_{2}}\mathcal U}
\le C\Bigl(\bigl|\partial_{\xi_{2}}\mathbf s\bigr|
+\bigl|\partial_{\xi_{2}}\sqrt{|D(\mathbf{s})|}\bigr|\Bigr).
\end{equation}
\item {\rm(Integrability.)}  
\begin{equation}\label{eq:simp_budget}
\int_{-\infty}^{+\infty}\norm{\mathcal U^{*}\partial_{\xi_{2}}\mathcal U} \dd{\xi_{2}}\le C .
\end{equation}
\end{enumerate}
\end{lemma}

\begin{proof}
\medskip\noindent\textbf{Step 1 (Set-up).}
Fix $(\xi_{1},\xi_{3})$ with $\xi_{1}\ne0$ and abbreviate
\begin{equation}\label{def_in_lem_simp_frame}
    \begin{aligned}
        r &\coloneqq\sqrt{\xi_{1}^{2}+\xi_{3}^{2}},
\qquad
\ell\coloneqq\sqrt{\Lambda^{2}\xi_{3}^{2}+(1-f) \xi_{1}^{2}},
\\
P_{0} &\coloneqq r^{2}\bigl(\Lambda^{2}\xi_{3}^{2}-f \xi_{1}^{2}\bigr),
\qquad
P_{1}\coloneqq\ell^{2} \xi_{2}^{2}+P_{0},
    \end{aligned}
\end{equation}
so that, by \eqref{eq:simp_Did},
\begin{equation}\label{eq:simp_Dline}
D(\mathbf{s}(\xi_{2}))=\frac{\ell^{2} \xi_{2}^{2}+P_{0}}{|\bm{\xi}|^{4}},
\qquad
\Lambda r\le\ell\le\sqrt{1-f} r,
\qquad
\ell^{2}-(1-f)^{2} \xi_{1}^{2}=f(1-f) r^{2},
\end{equation}
the identity holding because $\Lambda^{2}=f(1-f)$. The entries of $A(\mathbf s(\xi_{2}))$ relevant below are
\begin{equation}\label{eq:simp_entries}
A_{21}=-\frac{f r^{2}}{|\bm{\xi}|^{2}},
\qquad
N_{11}=-(1-f) \frac{\xi_{1}\xi_{2}}{|\bm{\xi}|^{2}}.
\end{equation}

\medskip\noindent\textbf{Step 2 (Construction of the unitary frame).}
We set $\tilde{\mathbf{q}}_{+}\coloneqq(\sqrt{D}+N_{11},\ A_{21})^{\intercal}$; using \eqref{eq:simp_a12a21},
\[
\begin{aligned}
(N-\sqrt{D}) \tilde{\mathbf{q}}_{+}
&=\begin{pmatrix}
(N_{11}-\sqrt{D})(\sqrt{D}+N_{11})+A_{12}A_{21}\\[2pt]
A_{21} (\sqrt{D}+N_{11})-(N_{11}+\sqrt{D}) A_{21}
\end{pmatrix}\\
&=\begin{pmatrix}
N_{11}^{2}-D+A_{12}A_{21}\\
0
\end{pmatrix}
=\begin{pmatrix}0\\0\end{pmatrix} .
\end{aligned}
\]
The second row vanishes identically because $N_{22}=-N_{11}$, and the first row vanishes by \eqref{eq:simp_a12a21}. 

Hence, $N\tilde{\mathbf{q}}_{+}=\sqrt{D} \tilde{\mathbf{q}}_{+}$. One can also verify vector $\tilde{\mathbf{q}}_{-}\coloneqq(-\sqrt{D}+N_{11},A_{21})^{\intercal}$ satisfies $(N+\sqrt{D})\tilde{\mathbf{q}}_{-}=0$. And 
\[
|\tilde{\mathbf{q}}_{+}| \ge |A_{21}|=\frac{f r^{2}}{|\bm{\xi}|^{2}}>0 .
\]
The term $A_{21}$ vanishes only as $|\xi_{2}|\to\infty$. Normalize the eigenvectors, one has
\begin{equation}\label{def_of_q}
    \mathbf q_{1}\coloneqq\frac{\tilde{\mathbf{q}}_{+}}{|\tilde{\mathbf{q}}_{+}|},
\qquad
\mathbf q_{2}\coloneqq\bigl(-\overline{(\mathbf q_{1})_{2}},\ \overline{(\mathbf q_{1})_{1}}\bigr)^{\intercal},
\qquad
\mathcal U\coloneqq(\mathbf q_{1},\mathbf q_{2}).
\end{equation}
Hence, by definition, $\mathcal U$ is unitary, continuous on $\mathbb R$, and $C^{1}$ wherever $D\ne0$. 

Since $N\mathbf q_{1}=\sqrt{D} \mathbf q_{1}$ and $N$ being traceless,
\[
\mathcal U^{*}A\mathcal U
=s_{1}s_{2} I_{2}+\mathcal U^{*}N\mathcal U
=\begin{pmatrix}s_{1}s_{2}+\sqrt{D}&\chi\\0&s_{1}s_{2}-\sqrt{D}\end{pmatrix}
=\begin{pmatrix}\lambda_{+}&\chi\\0&\lambda_{-}\end{pmatrix}
\]
with $|\chi|\le\norm{N}$, and $\norm{N}\le C$ by \eqref{eq:simp_Nbdd}. This concludes \eqref{eq:simp_tri}.

\medskip\noindent\textbf{Step 3 (Estimates of the derivative).}
Write $\dot{\phantom{a}}\coloneqq\frac d{d\xi_{2}}$ and $E\coloneqq\mathcal U^{*}\dot{\mathcal U}$. The quotient rule applied to $\mathbf q_{1}=\tilde{\mathbf q}_{+}/|\tilde{\mathbf q}_{+}|$ gives
\[
\begin{aligned}
\dot{\mathbf q}_{1}
&=\frac{\dot{\tilde{\mathbf q}}_{+}}{|\tilde{\mathbf q}_{+}|}
-\frac{\tilde{\mathbf q}_{+}}{|\tilde{\mathbf q}_{+}|^{3}} 
\operatorname{Re}\bigl(\tilde{\mathbf q}_{+}^{*}\dot{\tilde{\mathbf q}}_{+}\bigr)\\
|\dot{\mathbf q}_{1}|
&\le\frac{|\dot{\tilde{\mathbf q}}_{+}|}{|\tilde{\mathbf q}_{+}|}
+\frac{|\tilde{\mathbf q}_{+}| \bigl|\tilde{\mathbf q}_{+}^{*}\dot{\tilde{\mathbf q}}_{+}\bigr|}{|\tilde{\mathbf q}_{+}|^{3}}
\le\frac{2 |\dot{\tilde{\mathbf q}}_{+}|}{|\tilde{\mathbf q}_{+}|}.
\end{aligned}
\]  
The flip-conjugation defining $\mathbf q_{2}$ commutes with differentiation and preserves norms, whence
\[
\dot{\mathbf q}_{2}
=\bigl(-\overline{(\dot{\mathbf q}_{1})_{2}},\ \overline{(\dot{\mathbf q}_{1})_{1}}\bigr)^{\intercal},
\qquad
|\dot{\mathbf q}_{2}|=|\dot{\mathbf q}_{1}|.
\]
Finally, differentiating the components of $\tilde{\mathbf q}_{+}$ and using the triangle inequality,
\[
|\dot{\tilde{\mathbf q}}_{+}|^{2}
=|\partial_{\xi_{2}}\sqrt{D}+\dot N_{11}|^{2}+|\dot A_{21}|^{2}
\le\bigl(|\partial_{\xi_{2}}\sqrt{D}|+|\dot N_{11}|+|\dot A_{21}|\bigr)^{2},
\]
with
\[
\dot N_{11}=-(1-f) \bigl(\dot s_{1}s_{2}+s_{1}\dot s_{2}\bigr),
\qquad
\dot A_{21}=-2f \bigl(s_{1}\dot s_{1}+s_{3}\dot s_{3}\bigr),
\]
so that, one has $|\dot N_{11}|\le2(1-f) |\dot{\mathbf s}|$ and $|\dot A_{21}|\le2f |\dot{\mathbf s}|$; therefore
\[
|\dot{\tilde{\mathbf q}}_{+}|
\le|\partial_{\xi_{2}}\sqrt{D}|+|\dot N_{11}|+|\dot A_{21}|
\le|\partial_{\xi_{2}}\sqrt{D}|+2 |\dot{\mathbf s}|.
\]

By the quotient rule, each of the two normalized derivatives splits into a tangential part and a radial part,
\[
\dot{\mathbf q}_{1}
=\frac{\dot{\tilde{\mathbf q}}_{+}}{|\tilde{\mathbf q}_{+}|}
-\mathbf q_{1} \frac{\operatorname{Re}\bigl(\tilde{\mathbf q}_{+}^{*}\dot{\tilde{\mathbf q}}_{+}\bigr)}{|\tilde{\mathbf q}_{+}|^{2}},
\qquad
\dot{\mathbf q}_{2}
=\frac{\dot{\tilde{\mathbf q}}_{+}^{\perp}}{|\tilde{\mathbf q}_{+}|}
-\mathbf q_{2} \frac{\operatorname{Re}\bigl((\tilde{\mathbf q}_{+}^{\perp})^{*}\dot{\tilde{\mathbf q}}_{+}^{\perp}\bigr)}{|\tilde{\mathbf q}_{+}|^{2}},
\]
the radial parts being parallel to $\mathbf q_{1}$ and $\mathbf q_{2}$ respectively; by Cauchy--Schwarz, the radial part of $\dot{\mathbf q}_{1}$ has modulus at most $|\dot{\tilde{\mathbf q}}_{+}|/|\tilde{\mathbf q}_{+}|$, so $|\dot{\mathbf q}_{1}|\le2 |\dot{\tilde{\mathbf q}}_{+}|/|\tilde{\mathbf q}_{+}|$.

Next, we investigate the estimates for \( E \). Since the entries of $E=\mathcal U^{*}\dot{\mathcal U}$ are the matrix products $E_{jk}=\mathbf q_{j}^{*}\cdot\dot{\mathbf q}_{k}$; Cauchy--Schwarz and $|\mathbf q_{j}|=1$ give $|E_{jk}|\le\bigl|\mathbf q_{j}^{*}\bigr| \bigl|\dot{\mathbf q}_{k}\bigr|=\bigl|\dot{\mathbf q}_{k}\bigr|$ for all $j,k$; since $\norm{E}\le\norm{E}_{F}$ and
$\norm{E}_{F}^{2}=\sum_{j,k}\bigl|E_{jk}\bigr|^{2}
\le2 |\dot{\mathbf q}_{1}|^{2}+2 |\dot{\mathbf q}_{2}|^{2}
=4 |\dot{\mathbf q}_{1}|^{2}$;
and $|\partial_{\xi_{2}}\sqrt{D}|=|\partial_{\xi_{2}}\sqrt{|D(\mathbf{s})|}|$ whether $D$ is positive or negative,
\begin{equation}\label{eq:simp_Ecrude}
\norm{E} \le 2 |\dot{\mathbf q}_{1}|
 \le \frac{4 \bigl(|\partial_{\xi_{2}}\sqrt{D}|+2 |\dot{\mathbf s}|\bigr)}{|\tilde{\mathbf{q}}_{+}|}.
\end{equation}

Hereafter, we offer another estimate for \( E \). Wherever $D>0$, the vector $\tilde{\mathbf{q}}_{+}$ is real. Furthermore, by \eqref{def_of_q}, the rotated vector $\tilde{\mathbf q}_{+}^{\perp}\coloneqq(-A_{21},\ \sqrt{D}+N_{11})^{\intercal}$ gives the column vector $\mathbf q_{2}=\tilde{\mathbf q}_{+}^{\perp}/|\tilde{\mathbf q}_{+}|$; the matrix $E$ is thus real antisymmetric, so $\norm{E}=|E_{12}|$, and  
\begin{equation}\label{eq:simp_Ecross}
\norm{E}=\frac{\bigl|(\sqrt{D}+N_{11}) \dot A_{21}-A_{21} (\partial_{\xi_{2}}\sqrt{D}+\dot N_{11})\bigr|}{|\tilde{\mathbf{q}}_{+}|^{2}} .
\end{equation}
Indeed, the entries of $E$ are $E_{jk}=\mathbf q_{j}^{*}\cdot\dot{\mathbf q}_{k}$, and explicit calculations leads to
\[
E_{11}=\mathbf q_{1}^{*}\cdot\dot{\mathbf q}_{1}
=\tfrac12 \tfrac d{d\xi_{2}} |\mathbf q_{1}|^{2}=0,
\qquad
E_{22}=0,
\qquad
E_{21}=-E_{12};
\]
moreover,
\[
\begin{aligned}
E_{12}
=\mathbf q_{1}^{*}\cdot\dot{\mathbf q}_{2}
&=\frac{\tilde{\mathbf q}_{+}^{*}\cdot\dot{\tilde{\mathbf q}}_{+}^{\perp}}{|\tilde{\mathbf q}_{+}|^{2}}
=\frac{A_{21} (\partial_{\xi_{2}}\sqrt{D}+\dot N_{11})-(\sqrt{D}+N_{11}) \dot A_{21}}{|\tilde{\mathbf q}_{+}|^{2}}.
\end{aligned}
\]
The two estimates divide the work in Step 4: \eqref{eq:simp_Ecrude} is used in the interior, where $|\tilde{\mathbf q}_{+}|$ has a uniform lower bound, while \eqref{eq:simp_Ecross} takes over in the tail.

\medskip\noindent\textbf{Step 4 (Interior and tail estimates).}
The proof distinguishes the two regions $\{|\xi_{2}|\le\vartheta r\}$ and $\{|\xi_{2}|>\vartheta r\}$, where $\vartheta\coloneqq\sqrt{2/(1-f)}$.

The choice of $\vartheta$ ensures that the term of the numerator in \eqref{eq:simp_Dline} shares the following relation when \( |\xi_{2}|>\vartheta r \)
\[
\begin{aligned}
\frac{|P_{0}|}{\ell^{2} \xi_{2}^{2}}
&=\frac{r^{2} \bigl|\Lambda^{2}\xi_{3}^{2}-f \xi_{1}^{2}\bigr|}{\ell^{2} \xi_{2}^{2}}
 \le \frac{f r^{4}}{f(1-f) r^{2} \xi_{2}^{2}}
\\
&=\frac{1}{1-f}\cdot\frac{r^{2}}{\xi_{2}^{2}}
\ <\ \frac{1}{(1-f) \vartheta^{2}}
 = \frac12.
\end{aligned}
\] 
So that \( D > 0 \), and thus \( E \) can be estimated with the finer bound \eqref{eq:simp_Ecross}.

For the estimate in the interior $|\xi_{2}|\le\vartheta r$, we have $|\bm{\xi}|^{2}\le(1+\vartheta^{2}) r^{2}$, so
\[
|\tilde{\mathbf{q}}_{+}| \ge |A_{21}|=\frac{f r^{2}}{|\bm{\xi}|^{2}} \ge \frac{f(1-f)}{3-f}\ >\ 0 ,
\]
and the crude bound \eqref{eq:simp_Ecrude} yields \eqref{eq:simp_adia} on this region.

Next, we estimate the tail $|\xi_{2}|>\vartheta r$.
Here the second component $A_{21}$ decays, so the crude bound \eqref{eq:simp_Ecrude} is useless, and the refinement \eqref{eq:simp_Ecross} applies. The argument rests on three uniform facts.
\begin{enumerate}[label=\textup{(\alph*)}]
\item  First, we estimate the size of $P_{1}$. Recalling \eqref{def_in_lem_simp_frame}, we have
\begin{equation}\label{estimate_lem_1}
    \tfrac12 \ell^{2} \xi_{2}^{2} \le P_{1} \le \bigl(\ell |\xi_{2}|+\sqrt{f} r^{2}\bigr)^{2} \le 2 r^{2} \xi_{2}^{2},
\qquad
|\bm{\xi}|^{2}\le\tfrac32 \xi_{2}^{2},
\end{equation}
using $P_{1}\ge\ell^{2}\xi_{2}^{2}-|P_{0}|$, $|P_{0}|\le f r^{4}$, $r^{2}\le r |\xi_{2}|$ and $\bigl(\sqrt{1-f}+\sqrt f\bigr)^{2}\le2$.
\item Second, we consider the lower bound for \( |\tilde{\mathbf{q}}_{+ }| \). By \eqref{eq:simp_entries} and $\sqrt{D}=\sqrt{P_{1}}/|\bm{\xi}|^{2}$ as above, the first component of $\tilde{\mathbf q}_{+}$ is
$\sqrt{D}+N_{11}=\bigl(\sqrt{P_{1}}-(1-f)\xi_{1}\xi_{2}\bigr)/|\bm{\xi}|^{2}$
and has a lower bound as follows: by \eqref{eq:simp_Dline} 
\[
P_{1}-(1-f)^{2} \xi_{1}^{2} \xi_{2}^{2}
=\bigl(\ell^{2}-(1-f)^{2}\xi_{1}^{2}\bigr) \xi_{2}^{2}+P_{0}
=f(1-f) r^{2} \xi_{2}^{2}+P_{0}
 \ge \tfrac12 f(1-f) r^{2} \xi_{2}^{2},
\]
where we used $|P_{0}|\le f r^{4}$ and $|\xi_{2}|>\vartheta r$. Then, using the crude upper bound 
\[
    \bigl|\sqrt{P_{1}}+(1-f)\xi_{1}\xi_{2}\bigr|\le\sqrt{P_{1}}+(1-f)|\xi_{1}||\xi_{2}|\le(\sqrt2+1-f) r |\xi_{2}|
\]  
on the denominator,
\[
\begin{aligned}
\bigl|\sqrt{D}+N_{11}\bigr|
&=\frac{\bigl|P_{1}-(1-f)^{2}\xi_{1}^{2}\xi_{2}^{2}\bigr|}{|\bm{\xi}|^{2} \bigl|\sqrt{P_{1}}+(1-f)\xi_{1}\xi_{2}\bigr|}
\\
&\ge \frac{\frac12 f(1-f) r^{2} \xi_{2}^{2}}{\frac32 \xi_{2}^{2}\cdot(\sqrt2+1-f) r |\xi_{2}|}
 \ge C \frac{r}{|\xi_{2}|},
\end{aligned}
\]
hence $|\tilde{\mathbf{q}}_{+}|\ge\bigl|\sqrt{D}+N_{11}\bigr|\ge C  r/|\xi_{2}|$.
\item Third, we give the bound of terms in \( E \). The coefficients entering \eqref{eq:simp_Ecross} obey
\begin{equation}
    \begin{aligned}
        \bigl|\sqrt{D}+N_{11}\bigr|
&\le|\sqrt{D}|+|N_{11}| 
\le\frac{\sqrt{P_{1}}}{|\bm{\xi}|^{2}}+\frac{(1-f) r}{|\xi_{2}|}
\le C \frac{r}{|\xi_{2}|},
\\
|A_{21}|&\le\frac{f r^{2}}{\xi_{2}^{2}},
\qquad
|\dot A_{21}|\le\frac{2f r^{2}}{|\xi_{2}|^{3}},
    \end{aligned}
\end{equation}
where we have used \eqref{eq:simp_Dline}, \eqref{eq:simp_entries}, \eqref{estimate_lem_1} with $|\xi_{1}|\le r$ and $|\bm{\xi}|^{2}\ge\xi_{2}^{2}$. 

And for the term \( |\partial_{\xi_{2}}\sqrt{D}| \) and \( \dot N_{11} \), we have
\[
|\partial_{\xi_{2}}\sqrt{D}|=\Bigl|\frac{\ell^{2} \xi_{2}}{|\bm{\xi}|^{2}\sqrt{P_{1}}}-\frac{2 \xi_{2} \sqrt{P_{1}}}{|\bm{\xi}|^{4}}\Bigr|
 \le \frac{\sqrt2 \ell}{\xi_{2}^{2}}+\frac{2 \sqrt{P_{1}}}{|\xi_{2}|^{3}}
 \le C \frac{r}{\xi_{2}^{2}},
\]
and
\[
|\dot N_{11}|=\Bigl|(1-f) \xi_{1} \frac{r^{2}-\xi_{2}^{2}}{|\bm{\xi}|^{4}}\Bigr|
\le\frac{(1-f) r}{\xi_{2}^{2}} ,
\]
the estimates using \eqref{estimate_lem_1}, $\ell\le\sqrt{1-f} r$ from \eqref{eq:simp_Dline}, $|\bm{\xi}|^{4}\ge\xi_{2}^{4}$, and $\xi_{2}^{2}>r^{2}$.
\end{enumerate}
Combining the three uniform estimate on the tail, we then have
\[
\begin{aligned}
&\bigl|(\sqrt{D}+N_{11}) \dot A_{21}-A_{21} (\partial_{\xi_{2}}\sqrt{D}+\dot N_{11})\bigr| \\
&\le \bigl|\sqrt{D}+N_{11}\bigr| \bigl|\dot A_{21}\bigr|
+|A_{21}| \bigl(|\partial_{\xi_{2}}\sqrt{D}|+|\dot N_{11}|\bigr) \le C \frac{r^{3}}{|\xi_{2}|^{4}} .
\end{aligned}
\] 
Dividing by $|\tilde{\mathbf{q}}_{+}|^{2}\ge C  r^{2}/\xi_{2}^{2}$ from (b), equation \eqref{eq:simp_Ecross} gives
\[
\norm{E} 
 \le C \frac{r}{\xi_{2}^{2}}
 = C \frac{|\bm{\xi}|^{2}}{\xi_{2}^{2}}\cdot\frac{r}{|\bm{\xi}|^{2}}
 \le \tfrac32 C \bigl|\partial_{\xi_{2}}\mathbf s\bigr| ,
\]
because $|\partial_{\xi_{2}}\mathbf s|=r/|\bm{\xi}|^{2}$ and $|\bm{\xi}|^{2}\le\frac32 \xi_{2}^{2}$ from \eqref{estimate_lem_1}.
Collecting the two regions proves \eqref{eq:simp_adia} on the whole line \( \mathbb{R} \).

\medskip\noindent\textbf{Step 5 (Integrability).}
Now, we deal with the integral \eqref{eq:simp_budget}. First,
\[
\int_{\mathbb R}|\partial_{\xi_{2}}\mathbf s| \dd{\xi_{2}}
=\int_{\mathbb R}\frac{r \dd{\xi_{2}}}{r^{2}+\xi_{2}^{2}}
=\pi ,
\]
an arctangent integral.

Second, $\sqrt{|D(\mathbf{s}(\cdot))|}$ has total variation at most $C$. Indeed, finiteness of the zeros of the derivative suffices: the function is then monotone on each of the finitely many intervals between consecutive zeros, and its total variation is the sum of the corresponding endpoint differences. 

Here, by \eqref{eq:simp_Dline}, $D(\mathbf{s}(\xi_{2}))=P_{1}/|\bm{\xi}|^{4}$ is quadratic over quartic with respect to \( \xi_2 \), so the numerator of $D'(\mathbf s(\xi_{2}))$, after clearing denominators, is a polynomial of degree at most five; adjoining the at most two zeros of $P_{1}$ gives at most seven zeros in total, at most eight monotone intervals, and hence the total variation
\[
\operatorname{TV}_{\mathbb R}\bigl(\sqrt{|D(\mathbf{s}(\cdot))|}\bigr)
 \le 8\cdot2\max_{\mathbb S^{2}}\sqrt{|D|}
 = 16\max_{\mathbb S^{2}}\sqrt{|D|}
 \le C.
\] 
Thus, integrating \eqref{eq:simp_adia} gives the explicit estimate
\[
\begin{aligned}
\int_{\mathbb R}\norm{\mathcal U^{*}\partial_{\xi_{2}}\mathcal U} \dd{\xi_{2}}
&\le C\int_{\mathbb R}\Bigl(|\partial_{\xi_{2}}\mathbf s|
+\bigl|\partial_{\xi_{2}}\sqrt{|D(\mathbf{s})|}\bigr|\Bigr) \dd{\xi_{2}}
\\
& 
\le C\Bigl(\pi+\operatorname{TV}_{\mathbb R}\bigl(\sqrt{|D(\mathbf{s}(\cdot))|}\bigr)\Bigr) \le C .
\end{aligned}
\]
 This proves \eqref{eq:simp_budget} and finishes the proof.
\end{proof} 

\begin{lemma}\label{lem:simp_backward}
Let $\lambda_{+},\lambda_{-},\chi$ be continuous complex-valued functions on $[s,t]$ with
$|\chi(\xi_{2}')|\le C$ and
$\operatorname{Re}\lambda_{-}(\xi_{2}')\le\operatorname{Re}\lambda_{+}(\xi_{2}')$ for all $\xi_{2}'\in[s,t]$, and let
\[
\mathcal{T}(\xi_{2}')
\coloneqq
\begin{pmatrix}\lambda_{+}(\xi_{2}')&\chi(\xi_{2}')\\0&\lambda_{-}(\xi_{2}')\end{pmatrix}
\qfor \xi_{2}'\in[s,t].
\]
Let $\Phi^{0}=\Phi^{0}(\xi_{2},\xi_{2}')\colon[s,t]^{2}\to\mathbb C^{2\times2}$ denote the transition matrix from $\xi_{2}'$ to $\xi_{2}$, that is, $\mathbf{\widehat u}_{h}(\xi_{2})=\Phi^{0}(\xi_{2},\xi_{2}') \mathbf{\widehat u}_{h}(\xi_{2}')$ for every solution $\mathbf{\widehat u}_{h}$ of the homogeneous system; equivalently, $\Phi^{0}$ solves
\[
\partial_{\xi_{2}}\Phi^{0}(\xi_{2},\xi_{2}')
=\frac1{\xi_{1}}\bigl(\lambda+\nu|\bm{\xi}(\xi_{2})|^{2}I_{2}-\mathcal{T}(\xi_{2})\bigr)\Phi^{0}(\xi_{2},\xi_{2}'),
\qquad
\Phi^{0}(\xi_{2},\xi_{2})=I_{2}\qfor \xi_{2}\in[s,t],
\]
the derivative acting on the first variable, the target of the propagation; here $\lambda\in\mathbb C$, $\xi_{1}>0$ are parameters and $\bm{\xi}(\xi_{2}')=(\xi_{1},\xi_{2}',\xi_{3})$. Then for $s<t$,
\begin{equation}\label{eq:simp_phiseg}
\norm{\Phi^{0}(s,t)}
\le C\Bigl(1+\frac{t-s}{\xi_{1}}\Bigr)
\exp\Bigl(-\frac1{\xi_{1}}\int_{s}^{t}
\operatorname{Re}\bigl[\lambda+\nu|\bm{\xi}|^{2}-\lambda_{+}\bigr] \dd{\xi_{2}'}\Bigr).
\end{equation}
\end{lemma}
\begin{proof}
Since $\mathcal T$ is upper triangular, so is every transition matrix of the system; write
\[
\Phi^{0}(s,t)=
\begin{pmatrix}\Upsilon(s,t)&\Psi(s,t)\\0&\Xi(s,t)\end{pmatrix}.
\]
The diagonal entries are pure quadratures,
\[
\Upsilon(s,t)=\exp\Bigl(-\frac1{\xi_{1}}\int_{s}^{t}
\bigl(\lambda+\nu|\bm{\xi}(\xi_{2}')|^{2}-\lambda_{+}(\xi_{2}')\bigr) \dd{\xi_{2}'}\Bigr),
\]
\[
\Xi(s,t)=\exp\Bigl(-\frac1{\xi_{1}}\int_{s}^{t}
\bigl(\lambda+\nu|\bm{\xi}(\xi_{2}')|^{2}-\lambda_{-}(\xi_{2}')\bigr) \dd{\xi_{2}'}\Bigr),
\]
whence, using the hypothesis $\operatorname{Re}\lambda_{-}\le\operatorname{Re}\lambda_{+}$ pointwise,
\begin{equation}\label{estimate_lem_2}
    |\Xi(s,t)|\le|\Upsilon(s,t)|
=\exp\Bigl(-\frac1{\xi_{1}}\int_{s}^{t}
\operatorname{Re}\bigl[\lambda+\nu|\bm{\xi}(\xi_{2}')|^{2}-\lambda_{+}(\xi_{2}')\bigr] \dd{\xi_{2}'}\Bigr).
\end{equation}
Variation of constants on the first component gives that
\[
\Psi(s,t)
=-\frac{\Upsilon(s,t)}{\xi_{1}}\int_{s}^{t}\chi(\sigma) 
\frac{\Xi(\sigma,t)}{\Upsilon(\sigma,t)} \dd{\sigma} ,
\]
and, taking moduli in
$\Xi(\sigma,t)/\Upsilon(\sigma,t)=\exp\bigl(\tfrac1{\xi_{1}}\int_{\sigma}^{t}\bigl(\lambda_{+}(\xi_{2}')-\lambda_{-}(\xi_{2}')\bigr) \dd{\xi_{2}'}\bigr)$,
\[
\begin{aligned}
|\Psi(s,t)|
&\le\frac{|\Upsilon(s,t)|}{\xi_{1}}\int_{s}^{t}|\chi(\sigma)| 
\exp\Bigl(\frac1{\xi_{1}}\int_{\sigma}^{t}
\operatorname{Re}\bigl(\lambda_{+}(\xi_{2}')-\lambda_{-}(\xi_{2}')\bigr) \dd{\xi_{2}'}\Bigr) \dd{\sigma}\\
&\le\frac{|\Upsilon(s,t)| |\chi| (t-s)}{\xi_{1}} 
\exp\Bigl(\frac1{\xi_{1}}\int_{s}^{t}
\operatorname{Re}\bigl(\lambda_{+}(\xi_{2}')-\lambda_{-}(\xi_{2}')\bigr) \dd{\xi_{2}'}\Bigr).
\end{aligned}
\]
The exponential admits the following estimate by \eqref{estimate_lem_2}
\[
|\Upsilon(s,t)| 
\exp\Bigl(\frac1{\xi_{1}}\int_{s}^{t}
\operatorname{Re}\bigl(\lambda_{+}(\xi_{2}')-\lambda_{-}(\xi_{2}')\bigr) \dd{\xi_{2}'}\Bigr)
=|\Xi(s,t)| \le |\Upsilon(s,t)|,
\]
whence
\[
|\Psi(s,t)| \le \frac{|\chi| (t-s)}{\xi_{1}} |\Upsilon(s,t)|.
\]
The above estimate, together with \eqref{estimate_lem_2} and $|\chi|\le C$, gives \eqref{eq:simp_phiseg}.
\end{proof}

\begin{proof}[Proof of \autoref{thm:simp_growth}]
\medskip\noindent\textbf{Step 1 (lower bound).}
By \autoref{theorem_spectral_instability}, one has $\omega_{0}(\mathscr L)\ge\alpha(\mathscr L)$. The lower bound $\alpha(\mathscr L)\ge\Lambda$ is deferred to \autoref{coro_spec_lower_bound} in \autoref{sec_linear_instability}, which establishes the spectral inclusion $(-\infty,\Lambda]\subset\sigma(\mathscr L)$.

\medskip\noindent\textbf{Step 2 (Gearhart--Pr\"uss reduction).}
Since $\mathscr L$ generates a $C_{0}$-semigroup on $ [L^{2}(\mathbb R^{3})]^{3}$, for every $\omega\in\mathbb R$ the rescaled operator $\mathscr L-\omega$ generates $e^{t(\mathscr L-\omega)}=e^{-\omega t}e^{t\mathscr L}$. Choosing $\omega$ large enough that this semigroup is bounded, \autoref{thm:gearhart_pruss} applied to $\mathscr L-\omega$ yields
\begin{equation}\label{eq:simp_gp}
\omega_{0}(\mathscr L)
=\inf\Bigl\{\omega:\sup_{\Re\lambda>\omega}\norm{(\lambda-\mathscr L)^{-1}}<\infty\Bigr\},
\end{equation}
so that the growth-bound assertion of the theorem follows from the resolvent estimate \eqref{eq:simp_est}. The rest of the proof is therefore devoted to \eqref{eq:simp_est}. 

Similar to the proof of \autoref{lemma_no_point_spectrum}, the resolvent estimate can be considered slice-wise: on a slice with \( \xi_1\neq 0\) the resolvent problem \eqref{eq:simp_passive} has the form
\begin{equation}\label{eq:simp_ode}
\begin{aligned}
\partial_{\xi_{2}}\mathbf{\widehat{u}}_{h}
&=\frac1{\xi_{1}}\Bigl(\lambda+\nu|\bm{\xi}|^{2}-A(\mathbf{s}(\xi_{2}))\Bigr)\mathbf{\widehat{u}}_{h}
-\frac1{\xi_{1}}\mathbf{\widehat{F}}_{h} ,
\\
\partial_{\xi_{2}}\widehat{u}_{3}
&=\frac1{\xi_{1}}\bigl(\lambda+\nu|\bm{\xi}|^{2}\bigr)\widehat{u}_{3}
-\frac1{\xi_{1}}\bigl(fs_{2}s_{3}\widehat{u}_{1}+(2-f)s_{1}s_{3}\widehat{u}_{2}+\widehat{F}_{3}\bigr).
\end{aligned}
\end{equation}
It remains to solve \( \eqref{eq:simp_ode}_1 \) in $L^{2}(\mathbb R_{\xi_{2}};\mathbb C^{2})$, uniformly in $(\xi_{1},\xi_{3},\lambda)$; then \( \eqref{eq:simp_ode}_2 \) determines $\widehat{u}_{3}$. In Steps 3--4 we construct the bounded slice-wise inverse, which assembles into the inverse appearing in \eqref{eq:simp_est}.
  
\medskip\noindent\textbf{Step 3 (the slice $\xi_{1}=0$).}
Here the transport \( \xi_1 \partial_{\xi_2}\) drops and $\mathscr L$ acts as multiplication by a matrix whose spectral abscissa is at most $\Lambda$. Consequently, one has the resolvent estimate
\begin{equation}
    \norm{(\lambda-\mathscr{\widehat{L}})^{-1}}\le C\bigl(1+\varepsilon^{-1}\bigr)^{C}
\qq{for any $\lambda$ with}  \operatorname{Re}\lambda\ge\Lambda+\varepsilon.
\end{equation}
However, the slice \( \xi_1 = 0\) is of measure zero in \( \mathbb{R}^3 \).

\medskip\noindent\textbf{Step 4 (the slices \( \xi_1\neq 0\)).}
Recall that we abbreviate $r = \sqrt{\xi_{1}^{2}+\xi_{3}^{2}}$ as in \eqref{def_in_lem_simp_frame}. And the direction $\mathbf{s}=\bm{\xi}/|\bm{\xi}|$ from \eqref{eq:simp_A2}, when fixed \( (\xi_1,\xi_3) \), becomes a function of $\xi_{2}$:
\[
\mathbf{s}(\xi_{2})
=\frac{\bm{\xi}(\xi_{2})}{|\bm{\xi}(\xi_{2})|},
\qquad
|\bm{\xi}(\xi_{2})|=\sqrt{r^{2}+\xi_{2}^{2}},
\]
with $\mathbf{s}(\xi_{2})\to(0,\pm1,0)$, as $\xi_{2}\to\pm\infty$. 

Hereafter, we fix some slice \( \xi_1 > 0, \xi_3 \in \mathbb{R}\), and solve \( \eqref{eq:simp_ode}_1 \).

Let $\Phi(t,s)$ denote the transition matrix of the homogeneous part of \( \eqref{eq:simp_ode}_1 \), and define the Green operator
\begin{equation}\label{eq:simp_Kdef}
(\mathcal G\mathbf{\widehat{F}}_{h})(t)\coloneqq\frac1{\xi_{1}}\int_{t}^{+\infty}\Phi(t,s) \mathbf{\widehat{F}}_{h}(s) \dd{s}  .
\end{equation}
The goal of this step is the following pair of estimates. 
\begin{enumerate}[(a)]
    \item Decay of transition matrix: 
    \begin{equation}\label{eq:simp_kernel}
\norm{\Phi(s,t)}
\le C\bigl(1+\varepsilon^{-1}\bigr)e^{C\varepsilon^{-1}}
\exp\Bigl(-\frac{\varepsilon}{2\xi_{1}}(t-s)
-\frac{\nu}{3\xi_{1}}\bigl(t^{3}-s^{3}\bigr)\Bigr)
\end{equation}
for \( s<t \), uniformly over $(\xi_{1},\xi_{3},\lambda)$ with $\operatorname{Re}\lambda\ge \Lambda+\varepsilon$.
\item The associated Green operator is bounded on $L^{2}(\mathbb R_{\xi_{2}};\mathbb C^{2})$:
\begin{equation}\label{eq:simp_Kbnd}
\norm{\mathcal G\mathbf{\widehat{F}}_{h}}_{L^{2}}\le C\varepsilon^{-1}\bigl(1+\varepsilon^{-1}\bigr)e^{C\varepsilon^{-1}}\norm{\mathbf{\widehat{F}}_{h}}_{L^{2}} .
\end{equation}
That is the resolvent estimate.
\end{enumerate}

We first consider the estimate \eqref{eq:simp_kernel}. Let $\Phi$ denote the transition matrix of the homogeneous part of \( \eqref{eq:simp_ode}_1 \), 
\[
\partial_{s}\Phi(s,t)
=\frac1{\xi_{1}}\Bigl(\lambda+\nu |\bm{\xi}(s)|^{2}I_{2}-A(\mathbf{s}(s))\Bigr)\Phi(s,t),
\qquad
\Phi(t,t)=I_{2},
\qquad s\le t,
\]
namely $\mathbf{\widehat u}_{h}(s)=\Phi(s,t) \mathbf{\widehat u}_{h}(t)$.
And, by \autoref{lem:simp_frame}, $\mathbf{V}\coloneqq \mathcal U^{*}\mathbf{\widehat{u}}_{h}$ obeys
\begin{equation}\label{eq:simp_Wrot}
\partial_{\xi_{2}}\mathbf{V}
=\Bigl[\tfrac1{\xi_{1}}\bigl(\lambda+\nu|\bm{\xi}|^{2}I_{2}-\mathcal{T}\bigr)-E\Bigr]\mathbf{V},
\qquad
\mathcal{T}\coloneqq \mathcal U^{*}A\mathcal U,
\quad
E\coloneqq \mathcal U^{*}\partial_{\xi_{2}}\mathcal U .
\end{equation}
Let $\Phi^{\mathcal U}(t,s)$ denote the transition matrix of \eqref{eq:simp_Wrot}, defined by
\[
\partial_{s}\Phi^{\mathcal U}(s,t)
=\Bigl[\tfrac1{\xi_{1}}\bigl(\lambda+\nu|\bm{\xi}(s)|^{2}I_{2}-\mathcal{T}(s)\bigr)-E(s)\Bigr]\Phi^{\mathcal U}(s,t),
\quad
\Phi^{\mathcal U}(t,t)=I_{2},
\qquad s\le t,
\]
so that every solution of \eqref{eq:simp_Wrot} obeys $\mathbf{V}(s)=\Phi^{\mathcal U}(s,\xi_{2})\mathbf{V}(\xi_{2})$.
Since $\mathcal U$ is unitary,
\begin{equation}\label{eq:simp_unitary}
\norm{\Phi(s,t)}=\norm{\Phi^{\mathcal U}(s,t)}
\qfor s<t,
\end{equation}
Now, it suffices to bound $\Phi^{\mathcal U}$.  

By virtue of \eqref{eq:simp_tri}, we have
\[
\mathcal{T}=\mathcal U^{*}A\mathcal U
=\begin{pmatrix}\lambda_{+}&\chi\\0&\lambda_{-}\end{pmatrix}
\quad\text{is upper triangular, with }
\operatorname{Re}\lambda_{-}\le\operatorname{Re}\lambda_{+},\quad
|\chi|\le C
\]
a.e.\ on $\mathbb R$; moreover, $\mathcal{T}$ is continuous. Hence, setting $E= 0$ in \eqref{eq:simp_Wrot} and applying \autoref{lem:simp_backward} on any interval $[s,t]$ gives
\begin{equation}\label{eq:simp_Efree}
\norm{\Phi^{0}(s,t)}
\le C\Bigl(1+\frac{t-s}{\xi_{1}}\Bigr)
\exp\Bigl(-\frac1{\xi_{1}}\int_{s}^{t}
\operatorname{Re}\bigl[\lambda+\nu|\bm{\xi}|^{2}-\lambda_{+}\bigr] \dd{\xi_{2}'}\Bigr)
\qfor s<t,
\end{equation}
where $\Phi^{0}$ is the transition matrix of \eqref{eq:simp_Wrot} with $E = 0$.

Next, we apply Duhamel principle for the perturbation $E$. The two transition matrices obey
\[
\begin{aligned}
\partial_{\xi_{2}'}\Phi^{\mathcal U}(s,\xi_{2}')
&=-\Phi^{\mathcal U}(s,\xi_{2}')
\Bigl[\tfrac1{\xi_{1}}\bigl(\lambda+\nu|\bm{\xi}(\xi_{2}')|^{2}I_{2}-\mathcal{T}(\xi_{2}')\bigr)-E(\xi_{2}')\Bigr],
\\
\partial_{\xi_{2}'}\Phi^{0}(\xi_{2}',t)
&=\tfrac1{\xi_{1}}\bigl(\lambda+\nu|\bm{\xi}(\xi_{2}')|^{2}I_{2}-\mathcal{T}(\xi_{2}')\bigr)\Phi^{0}(\xi_{2}',t).
\end{aligned}
\]
Hence, the mixed terms cancel in
\[
    \partial_{\xi_{2}'}\bigl[\Phi^{\mathcal U}(s,\xi_{2}')\Phi^{0}(\xi_{2}',t)\bigr]
=\Phi^{\mathcal U}(s,\xi_{2}') E(\xi_{2}') \Phi^{0}(\xi_{2}',t),
\] 
and integrating $\xi_{2}'$ from $s$ to $t$ yields, for $s\le t$,
\begin{equation}\label{eq:simp_duhamel}
\Phi^{\mathcal U}(s,t)
=\Phi^{0}(s,t)
+\int_{s}^{t}\Phi^{\mathcal U}(s,\xi_{2}') 
E(\xi_{2}') \Phi^{0}(\xi_{2}',t) \dd{\xi_{2}'} .
\end{equation}
From \eqref{eq:simp_duhamel} we now bound \(\Phi^{\mathcal U}\). Since \(\operatorname{Re}\lambda\ge\Lambda+\varepsilon\), inequalities \eqref{eq:simp_excesspt}--\eqref{eq:simp_excessint} give
\begin{equation}\label{eq:simp_re_lower}
\begin{aligned}
\operatorname{Re}\bigl[\lambda+\nu|\bm{\xi}|^{2}-\lambda_{+}\bigr] 
&\ge\Lambda+\varepsilon+\nu\bigl(r^{2}+\xi_{2}^{2}\bigr)-\operatorname{Re}\lambda_{+} 
\\
&\ge\varepsilon+\nu\bigl(r^{2}+\xi_{2}^{2}\bigr)-\bigl[\operatorname{Re}\lambda_{+}-\Lambda\bigr]_{+}
\\
&=\varepsilon+\nu\bigl(r^{2}+\xi_{2}^{2}\bigr)-\beta(\mathbf{s}),
\end{aligned}
\end{equation}
 therefore, for \(s\le t\),
\begin{equation}\label{eq:simp_int_lower}
\begin{aligned}
\frac1{\xi_{1}}\int_{s}^{t}
\operatorname{Re}\bigl[\lambda+\nu|\bm{\xi}|^{2}-\lambda_{+}\bigr] \dd{\xi_{2}'}
&\ge\frac{\varepsilon}{\xi_{1}}(t-s)
+\frac{\nu}{\xi_{1}}\int_{s}^{t}\bigl(r^{2}+\xi_{2}'^{2}\bigr)\dd{\xi_{2}'}
-\frac1{\xi_{1}}\int_{s}^{t}\beta\bigl(\mathbf{s}(\xi_{2}')\bigr)\dd{\xi_{2}'}
\\
&\ge\frac{\varepsilon}{\xi_{1}}(t-s)
+\frac{\nu}{3\xi_{1}}\bigl(t^{3}-s^{3}\bigr)
-C,
\end{aligned}
\end{equation}
where the last step used \eqref{eq:simp_excessint}. Inserting \eqref{eq:simp_int_lower} into \eqref{eq:simp_Efree} yields
\begin{equation}\label{eq:simp_Phi0_abs}
\begin{aligned}
\norm{\Phi^{0}(s,t)}
&\le C\Bigl(1+\frac{t-s}{\xi_{1}}\Bigr)
\exp\Bigl(-\frac{\varepsilon}{\xi_{1}}(t-s)-\frac{\nu}{3\xi_{1}}\bigl(t^{3}-s^{3}\bigr)+C\Bigr)
\\
&\le C\bigl(1+\varepsilon^{-1}\bigr)\exp\Bigl(-\frac{\varepsilon}{2\xi_{1}}(t-s)-\frac{\nu}{3\xi_{1}}\bigl(t^{3}-s^{3}\bigr)\Bigr),
\end{aligned}
\end{equation}
the second inequality following from
\[
\bigl(1+u\bigr)e^{-\varepsilon u/2}\le C\bigl(1+\varepsilon^{-1}\bigr)
\qfor u=(t-s)/\xi_{1}\ge 0 .
\]

Set
\begin{equation}\label{eq:simp_gamma}
\gamma(s,t)
\coloneqq\frac{\varepsilon}{2\xi_{1}}(t-s)+\frac{\nu}{3\xi_{1}}\bigl(t^{3}-s^{3}\bigr),
\end{equation}
so that \(\gamma(s,t)=\gamma(s,\xi)+\gamma(\xi,t)\) whenever \(s\le\xi\le t\), and \(\norm{\Phi^{0}(s,t)}\le C\bigl(1+\varepsilon^{-1}\bigr)e^{-\gamma(s,t)}\) by \eqref{eq:simp_Phi0_abs}. Then \eqref{eq:simp_duhamel} and \eqref{eq:simp_Phi0_abs} give
\begin{equation}
\begin{aligned}
e^{\gamma(s,t)}\norm{\Phi^{\mathcal U}(s,t)}
&\le e^{\gamma(s,t)}\norm{\Phi^{0}(s,t)}
+\int_{s}^{t}e^{\gamma(s,t)}\norm{\Phi^{\mathcal U}(s,\xi)}\norm{E(\xi)}\norm{\Phi^{0}(\xi,t)}\dd{\xi}
\\
&\le C\bigl(1+\varepsilon^{-1}\bigr)
+C\bigl(1+\varepsilon^{-1}\bigr)\int_{s}^{t}e^{\gamma(s,\xi)}\norm{\Phi^{\mathcal U}(s,\xi)}\norm{E(\xi)}e^{\gamma(s,t)-\gamma(s,\xi)-\gamma(\xi,t)}\dd{\xi}
\\
&=C\bigl(1+\varepsilon^{-1}\bigr)
+C\bigl(1+\varepsilon^{-1}\bigr)\int_{s}^{t}e^{\gamma(s,\xi)}\norm{\Phi^{\mathcal U}(s,\xi)}\norm{E(\xi)}\dd{\xi}.
\end{aligned}
\end{equation}
Gr\"onwall's inequality and \eqref{eq:simp_budget} yield
\begin{equation}
\begin{aligned}
e^{\gamma(s,t)}\norm{\Phi^{\mathcal U}(s,t)}
&\le C\bigl(1+\varepsilon^{-1}\bigr)\exp\Bigl(C\bigl(1+\varepsilon^{-1}\bigr)\int_{s}^{t}\norm{E(\xi)}\dd{\xi}\Bigr)
\\
&\le C\bigl(1+\varepsilon^{-1}\bigr)\exp\Bigl(C\bigl(1+\varepsilon^{-1}\bigr)\Bigr)
\\
&\le C\bigl(1+\varepsilon^{-1}\bigr)e^{C\varepsilon^{-1}}.
\end{aligned}
\end{equation}
Hence
\begin{equation}\label{eq:simp_kernel_sharp}
\norm{\Phi^{\mathcal U}(s,t)}
\le C\bigl(1+\varepsilon^{-1}\bigr)e^{C\varepsilon^{-1}}\exp\Bigl(-\frac{\varepsilon}{2\xi_{1}}(t-s)-\frac{\nu}{3\xi_{1}}\bigl(t^{3}-s^{3}\bigr)\Bigr),
\qquad s<t,
\end{equation}
and \eqref{eq:simp_kernel} follows from \eqref{eq:simp_unitary}.

Next, we establish \eqref{eq:simp_Kbnd}. From \eqref{eq:simp_Kdef} and \eqref{eq:simp_kernel}, discarding the non-positive cubic viscosity contribution in the exponential,
\begin{equation}
\begin{aligned}
\bigl|(\mathcal G\mathbf{\widehat{F}}_{h})(\xi_{2})\bigr|
&\le\frac1{\xi_{1}}\int_{\xi_{2}}^{+\infty}\norm{\Phi(\xi_{2},\xi_{2}')}
\bigl|\mathbf{\widehat{F}}_{h}(\xi_{2}')\bigr|\dd{\xi_{2}'}
\\
&\le C\bigl(1+\varepsilon^{-1}\bigr)e^{C\varepsilon^{-1}}
\frac1{\xi_{1}}\int_{\xi_{2}}^{+\infty}
\exp\Bigl(-\frac{\varepsilon}{2\xi_{1}}(\xi_{2}'-\xi_{2})\Bigr)
\bigl|\mathbf{\widehat{F}}_{h}(\xi_{2}')\bigr|\dd{\xi_{2}'}.
\end{aligned}
\end{equation}
By Schur's test \cite[Appendix~A.1]{Grafakos2014}, \(\norm{\mathcal G}_{L^{2}\to L^{2}}\) is bounded by any common upper bound of
\begin{equation}
\begin{aligned}
\sup_{\xi_{2}}
\int_{\xi_{2}}^{+\infty}
\frac{C}{\xi_{1}}\bigl(1+\varepsilon^{-1}\bigr)e^{C\varepsilon^{-1}}
\exp\Bigl(-\frac{\varepsilon}{2\xi_{1}}(\xi_{2}'-\xi_{2})\Bigr)
\dd{\xi_{2}'},
\\
\sup_{\xi_{2}'}
\int_{-\infty}^{\xi_{2}'}
\frac{C}{\xi_{1}}\bigl(1+\varepsilon^{-1}\bigr)e^{C\varepsilon^{-1}}
\exp\Bigl(-\frac{\varepsilon}{2\xi_{1}}(\xi_{2}'-\xi_{2})\Bigr)
\dd{\xi_{2}}.
\end{aligned}
\end{equation}
Both term can be bounded by \(C\varepsilon^{-1}\bigl(1+\varepsilon^{-1}\bigr)e^{C\varepsilon^{-1}}\). Indeed, the change of variables \(u=(\xi_{2}'-\xi_{2})/\xi_{1}\) turns either integral into
\begin{equation}
C\bigl(1+\varepsilon^{-1}\bigr)e^{C\varepsilon^{-1}}\int_{0}^{\infty}e^{-\varepsilon u/2}\dd{u}
=C\varepsilon^{-1}\bigl(1+\varepsilon^{-1}\bigr)e^{C\varepsilon^{-1}}.
\end{equation}
Hence
\begin{equation}
\norm{\mathcal G\mathbf{\widehat{F}}_{h}}_{L^{2}}
\le C\varepsilon^{-1}\bigl(1+\varepsilon^{-1}\bigr)e^{C\varepsilon^{-1}}
\norm{\mathbf{\widehat{F}}_{h}}_{L^{2}},
\end{equation}
which is \eqref{eq:simp_Kbnd}, uniformly in \((\xi_{1},\xi_{3},\lambda)\).

It remains to check that \(\mathcal G\mathbf{\widehat{F}}_{h}\) indeed solves the resolvent problem \(\eqref{eq:simp_ode}_{1}\).
Define the horizontal operator by
\begin{equation}\label{eq:simp_slice_h}
\begin{aligned}
\mathcal{D}_{h}
&\coloneqq\Bigl\{\mathbf{\widehat{u}}_{h}\in L^{2}(\mathbb R_{\xi_{2}};\mathbb C^{2}):
|\bm{\xi}|^{2}\mathbf{\widehat{u}}_{h}, \xi_{1}\partial_{\xi_{2}}\mathbf{\widehat{u}}_{h}\in L^{2}\Bigr\},
\\
(\lambda-\mathscr L_{\xi_{1}}^{(h)})\mathbf{\widehat{u}}_{h}
&\coloneqq
(\lambda+\nu|\bm{\xi}|^{2}-A)\mathbf{\widehat{u}}_{h}
-\xi_{1}\partial_{\xi_{2}}\mathbf{\widehat{u}}_{h},
\end{aligned}
\end{equation}
so that \(\eqref{eq:simp_ode}_{1}\) is equivalent to
\begin{equation}
(\lambda-\mathscr L_{\xi_{1}}^{(h)})\mathbf{\widehat{u}}_{h}
=\mathbf{\widehat{F}}_{h}.
\end{equation}
For \(\mathbf{\widehat{F}}_{h}\in\mathcal S\), differentiation under the integral in \eqref{eq:simp_Kdef} yields 
\begin{equation}\label{eq:simp_g_right}
\begin{aligned}
(\lambda-\mathscr L_{\xi_{1}}^{(h)})\mathcal G\mathbf{\widehat{F}}_{h}
&=\mathbf{\widehat{F}}_{h}, \quad
\mathcal G\mathbf{\widehat{F}}_{h}
\in\mathcal{D}_{h}.
\end{aligned}
\end{equation}
By density of \(\mathcal S\) in \(L^{2}(\mathbb R_{\xi_{2}};\mathbb C^{2})\) and the bound \eqref{eq:simp_Kbnd}, the identities \eqref{eq:simp_g_right} extend to every \(\mathbf{\widehat{F}}_{h}\in L^{2}(\mathbb R_{\xi_{2}};\mathbb C^{2})\). Hence
\begin{equation}
(\lambda-\mathscr L_{\xi_{1}}^{(h)})\mathcal G
=I
\quad\text{on }L^{2}(\mathbb R_{\xi_{2}};\mathbb C^{2}),
\qquad
\norm{\mathcal G}_{L^{2}\to L^{2}}
\le C\varepsilon^{-1}\bigl(1+\varepsilon^{-1}\bigr)e^{C\varepsilon^{-1}}.
\end{equation}
If \((\lambda-\mathscr L_{\xi_{1}}^{(h)})\mathbf{\widehat{u}}_{h}=0\) with \(\mathbf{\widehat{u}}_{h}\in\mathcal{D}_{h}\), the integrating factor of \(\eqref{eq:simp_ode}_{1}\) forces the same super-exponential growth as in \autoref{lemma_no_point_spectrum}, whence \(\mathbf{\widehat{u}}_{h}\equiv0\). Therefore
\begin{equation}\label{eq:simp_invbdd}
(\lambda-\mathscr L_{\xi_{1}}^{(h)})^{-1}
=\mathcal G,
\qquad
\norm{(\lambda-\mathscr L_{\xi_{1}}^{(h)})^{-1}}
\le C\varepsilon^{-1}\bigl(1+\varepsilon^{-1}\bigr)e^{C\varepsilon^{-1}},
\end{equation}
uniformly in \(\xi_{1}\ne0\).

The equation \( \eqref{eq:simp_ode}_{2} \) is  
\begin{equation}\label{eq:simp_ode_3}
(\lambda+\nu|\bm{\xi}|^{2})\widehat{u}_{3}-\xi_{1}\partial_{\xi_{2}}\widehat{u}_{3}
=g,
\qquad
g\coloneqq fs_{2}s_{3}\widehat{u}_{1}+(2-f)s_{1}s_{3}\widehat{u}_{2}+\widehat{F}_{3}.
\end{equation}
Since \(|s_{j}|\le1\),
\begin{equation}
|g|
\le C\bigl(|\mathbf{\widehat{u}}_{h}|+|\widehat{F}_{3}|\bigr).
\end{equation}
The associated integrating factor is
\begin{equation}
\Phi_{3}(t,s)
\coloneqq\exp\Bigl(-\frac1{\xi_{1}}\int_{t}^{s}\bigl(\lambda+\nu|\bm{\xi}(\xi_{2}')|^{2}\bigr) \dd{\xi_{2}'}\Bigr)
\qq{for}
t\le s,
\end{equation}
and the Green formula
\begin{equation}\label{eq:simp_g3}
\widehat{u}_{3}(t)
=\frac1{\xi_{1}}\int_{t}^{+\infty}\Phi_{3}(t,s)\,g(s) \dd{s}
\end{equation}
recovers the unique \(L^{2}\) solution of \eqref{eq:simp_ode_3} once \(g\in L^{2}\). From \(\operatorname{Re}\lambda\ge\Lambda+\varepsilon\) one has
\begin{equation}
|\Phi_{3}(t,s)|
\le\exp\Bigl(-\frac{\Lambda+\varepsilon}{\xi_{1}}(s-t)\Bigr),
\end{equation}
after discarding the non-positive viscosity contribution. Hence
\begin{equation}
|\widehat{u}_{3}(t)|
\le\frac1{\xi_{1}}\int_{t}^{+\infty}
\exp\Bigl(-\frac{\Lambda+\varepsilon}{\xi_{1}}(s-t)\Bigr)
|g(s)| \dd{s}.
\end{equation}
Schur's test as in the proof of \eqref{eq:simp_Kbnd}, with the change of variables \(u=(s-t)/\xi_{1}\), yields
\begin{equation}\label{eq:simp_pas_est}
\norm{\widehat{u}_{3}}_{L^{2}}
\le\frac{C}{\Lambda+\varepsilon}\norm{g}_{L^{2}}
\le\frac{C}{\Lambda+\varepsilon}\Bigl(\norm{\mathbf{\widehat{u}}_{h}}_{L^{2}}+\norm{\widehat{F}_{3}}_{L^{2}}\Bigr),
\end{equation}
which completes \eqref{eq:simp_est} when \( \xi_1\neq 0 \).

\medskip\noindent\textbf{Step 5 (assembly and conclusion).}
Define the global operator \(\mathcal G\) on \([L^{2}(\mathbb R^{3})]^{3}\) by assembling the slice inverses of Steps~3--4:
\begin{equation}\label{eq:simp_g_asm}
\begin{aligned}
(\mathcal G\mathbf{\widehat{F}})_{h}(\xi_{2};\xi_{1},\xi_{3})
&=\frac1{\xi_{1}}\int_{\xi_{2}}^{+\infty}\Phi(\xi_{2},s)\mathbf{\widehat{F}}_{h}(s) \dd{s},
\\
(\mathcal G\mathbf{\widehat{F}})_{3}(\xi_{2};\xi_{1},\xi_{3})
&=\frac1{\xi_{1}}\int_{\xi_{2}}^{+\infty}\Phi_{3}(\xi_{2},s)g(s) \dd{s},
\end{aligned}
\qquad\xi_{1}\neq0,
\end{equation}
where \(g\) is as in \eqref{eq:simp_ode_3} with \(\mathbf{\widehat{u}}_{h}=(\mathcal G\mathbf{\widehat{F}})_{h}\), cf.\ \eqref{eq:simp_Kdef} and \eqref{eq:simp_g3}, and
\begin{equation}
(\mathcal G\mathbf{\widehat{F}})(\bm{\xi})
=(\lambda-\mathscr{\widehat{L}}(\bm{\xi}))^{-1}\mathbf{\widehat{F}}(\bm{\xi})
\qquad\text{when }\xi_{1}=0.
\end{equation}
The slice-wise bounds \eqref{eq:simp_invbdd} and \eqref{eq:simp_pas_est} yield
\begin{equation}
\norm{\mathcal G}_{L^2\to L^2}
\le C\varepsilon^{-1}\bigl(1+\varepsilon^{-1}\bigr)e^{C\varepsilon^{-1}}
\qq{whenever}
\operatorname{Re}\lambda\ge\Lambda+\varepsilon.
\end{equation}
On almost every slice the operators of Steps~3--4 are invertible. Hence, for every \(\mathbf{\widehat{F}}\in[L^{2}(\mathbb R^{3})]^{3}\) and every \(\mathbf{\widehat{u}}\in\mathcal D(\mathscr L)\),
\begin{equation}\label{eq:simp_two_sided}
\begin{aligned}
(\lambda-\mathscr L)\mathcal G\mathbf{\widehat{F}}
&=\mathbf{\widehat{F}}, 
\quad \mathcal G(\lambda-\mathscr L)\mathbf{\widehat{u}}
&=\mathbf{\widehat{u}}.
\end{aligned}
\end{equation}
In particular
\begin{equation}
(\lambda-\mathscr L)^{-1}
=\mathcal G,
\qquad
\norm{(\lambda-\mathscr L)^{-1}}
\le C\varepsilon^{-1}\bigl(1+\varepsilon^{-1}\bigr)e^{C\varepsilon^{-1}}
\qq{for}
\operatorname{Re}\lambda\ge\Lambda+\varepsilon,
\end{equation}
which is \eqref{eq:simp_est}. Combined with \eqref{eq:simp_gp} one obtains \(\omega_{0}(\mathscr L)\le\Lambda\). And, Step~1 supplies \(\omega_{0}(\mathscr L)\ge\Lambda\), whence \(\omega_{0}(\mathscr L)=\Lambda\). In particular, for every \(\varepsilon>0\) there exists \(C_{\varepsilon}>0\) such that \eqref{eq:simp_semigroup_bnd} holds.
The proof is completed.
\end{proof}

\section{The linear instability}\label{sec_linear_instability}
In this section we construct a linearly unstable solution of the system \eqref{Linear_eq_abstract_form} under the assumption that \( 0<f < 1\), i.e. \( \R(f) < 0 \). In this case it is known that the \( x \)-independent system is unstable.

Assume \( \mathbf{u} \) is \( x \)-independent. Recalling \eqref{definition_of_L}, the Fourier transform of \eqref{Linear_eq_abstract_form} gives that
\begin{equation}\label{Linear_system_x_independent_Fourier_form}
	\begin{aligned}
		\partial_t \mathbf{\widehat{u}}
		=
		\mathscr{\hat{L}}_{0} (\xi_2,\xi_3) \mathbf{\widehat{u}}
		\coloneqq
		 \begin{pmatrix}
       - \nu \abs{\bm{\xi}}^2   & f - 1 & 0 \\ 
       -f + f \frac{\xi_2^2}{ \abs{\bm{\xi}}^2 }   &  - \nu \abs{\bm{\xi}}^2  & 0\\ 
       f \frac{\xi_2 \xi_3}{ \abs{\bm{\xi}}^2 } & 0 & - \nu \abs{\bm{\xi}}^2  \\
    \end{pmatrix} 
		\mathbf{\widehat{u}},
	\end{aligned}
\end{equation}
where $\abs{\bm{\xi}}^2=\xi_2^2+\xi_3^2$ and \( \mathscr{\hat{L}}_{0}(\xi_2,\xi_3)  \equiv \mathscr{\hat{L}}(\bm{\xi})|_{\xi_1 = 0} \). 
The eigenvalues and eigenfunctions are as follows
\begin{equation}\label{streak_unstable_system}
	\begin{aligned}
		\lambda_1(\xi_2,\xi_3) & =  - \nu (\xi_2^2 + \xi_3^2) + \sqrt{f(1 - f)  \frac{\xi_3^2}{ (\xi_2^2 + \xi_3^2)} } , \\ \mathbf{\widehat{u}}_1(\xi_2,\xi_3) &= \qty(-\frac{ \abs{\bm{\xi}}  \sqrt{f(1 - f) \xi_3^2 }}{f }, - \xi_3^2,\xi_2 \xi_3),   \\
		\lambda_2(\xi_2,\xi_3) & =  - \nu (\xi_2^2 + \xi_3^2) - \sqrt{f(1 - f)  \frac{\xi_3^2}{ (\xi_2^2 + \xi_3^2)} } , \\\mathbf{\widehat{u}}_2(\xi_2,\xi_3) &= \qty( \frac{ \abs{\bm{\xi}}   \sqrt{ f(1 - f) \xi_3^2}}{f }, - \xi_3^2, \xi_2 \xi_3), \\ 
        \lambda_3(\xi_2,\xi_3) &= - \nu (\xi_2^2 + \xi_3^2), \quad \mathbf{\widehat{u}}_3(\xi_2,\xi_3) = \qty(0,0,1).
	\end{aligned}
\end{equation}
It is clear that if 
$
	f(1 - f) > \nu^2 \frac{(\xi_2^2 + \xi_3^2)^3}{\xi_3^2},
$
the eigenvalue \( \lambda_1 > 0 \), hence   the system admits exponentially growing solution.
Since the eigenfunctions are \( x \)-independent, they are not integrable in \( \mathbb{R}^3 \). Thus, the unstable eigenfunctions above are not eigenfunctions in \( [L^2(\mathbb{R}^3)]^3 \), and they do not imply the existence of unstable eigenvalues for \( \mathscr{L} \). However, the above unstable eigenfunction assists us to prove the instability of the 3D system.

For the general case when \( \mathbf{u} \) is \( x \)-dependent, the following statement for \( \varepsilon \)-pseudo-spectra of \( \mathscr{L} \) is true.

\begin{lemma}\label{lemma_pseudo_spectrum_of_Linear_Op}
	Let \( \lambda_{j}(\xi_2,\xi_3) \), \( j = 1,2,3\) be the eigenvalues of the matrix \( \mathscr{\hat{L}}_{0}(\xi_2,\xi_3) \) defined in \eqref{Linear_system_x_independent_Fourier_form}--\eqref{streak_unstable_system}. If \( 0 < f < 1\), then for any \( \varepsilon > 0 \) and \( (\xi_2,\xi_3) \in \mathbb{R}^2 \setminus \qty{\bm{0}} \), the eigenvalues \( \lambda_{j}(\xi_2,\xi_3) \), \( j = 1,2,3\) lie in the \( \varepsilon \)-pseudo-spectrum of \( \mathscr{L} \), i.e. \( \lambda_{j}(\xi_2,\xi_3) \in \sigma_{\varepsilon}(\mathscr{L}) \) for \( j = 1 ,2,3\).
	Moreover,  \( \lambda_{j}(\xi_2,\xi_3) \in \sigma(\mathscr{L}) \) for \( j = 1,2,3\).
\end{lemma}
\begin{proof}

	The proof is the same for  \( j = 1 ,2,3\).  We provide the details  for \( \lambda_{1} (\xi_2,\xi_3) \).

	 By virtue of Plancherel's identity, for any \( \mathbf{u} \in \mathcal{D}(\mathscr{L}) \subset [L^2(\mathbb{R}^3)]^3 \),
	\begin{equation}
		\norm{ \qty(\lambda_1(\xi_2,\xi_3)  - \mathscr{L}) \mathbf{u} }_{L^2} =  \norm{(\lambda_1(\xi_2,\xi_3)  - \hat{\mathscr{L}}) \hat{\mathbf{u}} }_{L^2}.
	\end{equation}
	Hence, we can show \( \lambda_1 (\xi_2,\xi_3) \) is an \( \varepsilon \)-pseudo-eigenvalue by finding a vector-valued function  \( \hat{\mathbf{u}}_{\varepsilon} (\bm{\xi}) \in [L^2(\mathbb{R}^3)]^3 \) satisfying \( \norm{\hat{\mathbf{u}}_{\varepsilon}}_{L^2} = 1   \) and \( \mathcal{F}^{-1}\hat{\mathbf{u}}_{\varepsilon} \in \mathcal{D} (\mathscr{L}) \) such that
	\begin{equation}\label{linear_analysis_1}
		\norm{(\lambda_1(\xi_2,\xi_3)  - \hat{\mathscr{L}}) \hat{\mathbf{u}}_{\varepsilon} }_{L^2} < \varepsilon .
	\end{equation}

	Apply the Fourier transformation to \eqref{Linear_eq_abstract_form}, one obtains
	\begin{equation}
		\begin{aligned}
			\partial_t \mathbf{\widehat{u}}
			=
			 \begin{pmatrix}
        - \nu \abs{\bm{\xi}}^2 + \xi_1 \partial_{\xi_2} + f \frac{\xi_1 \xi_2}{ \abs{\bm{\xi}}^2 }& f - 1 + (2 - f) \frac{\xi_1^2}{ \abs{\bm{\xi}}^2 } & 0\\ 
        - f + f \frac{\xi_2^2}{ \abs{\bm{\xi}}^2 }& - \nu \abs{\bm{\xi}}^2 + \xi_1 \partial_{\xi_2} + (2 - f) \frac{\xi_1 \xi_2}{ \abs{\bm{\xi}}^2 }& 0\\ 
        f \frac{\xi_2 \xi_3}{ \abs{\bm{\xi}}^2 }& (2 - f) \frac{\xi_1 \xi_3}{ \abs{\bm{\xi}}^2 }& - \nu \abs{\bm{\xi}}^2 + \xi_1 \partial_{\xi_2}
     \end{pmatrix} 
            \widehat{\mathbf{u}} .
		\end{aligned}
	\end{equation}
	One introduces  fixed frequencies
	\begin{equation}\label{frequencies_1}
		\bm{\xi}_1^{*} = (0,\xi_2^{*},\xi_3^{*}) ,\quad  \bm{\xi}_2^{*} = (0, -\xi_2^{*},\xi_3^{*}) ,\quad \bm{\xi}_3^{*} = (0,\xi_2^{*}, -\xi_3^{*}) ,\quad \bm{\xi}_4^{*} = (0, -\xi_2^{*}, -\xi_3^{*}) ,
	\end{equation}
	such that \(\xi_2^{*}, \xi_3^{*} > 0\). For the degenerate case \( \xi_2^{*} \xi_3^{*} = 0\) and \( (\xi_2^{*},\xi_3^{*}) \neq (0,0) \), we only need to fix two frequencies. For example, if \( \xi_2^{*} = 0\) and \( \xi_3^{*} \neq 0 \), we define
	\begin{equation}\label{frequencies_2}
		\bm{\xi}_1^{*} = (0,0,\xi_3^{*}) ,\quad  \bm{\xi}_2^{*} = (0,0, -\xi_3^{*}) .
	\end{equation}
	Hereafter, we only prove the non-degenerate case, the proof of the degenerate case is similar.

    
    Let \( \theta (x) \) be the standard mollifier  with compact support in the interval \( [ - 1,1]\), and \( \theta_{\delta}(x) \) be defined by
	\begin{equation}\label{1D_mollifier}
		\theta_{\delta}(x) \coloneqq \frac{1}{\delta} \theta \qty(\frac{x}{\delta}), \quad \delta \in (0, 1].
	\end{equation}
	Then
	\begin{equation}
		\eta(\bm{\xi}) \coloneqq\theta (\xi_1)\theta (\xi_2)\theta (\xi_3),\quad \eta_{\delta,\delta'}(\bm{\xi}) \coloneqq \theta_{\delta}  (\xi_1)\theta_{\delta'}  (\xi_2) \theta_{\delta} (\xi_3), \quad \delta' \in (0,1], 
	\end{equation}
	and
	\begin{equation}\label{construction_of_eigenfunction}
		\hat{\mathbf{u}}_{\delta,\delta'}(\bm{\xi};\xi_2^{*},\xi_3^{*}) \coloneqq \sum_{j = 1}^{4} \eta_{\delta,\delta'} (\bm{\xi} - \bm{\xi}_{j}^{*})\hat{\mathbf{u}}_{1}(\xi_2^{*},\xi_3^{*}) ,
	\end{equation}
	where \( \hat{\mathbf{u}}_{1}(\xi_2^{*},\xi_3^{*})  \) is the eigenvector of the \( \mathscr{\hat{L}}_{0} \) corresponding to the eigenvalue \( \lambda_1(\xi_2^{*},\xi_3^{*}) \) given by \eqref{streak_unstable_system}. Hereafter we use the shorthand notations \( \lambda_1 = \lambda_1(\xi_2^{*},\xi_3^{*}) \), \( \hat{\mathbf{u}}_{1}  =\hat{\mathbf{u}}_{1}(\xi_2^{*},\xi_3^{*}) \) and \( \hat{\mathbf{u}}_{\delta,\delta'} (\bm{\xi})= \hat{\mathbf{u}}_{\delta,\delta'}(\bm{\xi};\xi_2^{*},\xi_3^{*}) \), with the understanding that the frequency \( (\xi_2^{*},\xi_3^{*}) \) is fixed.

	Notice that by construction  \(  \hat{\mathbf{u}}_{\delta,\delta'}(\bm{\xi}) \) is even in variables \( \xi_{j} \), \( j = 1,2,3\). Hence   \( \mathbf{u}_{\delta,\delta'} (\mathbf{x}) = \mathcal{F}^{-1} \hat{\mathbf{u}}_{\delta,\delta'} (\mathbf{x}) \) is a real function. Moreover, since \( \hat{\mathbf{u}}_{\delta,\delta'} \) is compactly supported and smooth, \( \mathbf{u}_{\delta,\delta'}(\mathbf{x}) \) lies in the Schwarz space \( [\mathcal{S}(\mathbb{R}^3)]^3 \), and thus in \( \mathcal{D}(\mathscr{L}) \).

		We verify that  \( \mathbf{u}_{\delta,\delta'}(\mathbf{x}) \)  is indeed an \( \varepsilon \)-pseudo-eigenfunction corresponding to \( \lambda_1 \) for suitably small \( \delta \) and \( \delta' \), namely \eqref{linear_analysis_1} holds for \( \hat{\mathbf{u}}_{\delta,\delta'}(\bm{\xi}) \).

	We assume without loss of generality that for any \( \delta,\delta' \in (0,1] \), the support of \( \eta_{\delta,\delta'} (\bm{\xi} - \bm{\xi}_{j}^{*}) \) is mutually disjoint. Then, one gets by change of variables
	\begin{equation}
		\norm{ \hat{\mathbf{u}}_{\delta,\delta'} }_{L^2}^2 = \sum_{j = 1}^{4} \int_{\mathbb{R}^3}^{} \abs{\eta_{\delta,\delta'} (\bm{\xi} - \bm{\xi}_{j}^{*}) \hat{\mathbf{u}}_{1} }^2  \dd{\bm{\xi}} = \frac{4}{\delta^2 \delta'}  \abs{\hat{\mathbf{u}}_{1} }^2 \norm{\eta(\bm{\xi})}_{L^2}^2.
	\end{equation}

	Furthermore, recalling \eqref{definition_of_L}, a direct calculation shows for \( \lambda_1 \) that
	\begin{equation}
		\begin{aligned}
		\frac{1}{\norm{\hat{\mathbf{u}}_{\delta,\delta'}}_{L^2} }	\norm{ (\lambda_1  - \hat{\mathscr{L}}) \hat{\mathbf{u}}_{\delta,\delta'} }_{L^2}   \leq {} & \frac{1}{\norm{\hat{\mathbf{u}}_{\delta,\delta'}}_{L^2} }
        \sum_{j = 1}^{4}  \bigg( \norm{\eta_{\delta,\delta'}(\bm{\xi} - \bm{\xi}_{j}^{*}) \qty[ \qty(\lambda_1  - \mathscr{\widehat{L}}_0(\bm{\xi}) - \mathscr{\widehat{L}}_1(\bm{\xi}) )\hat{\mathbf{u}}_{1}] }_{L^2}     \\
        & {}  + \norm{\xi_1 \partial_{\xi_2}  \eta_{\delta,\delta'}(\bm{\xi} - \bm{\xi}_{j}^{*}) \hat{\mathbf{u}}_{1} }_{L^2}  \bigg) \\
			\eqqcolon                                                                                                                            & I_1 + I_2.
		\end{aligned}
	\end{equation}
	 Denote \( h (\bm{\xi}) = \qty(\lambda_1  - \mathscr{\widehat{L}}_0 (\bm{\xi}) - \mathscr{\widehat{L}}_1(\bm{\xi}) ) \frac{\hat{\mathbf{u}}_{1}}{ \abs{\hat{\mathbf{u}}_{1}} }\). Then, taking the symmetry of \( \eta_{\delta,\delta'}(\bm{\xi} ) \) and \( \abs{\bm{\xi}}^2  \) into account, we have
	\begin{equation}
		\begin{aligned}
			I_1
			\leq {} & \frac{1}{\norm{\hat{\mathbf{u}}_{\delta,\delta'}}_{L^2} }\sum_{j = 1}^{4} \norm{\eta_{\delta,\delta'} (\bm{\xi} - \bm{\xi}_{j}^{*})}_{L^2}   \norm{  \qty(\lambda_1  - \mathscr{\widehat{L}}_0 (\bm{\xi}) - \mathscr{\widehat{L}}_1(\bm{\xi}) )\hat{\mathbf{u}}_{1}}_{L^{ \infty}(\supp \eta_{\delta,\delta'}   )}                                                      \\
			=  {} & \frac{4}{\norm{\hat{\mathbf{u}}_{\delta,\delta'}}_{L^2} }  \norm{\eta_{\delta,\delta'} (\bm{\xi} )}_{L^2} \abs{\hat{\mathbf{u}}_{1}}   \norm{ h (\bm{\xi}) }_{L^{ \infty}( [ -\delta,\delta] \times [ \xi_2^{*} -\delta',\xi_2^{*} +\delta'] \times [ \xi_3^{*}-\delta,\xi_3^{*} +\delta]  )}                                                                \\
			=  {} &
			\frac{4}{\frac{2}{\delta \delta'^{\frac{1}{2}}}  \abs{\mathbf{\hat{u}}_{1} } \norm{\eta(\bm{\xi})}_{L^2}  } \frac{1}{\delta \delta'^{\frac{1}{2}} }\norm{\eta (\bm{\xi} )}_{L^2} \abs{\hat{\mathbf{u}}_{1}}  \norm{ h (\bm{\xi}) }_{L^{ \infty}( [ -\delta,\delta] \times [ \xi_2^{*} -\delta',\xi_2^{*} +\delta'] \times [ \xi_3^{*}-\delta,\xi_3^{*} +\delta]  )} \\
			=  {} &
			2 \norm{ h (\bm{\xi}) }_{L^{ \infty}( [ -\delta,\delta] \times [ \xi_2^{*} -\delta',\xi_2^{*} +\delta'] \times [ \xi_3^{*}-\delta,\xi_3^{*} +\delta]  )},
		\end{aligned}
	\end{equation}
    where we utilized the fact that
	\begin{equation}
		\norm{\eta_{\delta,\delta'} (\bm{\xi} - \bm{\xi}_{j}^{*})}_{L^2} = \frac{1}{\delta \delta'^{\frac{1}{2}} }\norm{\eta (\bm{\xi} )}_{L^2}.
	\end{equation}
	One observes that by the construction of \( \hat{\mathbf{u}}_{1} \), the function \( h (\bm{\xi}) \) is a smooth function in the neighborhood of \( \bm{\xi}_{j}^{*} \), \( j = 1,2,3,4\) and satisfies \( h (\bm{\xi}_{j}^{*}) = \bm{0} \). Thus, for any \( \varepsilon > 0\), there exists \( \delta' > 0 \) and \( \delta_0(\delta',\varepsilon) > 0\) such that for any \( 0 < \delta < \delta_0(\delta',\varepsilon) \),
	\begin{equation}
		I_1 < \frac{\varepsilon}{2} .
	\end{equation}

	For \( I_2 \), we have
	\begin{equation}
		\begin{aligned}
			I_2 =  {} & \frac{1}{\norm{\hat{\mathbf{u}}_{\delta,\delta'}}_{L^2} }\sum_{j = 1}^{4} \norm{\xi_1 \partial_{\xi_2}  \eta_{\delta,\delta'} (\bm{\xi} - \bm{\xi}_{j}^{*}) \hat{\mathbf{u}}_{1}}_{L^2}                      \\
			\leq  {} & \frac{4}{\frac{2}{\delta \delta'^{\frac{1}{2}}}  \abs{\hat{\mathbf{u}}_{1} } \norm{\eta(\bm{\xi})}_{L^2}  }\abs{\hat{\mathbf{u}}_{1}} \norm{\xi_1 \partial_{\xi_2}  \eta_{\delta,\delta'}(\bm{\xi} ) }_{L^2}, \\
		\end{aligned}
	\end{equation}
	where
	\begin{equation}
		\begin{aligned}
			&\norm{\xi_1 \partial_{\xi_2}  \eta_{\delta,\delta'} (\bm{\xi} ) }_{L^2}^2 
             =  \int_{\mathbb{R}^3}^{} \abs{\xi_1 \partial_{\xi_2} \qty(\frac{1}{\delta^2\delta'}  \eta \qty( \frac{\xi_1}{\delta}, \frac{\xi_2}{\delta'}, \frac{\xi_3}{\delta} ))}^2  \dd{\bm{\xi}}          \\
	   & = \int_{\mathbb{R}^3}^{} \abs{\frac{\xi_1}{\delta'}  \frac{1}{\delta^2\delta'}  \partial_{\xi_2}\eta \qty( \frac{\xi_1}{\delta}, \frac{\xi_2}{\delta'}, \frac{\xi_3}{\delta} )}^2  \dd{\bm{\xi}}  \leq \frac{1}{\delta^2\delta'} \abs{\frac{\delta}{\delta'} }^2 \int_{\mathbb{R}^3}^{} \abs{ \xi_1\partial_{\xi_2}  \eta \qty(\bm{\xi} )}^2  \dd{\bm{\xi}}.                                          \\
		\end{aligned}
	\end{equation}
	Hence
	\begin{equation}
		\begin{aligned}
			I_2 \leq {} & \frac{4}{\frac{2}{\delta \delta'^{\frac{1}{2}}}  \abs{\hat{\mathbf{u}}_{1} } \norm{\eta(\bm{\xi})}_{L^2}  } \abs{\hat{\mathbf{u}}_{1}}  \frac{1}{\delta\delta'^{\frac{1}{2}}} \frac{\delta}{\delta'}  \norm{\xi_1\partial_{\xi_2}  \eta \qty(\bm{\xi} )}_{L^2} \\
			\leq   {}  & \frac{2  \norm{\xi_1\partial_{\xi_2}  \eta \qty(\bm{\xi} )}_{L^2}  }{ \norm{\eta(\bm{\xi})}_{L^2}}    \frac{\delta}{\delta'}   \leq C  \frac{\delta}{\delta'},
		\end{aligned}
	\end{equation}
	for a constant \( C > 0\) determined by \( \eta \).
Therefore, for any \( \varepsilon > 0\), there exists \( \delta' > 0 \) and \( \delta_1(\varepsilon,\delta') > 0\) such that for any \( 0 < \delta < \delta_1(\varepsilon,\delta')\)
	\begin{equation}
		I_2 < \frac{\varepsilon}{2}.
	\end{equation}

	Combining the estimates for \( I_1 \) and \( I_2 \), we have for any \( \varepsilon > 0 \), there exists small enough \( \delta' > 0 \), such that for any \( \delta < \min \qty{\delta_1(\varepsilon,\delta'),\delta_2(\varepsilon,\delta')} \) 
	\begin{equation}
		\norm{(\lambda_1  - \hat{\mathscr{L}}) \frac{ \hat{\mathbf{u}}_{\delta,\delta'}}{ \norm{ \hat{\mathbf{u}}_{\delta,\delta'}}_{L^2}  } }_{L^2} < \varepsilon.
	\end{equation}
	This  implies that  \( \lambda_1(\xi_2^\ast,\xi_3^\ast) \) is an \( \varepsilon \)-pseudo-eigenvalue of \( \mathscr{L} \) with corresponding \( \varepsilon \)-pseudo-eigenfunction \( \mathbf{u}_{\delta,\delta'} \).

	Finally,  by \autoref{theorem_pspectral_properties}, one has \( \cap_{\varepsilon > 0} \sigma_{\varepsilon}(\mathscr{L}) = \sigma(\mathscr{L}) \). The proof is complete.
\end{proof}
 
\begin{remark}
The anisotropic scaling $I_{1}$ and $I_{2}$ is conceptually related to the technique used in the geometric-optics and short-wave instability theory of Lifschitz--Hameiri \cite{Lifschitz1991,Lifschitz1992} and Friedlander--Vishik \cite{Friedlander1991,Friedlander1992} and in the bicharacteristic-amplitude system and essential spectrum theory of Shvydkoy \cite{Shvydkoy2006,ShvydkoyLatushkin2005}. There the $x$-independent matrix $\mathscr{\hat L}_{0}$ is the amplitude equation along bicharacteristics, and localization near $\xi_{1}=0$ corresponds to wave packets almost orthogonal to the shear. Our $\delta\ll\delta'$ mollification makes this picture rigorous in the $L^{2}$ pseudo-spectral framework and explains the concentration near the invariant plane $\xi_{1}=0$ where $\mathscr{\widehat L}_{1}$ vanishes.
\end{remark}

\begin{corollary}\label{not_small_of_u2}
    Under the conditions in \autoref{lemma_pseudo_spectrum_of_Linear_Op}, for any \( (\xi_2^{*},\xi_3^{*}) \in  \mathbb{R}^2 \setminus \qty{\bm{0}} \), \( \xi_3^{*}\neq 0 \), there exist an \( \varepsilon \)-pseudo-eigenfunction of \( \lambda_1(\xi_2^{*},\xi_3^{*}) \), denoted by \( \mathbf{u}_\varepsilon = (u_{1,\varepsilon},u_{2,\varepsilon},u_{3,\varepsilon}) \), such that
\begin{equation}\label{ratio_of_u}
        u_{1,\varepsilon}(\mathbf{x}) =   \frac{ \sqrt{ (\xi_2^{*})^2 + (\xi_3^{*})^2}  }{ \abs{\xi_3^{*}}  } \frac{\sqrt{1 - f }}{\sqrt{f}}  u_{2,\varepsilon}(\mathbf{x}), \quad 
    u_{3,\varepsilon}(\mathbf{x}) =  - \frac{ \xi_2^{*}  }{ \xi_3^{*}  }  u_{2,\varepsilon}(\mathbf{x}),
    \end{equation}
    for all \( \mathbf{x} \in \mathbb{R}^3 \). Moreover, the support of \( \widehat{\mathbf{u}}_\varepsilon \) satisfies 
    \begin{equation}\label{supp_cond}
        \supp \mathbf{\widehat{u}}_\varepsilon \subseteq \bigcup_{j = 1}^{4} B \qty(\bm{\xi}_j^{*},\frac{1}{2})
    \end{equation}
    where \( \bm{\xi}_j^{*} ,j = 1,2,3,4\) are defined in \eqref{frequencies_1} and \( B \qty(\bm{\xi}_j^{*} ,\frac{1}{2}) \) is the ball centered at \(\bm{\xi}_j^{*}\) of radius \( \frac{1}{2} \). 

\end{corollary}
\begin{proof}
    The proof of \autoref{lemma_pseudo_spectrum_of_Linear_Op} shows that, there exists a small enough \( \delta' > 0 \) and \( \delta < \min \qty{\delta_1(\varepsilon,\delta'),\delta_2(\varepsilon,\delta')} \) such that \( \mathcal{F}^{-1} \hat{\mathbf{u}}_{\delta,\delta'}(\mathbf{x};\xi_2^{*},\xi_3^{*})\) is an \( \varepsilon \)-pseudo-eigenfunction of \( \lambda_1(\xi_2^{*},\xi_3^{*}) \). Denote the three components of \( \hat{\mathbf{u}}_{\delta,\delta'}(\bm{\xi};\xi_2^{*},\xi_3^{*}) \) by 
    \[
        \hat{\mathbf{u}}_{\delta,\delta'}(\bm{\xi};\xi_2^{*},\xi_3^{*}) = (\widehat{u}_{1,\delta,\delta'},\widehat{u}_{2,\delta,\delta'},\widehat{u}_{3,\delta,\delta'})(\bm{\xi};\xi_2^{*},\xi_3^{*}) .
    \]  
    Then by the construction \eqref{construction_of_eigenfunction} and \eqref{streak_unstable_system}, we have 
    \begin{equation}\label{linear_analysis_{12}}
        \begin{aligned}
     &\widehat{u}_{1,\delta,\delta'}(\bm{\xi}) =   \frac{ \sqrt{ (\xi_2^{*})^2 + (\xi_3^{*})^2}  }{ \abs{\xi_3^{*}}  } \frac{\sqrt{1 - f }}{\sqrt{f}}  \widehat{u}_{2,\delta,\delta'}(\bm{\xi}), \quad 
    \widehat{u}_{3,\delta,\delta'}(\bm{\xi}) =  - \frac{ \xi_2^{*}  }{ \xi_3^{*}  }  \widehat{u}_{2,\delta,\delta'}(\bm{\xi}), \\ 
\end{aligned}
    \end{equation}
    for all \( \bm{\xi} \in \mathbb{R}^3 \) when \( \xi_3^{*}\neq 0 \). 
    
    Then, from \eqref{linear_analysis_{12}} and the inverse Fourier transform, we have
    \begin{equation}\label{linear_analysis_{13}}
        u_{1,\delta,\delta'}  (\mathbf{x}) =   \frac{ \sqrt{ (\xi_2^{*})^2 + (\xi_3^{*})^2}  }{ \abs{\xi_3^{*}}  } \frac{\sqrt{1 - f }}{\sqrt{f}}  u_{2,\delta,\delta'}(\mathbf{x}), \quad 
    u_{3,\delta,\delta'}(\mathbf{x}) =  - \frac{ \xi_2^{*}  }{ \xi_3^{*}  }  u_{2,\delta,\delta'}(\mathbf{x}),
    \end{equation}
     for all \( \mathbf{x} \in \mathbb{R}^3 \). Moreover, \eqref{linear_analysis_{13}} is valid for all \( \delta,\delta' \).

     The support of \( \hat{\mathbf{u}}_{\delta,\delta'} \) satisfies
    \begin{equation}
   \supp \hat{\mathbf{u}}_{\delta,\delta'}\subseteq [ -\delta,\delta] \times A_{2,\delta'}\times A_{3,\delta}.
    \end{equation}
    where \(A_{2,\delta'}\) and \(A_{3,\delta}\) are defined by
    \[
    \begin{aligned}
   A_{2,\delta'} &= [ \xi_2^{*}-\delta',\xi_2^{*} +\delta'] \cup  [ -\xi_2^{*}-\delta', -\xi_2^{*} +\delta'] ,\\
    A_{3,\delta} &= [ \xi_3^{*}-\delta,\xi_3^{*} +\delta] \cup  [ -\xi_3^{*}-\delta, -\xi_3^{*} +\delta] .
    \end{aligned}
    \]
And for any \( 0 <\delta'_1 < \delta' \) and \( 0 < \delta_1 < \delta\), the function \( \mathcal{F}^{-1} \hat{\mathbf{u}}_{\delta_1,\delta'_1}(\mathbf{x};\xi_2^{*},\xi_3^{*}) \) is also the \( \varepsilon \)-pseudo-eigenfunction of \( \lambda_1 (\xi_2^{*},\xi_3^{*})\), hence \eqref{ratio_of_u} and \eqref{supp_cond} hold for \( \mathbf{u}_{\varepsilon}(\mathbf{x}) = \mathcal{F}^{-1} \hat{\mathbf{u}}_{\delta_1,\delta'_1}(\mathbf{x};\xi_2^{*},\xi_3^{*})  \) with small enough \( \delta_1 \) and \( \delta_1' \). This completes the proof.
\end{proof} 
\begin{corollary}\label{coro_spec_lower_bound}
	Suppose \( 0 < f < 1\), then for any real number \( \lambda \leq  \sqrt{f(1 - f) }  \), \( \lambda \in \sigma(\mathscr{L})  \).
\end{corollary}
\begin{proof}
	The \autoref{lemma_pseudo_spectrum_of_Linear_Op} shows that \( \lambda_1 \) defined in \eqref{streak_unstable_system} lies in the spectrum of \( \mathscr{L} \). By the definition of \( \lambda_1 \), one has
	\begin{equation}
		\overline{ \qty{ \lambda_1(\xi_2,\xi_3) | \xi_2,\xi_3 \in \mathbb{R}^2 \setminus \{\bm{0}\}  } } = ( - \infty, \sqrt{f(1 - f) } ].
	\end{equation}
	This completes the proof.
\end{proof} 
\begin{lemma}
	\label{lemma:linear_instability}
	Suppose \( 0 < f < 1\), then for any \( T > 0 \) and \( \epsilon > 0 \), there exists a $ \mathbf{u}_{T,\epsilon} = (u_{1,T,\epsilon},u_{2,T,\epsilon},u_{3,T,\epsilon}) \in [L^{2}(\mathbb{R}^{3})]^{3}$ such that
    \begin{enumerate}[{\rm (i)}]
        \item $\div\mathbf u_{T,\epsilon}=0$ and $\mathbf u_{T,\epsilon}\in[\mathcal S(\mathbb{R}^{3})]^{3}$ is real-valued;
        \item for all \( t\in[0,T] \)
    \begin{equation}\label{linear_analysis_7}
		\norm{ \qty(e^{t \Lambda } -e^{t \mathscr{L}})\mathbf{u}_{T,\epsilon}}_{L^2}  \leq  \epsilon\norm{\mathbf{u}_{T,\epsilon}}_{L^2} ;
	\end{equation}
    \item for all \( t\in[0,T] \), \( e^{t \mathscr{L} }\mathbf{u}_{T,\epsilon} \in [H^k(\mathbb{R}^3)]^3 \) for any \( t\in[0,T] \), \( k \geq 0 \); 
    \item \( \supp \widehat{\mathbf{u}_{T,\epsilon}} \subseteq B(0,1)  \), where \( B(0,1) \) is the unit ball at origin.
    \end{enumerate}  
\end{lemma}
\begin{proof}
    \noindent Throughout (i)--(iv), fix \(\xi_{3}^{*}\neq 0\) and set
    \begin{equation}\label{eq_freq_notation}
        r_{*}\coloneqq\abs{\xi_{3}^{*}},\qquad
        \widehat{\mathbf{u}}_{1}\coloneqq\widehat{\mathbf{u}}_{1}(0,\xi_{3}^{*}),\qquad
        \lambda_{1}\coloneqq\lambda_{1}(0,\xi_{3}^{*}),
    \end{equation}
    where \(\widehat{\mathbf{u}}_{1}\) and \(\lambda_{1}\) are given by \eqref{streak_unstable_system} with \(\xi_{2}^{*}=0\).
    All constants \(C,C_{k},\ldots\) below depend only on \(f\) and the mollifier \(\eta\), and are independent of \(r_{*}\), \(\delta\), \(\delta'\), \(T\), and \(\epsilon\).

    \medskip\noindent\textit{Proof of (i).}
    Let \(\mathbf{u}_{\delta,\delta'}\) be the \(\gamma\)-pseudo-eigenfunction from \autoref{lemma_pseudo_spectrum_of_Linear_Op} with frequencies \eqref{frequencies_2}. Then \(\mathbf{u}_{\delta,\delta'}\in[\mathcal{S}(\mathbb{R}^{3})]^{3}\) is real-valued, but is not a priori divergence-free. We project it onto the divergence-free subspace and control the error in eight steps.

    \medskip\noindent\textbf{Step 1 (Leray projection and error term).}
    The Leray projector \(\mathbb{P}\) has the multiplier form
    \begin{equation}
        \widehat{\mathbb{P}\mathbf{u}}(\bm{\xi})=m(\bm{\xi})\widehat{\mathbf{u}}(\bm{\xi}),\qquad
        m(\bm{\xi})=I-\frac{\bm{\xi}\bm{\xi}^{\intercal}}{\abs{\bm{\xi}}^{2}}.
    \end{equation}
    Set \(\mathbf{w}_{\delta,\delta'}\coloneqq\mathbb{P}\mathbf{u}_{\delta,\delta'}-\mathbf{u}_{\delta,\delta'}\). Then
    \begin{equation}\label{eq_w_hat}
        \widehat{\mathbf{w}}_{\delta,\delta'}(\bm{\xi})
        =\bigl(m(\bm{\xi})-I\bigr)\widehat{\mathbf{u}}_{\delta,\delta'}(\bm{\xi})
        =-\sum_{j}\eta_{\delta,\delta'}(\bm{\xi}-\bm{\xi}_{j}^{*}) \frac{\bm{\xi}\bm{\xi}^{\intercal}}{\abs{\bm{\xi}}^{2}}\widehat{\mathbf{u}}_{1}.
    \end{equation}

    \medskip\noindent\textbf{Step 2 (Multiplier bounds on the support).}
    For \(\delta<\frac{r_{*}}{2}\) and \(\delta'<\frac{r_{*}}{2}\), one has \(\abs{\bm{\xi}}\ge\frac{r_{*}}{2}\) on \(\supp\widehat{\mathbf{u}}_{\delta,\delta'}\), hence \(m\) is smooth there. Since \(m_{ij}=\delta_{ij}-\xi_{i}\xi_{j}/\abs{\bm{\xi}}^{2}\),
    \begin{equation}
        \partial_{\xi_{k}}m_{ij}
        =-\frac{\delta_{ik}\xi_{j}+\xi_{i}\delta_{jk}}{\abs{\bm{\xi}}^{2}}
        +\frac{2\xi_{i}\xi_{j}\xi_{k}}{\abs{\bm{\xi}}^{4}},
    \end{equation}
    and on the support \(\abs{\partial_{\xi_{k}}m_{ij}}\le 3/\abs{\bm{\xi}}\le 6/r_{*}\). Therefore
    \begin{equation}\label{eq_m_bounds}
        \norm{m}_{L^{\infty}(\supp\widehat{\mathbf{u}}_{\delta,\delta'})}\le 2,\qquad
        \norm{\partial_{\xi_{k}}m}_{L^{\infty}(\supp\widehat{\mathbf{u}}_{\delta,\delta'})}\le \frac{C}{r_{*}},\qquad k=1,2,3.
    \end{equation}

    \medskip\noindent\textbf{Step 3 (Exact identity for \((m-I)\widehat{\mathbf{u}}_{1}\)).}
    By \eqref{streak_unstable_system} and \eqref{frequencies_2},
    \begin{equation}
        \widehat{\mathbf{u}}_{1}=\qty(-\frac{r_{*}\sqrt{f(1-f) r_{*}^{2}}}{f},- r_{*}^{2},0),
    \end{equation}
    and therefore \(\bm{\xi}_{j}^{*}\cdot\widehat{\mathbf{u}}_{1}=0\) for \(\bm{\xi}_{j}^{*}=(0,0,\pm\xi_{3}^{*})\). Consequently
    \begin{equation}\label{eq_m_minus_I_u1}
        \begin{aligned}
            m(\bm{\xi}_{j}^{*})\widehat{\mathbf{u}}_{1}&=\widehat{\mathbf{u}}_{1}, 
            \quad
            \bigl(m(\bm{\xi})-I\bigr)\widehat{\mathbf{u}}_{1}
             =-\frac{\bigl((\bm{\xi}-\bm{\xi}_{j}^{*})\cdot\widehat{\mathbf{u}}_{1}\bigr)\bm{\xi}}{\abs{\bm{\xi}}^{2}},
        \end{aligned}
    \end{equation}
    On \(\supp\eta_{\delta,\delta'}( \cdot -\bm{\xi}_{j}^{*})\) one has \(\abs{\bm{\xi}-\bm{\xi}_{j}^{*}}\le\delta+\delta'\) and \(\abs{\bm{\xi}}\ge r_{*}/2\), hence
    \begin{equation}\label{eq_m_minus_I_u1_bound}
        \abs{\bigl(m(\bm{\xi})-I\bigr)\widehat{\mathbf{u}}_{1}}\le \frac{2(\delta+\delta')}{r_{*}}\abs{\widehat{\mathbf{u}}_{1}}.
    \end{equation}

    \medskip\noindent\textbf{Step 4 (\(L^{2}\) bounds for \(\mathbf{w}_{\delta,\delta'}\)).}
    Since \(\widehat{\mathbf{w}}_{\delta,\delta'}=(m-I)\widehat{\mathbf{u}}_{\delta,\delta'}\) and \(\abs{\widehat{\mathbf{u}}_{\delta,\delta'}}\le\abs{\widehat{\mathbf{u}}_{1}}\) on the support,
    \begin{equation}\label{eq_w_L2_bound}
        \begin{aligned}
            \norm{\mathbf{w}_{\delta,\delta'}}_{L^{2}}
            &\le \frac{C(\delta+\delta')}{r_{*}}\norm{\mathbf{u}_{\delta,\delta'}}_{L^{2}},\\
            \norm{\mathbb{P}\mathbf{u}_{\delta,\delta'}}_{L^{2}}
            &\ge \Bigl(1-\frac{C(\delta+\delta')}{r_{*}}\Bigr)\norm{\mathbf{u}_{\delta,\delta'}}_{L^{2}}.
        \end{aligned}
    \end{equation}

    \medskip\noindent\textbf{Step 5 (Regularity and reality of \(\mathbb{P}\mathbf{u}_{\delta,\delta'}\)).}
    Since \(m\) is smooth near \(\supp\widehat{\mathbf{u}}_{\delta,\delta'}\), \(\widehat{\mathbb{P}\mathbf{u}_{\delta,\delta'}}\) is smooth and compactly supported.
    The mollifier \(\theta\) is even, hence \(\eta_{\delta,\delta'}(-\bm{\xi})=\eta_{\delta,\delta'}(\bm{\xi})\); with \(\{\bm{\xi}_{j}^{*}\}=\{-\bm{\xi}_{j}^{*}\}\) from \eqref{frequencies_2},
    \begin{equation}
        \widehat{\mathbf{u}}_{\delta,\delta'}(-\bm{\xi})=\widehat{\mathbf{u}}_{\delta,\delta'}(\bm{\xi}),\qquad
        m(-\bm{\xi})=m(\bm{\xi}),
    \end{equation}
    and therefore \(\widehat{\mathbb{P}\mathbf{u}_{\delta,\delta'}}(-\bm{\xi})=\widehat{\mathbb{P}\mathbf{u}_{\delta,\delta'}}(\bm{\xi})\).
    Since \(m\) and \(\widehat{\mathbf{u}}_{\delta,\delta'}\) are real-valued, \(\mathbb{P}\mathbf{u}_{\delta,\delta'}\) is real-valued and \(\mathbb{P}\mathbf{u}_{\delta,\delta'}\in \mathcal{S}_\sigma(\mathbb{R}^{3})\).

    \medskip\noindent\textbf{Step 6 (\(H^{k}\) bounds for \(\mathbf{w}_{\delta,\delta'}\)).}
    By Plancherel's identity,
    \begin{equation}
        \norm{\mathbf{w}_{\delta,\delta'}}_{H^{k}}^{2}
        =\sum_{\abs{\bm{\alpha}}\le k}\norm{D^{\bm{\alpha}}\mathbf{w}_{\delta,\delta'}}_{L^{2}}^{2}
        =\sum_{\abs{\bm{\alpha}}\le k}\norm{\bm{\xi}^{\bm{\alpha}}\widehat{\mathbf{w}}_{\delta,\delta'}}_{L^{2}}^{2}.
    \end{equation}
    From \eqref{eq_m_minus_I_u1}, on each \(\supp\eta_{\delta,\delta'}( \cdot -\bm{\xi}_{j}^{*})\),
    \begin{equation}
        \widehat{\mathbf{w}}_{\delta,\delta'}(\bm{\xi})
        =-\eta_{\delta,\delta'}(\bm{\xi}-\bm{\xi}_{j}^{*}) 
        \frac{\bigl((\bm{\xi}-\bm{\xi}_{j}^{*})\cdot\widehat{\mathbf{u}}_{1}\bigr) \bm{\xi}}{\abs{\bm{\xi}}^{2}},
    \end{equation}
    so \(\widehat{\mathbf{w}}_{\delta,\delta'}\) is supported in a union of balls of radius \(\le\delta+\delta'\) centred at the \(\bm{\xi}_{j}^{*}\).
    Since \(\delta,\delta'<\frac{r_{*}}{2}\), one has \(\abs{\bm{\xi}}\le 2r_{*}\) on \(\supp\widehat{\mathbf{w}}_{\delta,\delta'}\).
    For each \(\abs{\bm{\alpha}}\le k\), therefore
    \begin{equation}
        \norm{\bm{\xi}^{\bm{\alpha}}\widehat{\mathbf{w}}_{\delta,\delta'}}_{L^{2}}
        \le (2r_{*})^{\abs{\bm{\alpha}}}\norm{\widehat{\mathbf{w}}_{\delta,\delta'}}_{L^{2}}.
    \end{equation}
    The case \(\abs{\bm{\alpha}}=0\) is \eqref{eq_w_L2_bound}.
    For \(\abs{\bm{\alpha}}\ge 1\), combining the display above with \eqref{eq_w_L2_bound} gives
    \begin{equation}
        \norm{\bm{\xi}^{\bm{\alpha}}\widehat{\mathbf{w}}_{\delta,\delta'}}_{L^{2}}
        \le C_{\abs{\bm{\alpha}}}(2r_{*})^{\abs{\bm{\alpha}}}\frac{\delta+\delta'}{r_{*}}\norm{\mathbf{u}_{\delta,\delta'}}_{L^{2}}.
    \end{equation}
    Summing over \(\abs{\bm{\alpha}}\le k\) yields
    \begin{equation}\label{eq_w_Hk_bound}
        \norm{\mathbf{w}_{\delta,\delta'}}_{H^{k}}
        \le C_{k}(2r_{*})^{k}\frac{\delta+\delta'}{r_{*}}\norm{\mathbf{u}_{\delta,\delta'}}_{L^{2}},\qquad k\ge 0.
    \end{equation}
    In particular, when \(r_{*}\le 1\) one has \((2r_{*})^{k}\le 2^{k}\), so the right-hand side is \(\le \tilde{C}_{k}(\delta+\delta')/r_{*} \norm{\mathbf{u}_{\delta,\delta'}}_{L^{2}}\) with \(\tilde{C}_{k}\) depending only on \(k\), \(f\), and \(\eta\).

    \medskip\noindent\textbf{Step 7 (Error estimate for \((\lambda_{1}-\mathscr{L})\mathbf{w}_{\delta,\delta'}\)).}
    Recalling \eqref{definition_of_L}, \(\widehat{\mathscr{L}}=\mathscr{\widehat{L}}_{0}+\mathscr{\widehat{L}}_{1}+\mathscr{\widehat{L}}_{2}\) with \(\mathscr{\widehat{L}}_{2}=\xi_{1}\partial_{\xi_{2}}\), we have
    \begin{equation}
        \bigl(\lambda_{1}-\widehat{\mathscr{L}}\bigr)\widehat{\mathbf{w}}_{\delta,\delta'}
        =\bigl(\lambda_{1}-\mathscr{\widehat{L}}_{0}(\bm{\xi})-\mathscr{\widehat{L}}_{1}(\bm{\xi})\bigr)\widehat{\mathbf{w}}_{\delta,\delta'}
        -\xi_{1}\partial_{\xi_{2}}\widehat{\mathbf{w}}_{\delta,\delta'},
    \end{equation}
    with
    \begin{equation}\label{eq_error_1}
        \begin{aligned}
            \xi_{1}\partial_{\xi_{2}}\widehat{\mathbf{w}}_{\delta,\delta'}
            &=\sum_{j}\xi_{1}\bigl(\partial_{\xi_{2}}m\bigr)\eta_{\delta,\delta'}(\bm{\xi}-\bm{\xi}_{j}^{*})\widehat{\mathbf{u}}_{1}\\
            &\qquad+\sum_{j}\xi_{1}\bigl(m-I\bigr)\partial_{\xi_{2}}\eta_{\delta,\delta'}(\bm{\xi}-\bm{\xi}_{j}^{*})\widehat{\mathbf{u}}_{1}.
        \end{aligned}
    \end{equation}
    By \eqref{eq_m_bounds} and \eqref{construction_of_eigenfunction},
    \begin{equation}
        \sum_{j}\eta_{\delta,\delta'}(\bm{\xi}-\bm{\xi}_{j}^{*})\widehat{\mathbf{u}}_{1}=\widehat{\mathbf{u}}_{\delta,\delta'},
    \end{equation}
    and on \(\supp\widehat{\mathbf{u}}_{\delta,\delta'}\) one has \(\abs{\xi_{1}}\le\delta\), therefore
    \begin{equation}\label{eq_error_1_partial_m_bound}
        \begin{aligned}
            \norm{\sum_{j}\xi_{1}\bigl(\partial_{\xi_{2}}m\bigr)\eta_{\delta,\delta'}(\bm{\xi}-\bm{\xi}_{j}^{*})\widehat{\mathbf{u}}_{1}}_{L^{2}}
            &\le\delta\norm{\partial_{\xi_{2}}m}_{L^{\infty}(\supp\widehat{\mathbf{u}}_{\delta,\delta'})}\norm{\widehat{\mathbf{u}}_{\delta,\delta'}}_{L^{2}}\\
            &\le \frac{C\delta}{r_{*}}\norm{\widehat{\mathbf{u}}_{\delta,\delta'}}_{L^{2}}.
        \end{aligned}
    \end{equation}
    For the second summand in \eqref{eq_error_1}, \(\partial_{\xi_{2}}\eta_{\delta,\delta'}\) is scalar and \(\widehat{\mathbf{u}}_{1}\) is independent of \(\bm{\xi}\), hence
    \begin{equation}
        \xi_{1}\bigl(m(\bm{\xi})-I\bigr)\partial_{\xi_{2}}\eta_{\delta,\delta'}(\bm{\xi}-\bm{\xi}_{j}^{*})\widehat{\mathbf{u}}_{1}
        =\xi_{1}\partial_{\xi_{2}}\eta_{\delta,\delta'}(\bm{\xi}-\bm{\xi}_{j}^{*}) \bigl(m(\bm{\xi})-I\bigr)\widehat{\mathbf{u}}_{1}.
    \end{equation}
    On \(\supp\eta_{\delta,\delta'}( \cdot -\bm{\xi}_{j}^{*})\), the same argument as for \eqref{eq_m_minus_I_u1_bound} yields, for every \(\mathbf{v}\in\mathbb{R}^{3}\),
    \begin{equation}
        \abs{\bigl(m(\bm{\xi})-I\bigr)\mathbf{v}}\le \frac{2(\delta+\delta')}{r_{*}}\abs{\mathbf{v}}.
    \end{equation}
    With \(\mathbf{v}=\partial_{\xi_{2}}\eta_{\delta,\delta'}(\bm{\xi}-\bm{\xi}_{j}^{*})\widehat{\mathbf{u}}_{1}\), and since \(\abs{\xi_{1}}\le\delta\) on \(\supp\widehat{\mathbf{u}}_{\delta,\delta'}\), therefore
    \begin{equation}
        \begin{aligned}
            \abs{\xi_{1}\bigl(m-I\bigr)\partial_{\xi_{2}}\eta_{\delta,\delta'}(\bm{\xi}-\bm{\xi}_{j}^{*})\widehat{\mathbf{u}}_{1}}
            &=\abs{\xi_{1}} \abs{\bigl(m-I\bigr)\partial_{\xi_{2}}\eta_{\delta,\delta'}(\bm{\xi}-\bm{\xi}_{j}^{*})\widehat{\mathbf{u}}_{1}}\\
            &\le \frac{2(\delta+\delta')}{r_{*}}\abs{\xi_{1}} 
            \abs{\partial_{\xi_{2}}\eta_{\delta,\delta'}(\bm{\xi}-\bm{\xi}_{j}^{*})\widehat{\mathbf{u}}_{1}}\\
            &\le \frac{2\delta(\delta+\delta')}{r_{*}}
            \abs{\partial_{\xi_{2}}\eta_{\delta,\delta'}(\bm{\xi}-\bm{\xi}_{j}^{*})\widehat{\mathbf{u}}_{1}}.
        \end{aligned}
    \end{equation}
    Integrating and summing over \(j\),
    \begin{equation}
        \begin{aligned}
            \norm{\sum_{j}\xi_{1}\bigl(m-I\bigr)\partial_{\xi_{2}}\eta_{\delta,\delta'}(\bm{\xi}-\bm{\xi}_{j}^{*})\widehat{\mathbf{u}}_{1}}_{L^{2}}
            &\le\sum_{j}\norm{\xi_{1}\bigl(m-I\bigr)\partial_{\xi_{2}}\eta_{\delta,\delta'}(\bm{\xi}-\bm{\xi}_{j}^{*})\widehat{\mathbf{u}}_{1}}_{L^{2}}\\
            &\le \frac{2\delta(\delta+\delta')}{r_{*}}\sum_{j}
            \norm{\partial_{\xi_{2}}\eta_{\delta,\delta'}(\bm{\xi}-\bm{\xi}_{j}^{*})\widehat{\mathbf{u}}_{1}}_{L^{2}}.
        \end{aligned}
    \end{equation}
    Since \(\widehat{\mathbf{u}}_{1}\) is independent of \(\bm{\xi}\),
    \begin{equation}
        \sum_{j}\norm{\partial_{\xi_{2}}\eta_{\delta,\delta'}(\bm{\xi}-\bm{\xi}_{j}^{*})\widehat{\mathbf{u}}_{1}}_{L^{2}}
        =\abs{\widehat{\mathbf{u}}_{1}}\sum_{j}\norm{\partial_{\xi_{2}}\eta_{\delta,\delta'}(\bm{\xi}-\bm{\xi}_{j}^{*})}_{L^{2}}.
    \end{equation}
    The same scaling as in the estimate of \(I_{2}\) in the proof of \autoref{lemma_pseudo_spectrum_of_Linear_Op} yields
    \begin{equation}
        \sum_{j}\norm{\partial_{\xi_{2}}\eta_{\delta,\delta'}(\bm{\xi}-\bm{\xi}_{j}^{*})}_{L^{2}}
        \le \frac{C}{\delta'}\norm{\eta(\bm{\xi})}_{L^{2}},
    \end{equation}
    and \(\abs{\widehat{\mathbf{u}}_{1}}\le C\norm{\widehat{\mathbf{u}}_{\delta,\delta'}}_{L^{2}}\), hence
    \begin{equation}\label{eq_error_1_mI_partial_eta_bound}
        \norm{\sum_{j}\xi_{1}\bigl(m-I\bigr)\partial_{\xi_{2}}\eta_{\delta,\delta'}(\bm{\xi}-\bm{\xi}_{j}^{*})\widehat{\mathbf{u}}_{1}}_{L^{2}}
        \le \frac{C\delta(\delta+\delta')}{r_{*}\delta'}\norm{\widehat{\mathbf{u}}_{\delta,\delta'}}_{L^{2}}.
    \end{equation}
    Combining \eqref{eq_error_1_partial_m_bound} and \eqref{eq_error_1_mI_partial_eta_bound} in \eqref{eq_error_1},
    \begin{equation}\label{eq_xi1_d2_w_bound}
        \norm{\xi_{1}\partial_{\xi_{2}}\widehat{\mathbf{w}}_{\delta,\delta'}}_{L^{2}}
        \le C\qty(\frac{\delta}{r_{*}}+\frac{\delta(\delta+\delta')}{r_{*}\delta'})\norm{\widehat{\mathbf{u}}_{\delta,\delta'}}_{L^{2}}.
    \end{equation}
    On \(\supp\widehat{\mathbf{u}}_{\delta,\delta'}\) one has \(\abs{\bm{\xi}}\le 2r_{*}\) for \(\delta,\delta'<r_{*}/2\), hence
    \begin{equation}\label{eq_h_infty}
        \norm{\lambda_{1}-\mathscr{\widehat{L}}_{0}(\bm{\xi})-\mathscr{\widehat{L}}_{1}(\bm{\xi})}_{L^{\infty}(\supp\widehat{\mathbf{u}}_{\delta,\delta'})}
        \le C\qty(\nu r_{*}^{2}+1+\abs{\lambda_{1}}).
    \end{equation}
    Thus
    \begin{equation}
        \begin{aligned}
            \norm{\bigl(\lambda_{1}-\mathscr{\widehat{L}}_{0}(\bm{\xi})-\mathscr{\widehat{L}}_{1}(\bm{\xi})\bigr)\widehat{\mathbf{w}}_{\delta,\delta'}}_{L^{2}}
            &\le C\qty(\nu r_{*}^{2}+1+\abs{\lambda_{1}})\norm{\mathbf{w}_{\delta,\delta'}}_{L^{2}}\\
            &\le \frac{C\qty(\nu r_{*}^{2}+1+\abs{\lambda_{1}})(\delta+\delta')}{r_{*}}\norm{\mathbf{u}_{\delta,\delta'}}_{L^{2}}.
        \end{aligned}
    \end{equation}
    Therefore, by \eqref{eq_xi1_d2_w_bound} and \eqref{eq_h_infty},
    \begin{equation}
        \norm{\bigl(\lambda_{1}-\mathscr{L}\bigr)\mathbf{w}_{\delta,\delta'}}_{L^{2}}
        \le C\qty(\frac{\delta}{r_{*}}+\frac{\delta(\delta+\delta')}{r_{*}\delta'}+\frac{(\nu r_{*}^{2}+1+\abs{\lambda_{1}})(\delta+\delta')}{r_{*}})\norm{\mathbf{u}_{\delta,\delta'}}_{L^{2}}.
    \end{equation}
    Given \(\gamma>0\), \autoref{lemma_pseudo_spectrum_of_Linear_Op} yields for small enough \(\delta,\delta'\) such that
    \begin{equation}
        \norm{\bigl(\lambda_{1}-\mathscr{L}\bigr)\mathbf{u}_{\delta,\delta'}}_{L^{2}}
        \le\frac{\gamma}{4}\norm{\mathbf{u}_{\delta,\delta'}}_{L^{2}}.
    \end{equation}
    Shrinking \(\delta'\) and then \(\delta\ll\delta'\) if necessary gives
    \begin{equation}
        C\qty(\frac{\delta}{r_{*}}+\frac{\delta(\delta+\delta')}{r_{*}\delta'}+\frac{(\nu r_{*}^{2}+1+\abs{\lambda_{1}})(\delta+\delta')}{r_{*}})\le\frac{\gamma}{4},\qquad
        \frac{C(\delta+\delta')}{r_{*}}\le\frac{1}{2},
    \end{equation}
    hence
    \begin{equation}
        \norm{\bigl(\lambda_{1}-\mathscr{L}\bigr)\mathbb{P}\mathbf{u}_{\delta,\delta'}}_{L^{2}}
        \le\frac{\gamma}{2}\norm{\mathbf{u}_{\delta,\delta'}}_{L^{2}}
        \le\gamma\norm{\mathbb{P}\mathbf{u}_{\delta,\delta'}}_{L^{2}}.
    \end{equation}

    \medskip\noindent\textbf{Step 8 (Normalization).}
    Define \(\mathbf{u}_{\gamma}\coloneqq\mathbb{P}\mathbf{u}_{\delta,\delta'}/\norm{\mathbb{P}\mathbf{u}_{\delta,\delta'}}_{L^{2}}\).
    Then \(\mathbf{u}_{\gamma}\) is a \(\gamma\)-pseudo-eigenfunction, \(\div\mathbf{u}_{\gamma}=0\), and \(\mathbf{u}_{\gamma}\in[\mathcal{S}(\mathbb{R}^{3})]^{3}\) is real-valued.
    We work with this \(\mathbf{u}_{\gamma}\) in the sequel.

    Proof of (ii):
    First, for any \( \mathbf{u} \in [L^2(\mathbb{R}^3)]^3 \) and \( \zeta  \in \mathbb{C} \),
	\begin{equation}
		\norm{e^{t \Lambda }\mathbf{u}- e^{t \mathscr{L}} \mathbf{u}}_{L^2}  \leq  \norm{e^{t \Lambda }\mathbf{u} - e^{t \zeta  }\mathbf{u}}_{L^2} + \norm{e^{t \zeta  }\mathbf{u}- e^{t \mathscr{L}} \mathbf{u}}_{L^2}  .
	\end{equation}
	Hence, the conclusion \eqref{linear_analysis_7} is valid if
	\begin{equation}\label{linear_analysis_8}
		\norm{e^{t \Lambda }\mathbf{u} - e^{t \zeta  }\mathbf{u}}_{L^2} \leq \frac{\epsilon}{2} \norm{\mathbf{u}}_{L^2} ,\quad  \norm{e^{t \zeta  }\mathbf{u}- e^{t \mathscr{L}} \mathbf{u}}_{L^2} \leq \frac{\epsilon}{2} \norm{\mathbf{u}}_{L^2}
	\end{equation}
	holds for all \( t\in[0,T] \) and for suitable \( \zeta  \) and \( \mathbf{u} \).

	Second, let us choose \(\xi_{3}^{*}=\xi_{3}^{\epsilon}\) and \(\lambda_{1}=\lambda_{1}(0,\xi_{3}^{\epsilon})\), where \(\xi_{3}^{\epsilon}\) is a small enough positive number satisfying
    \begin{equation}\label{linear_analysis_{11}}
        \abs{\xi_3^{\epsilon}} \leq  
        \begin{cases}
        \sqrt{ 
      -\frac{1}{\nu T} \ln\left(1 - \frac{\epsilon}{2e^{T \Lambda}}\right)  }, & \text{when} \quad  \epsilon < 2e^{T \Lambda}, \\ 
       + \infty, & \text{when} \quad  \epsilon \geq  2e^{T \Lambda}, 
        \end{cases}
    \end{equation}
   Then we have
	\begin{equation}\label{linear_analysis_{15}}
		e^{t \Lambda }- e^{t  \lambda_{1}} \leq \frac{\epsilon}{2}, \qq{for all $t\in[0,T]$,}
	\end{equation}
	hence
    \begin{equation}
        \norm{e^{t \Lambda }\mathbf{u} - e^{t \lambda_1 }\mathbf{u}}_{L^2} \leq \frac{\epsilon}{2} \norm{\mathbf{u}}_{L^2}.
    \end{equation}
    Hence, we now fix some \( \xi_3^{\epsilon }\) satisfying \eqref{linear_analysis_{11}}, then \(  \zeta = \lambda_{1} \) meets the first requirement in \eqref{linear_analysis_8} for all \( \mathbf{u} \in [L^2(\mathbb{R}^3)]^3 \).

	Third, with the help of \eqref{eq:simp_semigroup_bnd} from \autoref{thm:simp_growth}, we are able to estimate the difference of \( e^{t\mathscr{L}} \mathbf{u}_{\gamma} \) and \( e^{t \lambda_{1} } \mathbf{u}_{\gamma} \) with \( \mathbf{u}_{\gamma} \) being the \( \gamma \)-pseudo-eigenfunction corresponds to \( \lambda_{1}  \). Fix \(\varepsilon_0>0\). Recalling \eqref{streak_unstable_system}, we have \( \lambda_{1} =\lambda_{1}(0,\xi_3^{\epsilon}) \leq \Lambda \). We have  
	\begin{equation}
		\begin{aligned}
		     &\norm{(e^{t\mathscr{L}} - e^{t \lambda_{1}  }) \mathbf{u}_{\gamma}}_{L^2}
			  = \norm{e^{t \lambda_{1} } (e^{t(\mathscr{L}-  \lambda_{1} )} - 1) \mathbf{u}_{\gamma}}_{L^2}
			\\&=\norm{e^{ t \lambda_{1} } \int_{0}^{t} e^{s(\mathscr{L}-  \lambda_{1} )}(\mathscr{L} - \lambda_{1} )   \mathbf{u}_{\gamma} \dd{s}}_{L^2}                         \\
			&=  \norm{ \int_{0}^{t} e^{ (t- s) \lambda_{1}  }e^{s\mathscr{L}}(\mathscr{L} - \lambda_{1} )   \mathbf{u}_{\gamma} \dd{s}}_{L^2}                          
			 \\&\leq   \int_{0}^{t} e^{ (t- s) \lambda_{1}  } \norm{e^{s\mathscr{L}}}_{L^2 \to L^2}   \norm{(\mathscr{L} - \lambda_{1} )\mathbf{u}_{\gamma}}_{L^2}   \dd{s} \\
			 & \leq  C_{\varepsilon_0}\int_{0}^{t} e^{ t \lambda_{1}  } e^{s \qty(\Lambda+\varepsilon_0 - \lambda_{1} )}   \norm{(\mathscr{L} - \lambda_{1} )\mathbf{u}_{\gamma}}_{L^2}   \dd{s}             
			 \leq  C_{\varepsilon_0}\gamma \norm{\mathbf{u}_{\gamma}}_{L^2}   e^{ t \lambda_{1}  } \frac{1}{\Lambda+\varepsilon_0 - \lambda_{1} } e^{s \qty(\Lambda+\varepsilon_0 - \lambda_{1} )}   \Big|^{t}_{0}       \\                        
			  & \leq  C_{\varepsilon_0}\gamma \norm{\mathbf{u}_{\gamma}}_{L^2}   e^{ t \lambda_{1}  } \frac{1}{\Lambda+\varepsilon_0 - \lambda_{1} } \qty(e^{t \qty(\Lambda+\varepsilon_0 - \lambda_{1} )} - 1)   \\
		\end{aligned}
	\end{equation}
    Thus, using \eqref{linear_analysis_{15}}
    \begin{equation}
        \begin{aligned}
			 \norm{(e^{t\mathscr{L}} - e^{t \lambda_{1}  }) \mathbf{u}_{\gamma}}_{L^2} & 
             \le C_{\varepsilon_0}\gamma \norm{\mathbf{u}_{\gamma}}_{L^2} \frac{e^{t\lambda_1}(e^{t\varepsilon_0} - 1) + \frac{\epsilon}{2} e^{t\varepsilon_0}}{\varepsilon_0}.
        \end{aligned}
    \end{equation}
    Then we can choose \( \gamma \) small enough such that \( C_{\varepsilon_0}\gamma  \frac{e^{t\lambda_1}(e^{t\varepsilon_0} - 1) + \frac{\epsilon}{2} e^{t\varepsilon_0}}{\varepsilon_0} < \frac{\epsilon}{2}\)  for all \( t \in [0,T] \). Therefore, we have 
    \begin{equation}
        \norm{(e^{t\mathscr{L}} - e^{t \lambda_{1}  }) \mathbf{u}_{\gamma}}_{L^2}< \frac{\epsilon}{2} \norm{\mathbf{u}_{\gamma}}_{L^2}.
    \end{equation} 
    Then \eqref{linear_analysis_8} holds for \( \zeta = \lambda_{1}(0,\xi_3^{\epsilon}) \) and \( \mathbf{u}_\gamma \), which implies \eqref{linear_analysis_7} with \( \mathbf{u}_{T,\epsilon} = \mathbf{u}_{\gamma} \).  

      Proof of (iii): we prove \( e^{t \mathscr{L}} \mathbf{u}_{\gamma} \in \bigcap_{k \geq 0} [H^k(\mathbb{R}^3)]^3 \) for all \( t \in[0,T] \). For this matter, we recall that \( \mathcal{D}(\mathscr{L}^{\infty}) \) is an invariant set of the semigroup \( e^{t \mathscr{L}} \) for all \( t \geq 0 \) \cite[Chap.~II, Lemma~1.9]{Engel2008}. Then, by the fact that
     \begin{equation}
        \mathbf{u}_{\gamma} \in [\mathcal{S}(\mathbb{R}^3)]^3 \subseteq \mathcal{D}(\mathscr{L}^{\infty}) \subseteq \bigcap_{k \geq 0}[H^k(\mathbb{R}^3)]^3  ,
     \end{equation} 
     we have \( e^{t \mathscr{L}} \mathbf{u}_{\gamma} \in \bigcap_{k \geq 0}[H^k(\mathbb{R}^3)]^3   \) for all \( t \geq 0 \).

    Proof of (iv): In (ii) we fix \(\xi_{3}^{*}=\xi_{3}^{\epsilon}\) with \(\abs{\xi_{3}^{\epsilon}}\) small enough to satisfy \eqref{linear_analysis_{11}}, hence \(r_{*}=\abs{\xi_{3}^{*}}\le 1\).
    By choosing \(\epsilon\) small enough, one may also assume \(r_{*}<\frac{1}{2}\); then, with \(\delta,\delta'<\frac{r_{*}}{2}\), the support of \(\widehat{\mathbf{u}}_{\delta,\delta'}\) is contained in
    \begin{equation}
        \bigcup_{j}B(\bm{\xi}_{j}^{*},\tfrac{r_{*}}{2})\subseteq B\qty(0,\tfrac{3r_{*}}{2})\subseteq B(0,1).
    \end{equation}
    The Leray projector preserves this compact support and the \(H^{k}\) regularity, so (iv) holds for the divergence-free \(\mathbf{u}_{\gamma}\) as well. 


The proof is now complete.
\end{proof}

\noindent Proof of \autoref{thm:Linear_Instability}.
Fix \(T>0\) and \(\epsilon>0\).
Apply \autoref{lemma:linear_instability} with \(2\epsilon\) in place of \(\epsilon\), and let \(\mathbf{u}_{T,\epsilon}\) be the divergence-free field from part~(i) of \autoref{lemma:linear_instability}.
Set
\begin{equation}
    \mathbf{u}^{T,\epsilon}(t)\coloneqq e^{t\mathscr{L}}\mathbf{u}_{T,\epsilon}.
\end{equation}
Since \(\mathbf{u}_{T,\epsilon}\in[\mathcal{S}(\mathbb{R}^{3})]^{3}\subseteq\mathcal{D}(\mathscr{L})\) and \(\div\mathbf{u}_{T,\epsilon}=0\), \(\mathbf{u}^{T,\epsilon}\) is a real-valued solution of \eqref{linearized_perturbed_equation}--\eqref{boundary_conditions_for_perturbed_equation} with \(\div\mathbf{u}^{T,\epsilon}(t)=0\) and an associated pressure \(q^{T,\epsilon}\) determined by \eqref{linearized_perturbed_equation}.
By part~(iii) of \autoref{lemma:linear_instability}, \(e^{t\mathscr{L}}\mathbf{u}_{T,\epsilon}\in[H^{k}(\mathbb{R}^{3})]^{3}\) for every \(t\in[0,T]\) and \(k\ge 0\); since \(e^{t\mathscr{L}}\) is a \(C_{0}\)-semigroup on \([L^{2}(\mathbb{R}^{3})]^{3}\), one has \(\mathbf{u}^{T,\epsilon}\in C([0,T];[H^{k}(\mathbb{R}^{3})]^{3})\) and \(\nabla q^{T,\epsilon}\in C([0,T];[H^{k}(\mathbb{R}^{3})]^{3})\).

\medskip\noindent\textbf{Part (i).}
Since \(\mathbf{u}^{T,\epsilon}(0)=\mathbf{u}_{T,\epsilon}\)
\begin{equation}
    \begin{aligned}
        \norm{\mathbf{u}^{T,\epsilon}(t)}_{L^{2}}
        &\le \norm{e^{t \Lambda }\mathbf{u}_{T,\epsilon}}_{L^{2}}
        +\norm{\qty(e^{t\mathscr{L}}-e^{t \Lambda })\mathbf{u}_{T,\epsilon}}_{L^{2}}\\
        &= e^{t \Lambda }\norm{\mathbf{u}_{T,\epsilon}}_{L^{2}}
        +\norm{\qty(e^{t \Lambda }-e^{t\mathscr{L}})\mathbf{u}_{T,\epsilon}}_{L^{2}}\\
        &\le \qty(e^{t \Lambda }+2\epsilon)\norm{\mathbf{u}_{T,\epsilon}}_{L^{2}},
    \end{aligned}
\end{equation}
where the last step is \eqref{linear_analysis_7} with \(2\epsilon\) in place of \(\epsilon\).
The reverse triangle inequality and the same bound give
\begin{equation}
    \norm{\mathbf{u}^{T,\epsilon}(t)}_{L^{2}}
    \ge \qty(e^{t \Lambda }-2\epsilon)\norm{\mathbf{u}_{T,\epsilon}}_{L^{2}}.
\end{equation}
Renaming \(2\epsilon\) as \(\epsilon\) yields \autoref{thm:Linear_Instability}(i).

\medskip\noindent\textbf{Part (ii).}
By part~(iv) of \autoref{lemma:linear_instability}, \(\supp\widehat{\mathbf{u}}_{T,\epsilon}\subseteq B(0,1)\), hence \(\abs{\bm{\xi}}\le 1\) on \(\supp\widehat{\mathbf{u}}_{T,\epsilon}\).
For every \(\abs{\bm{\alpha}}\le k\),
\begin{equation}
    \norm{\bm{\xi}^{\bm{\alpha}}\widehat{\mathbf{u}}_{T,\epsilon}}_{L^{2}}
    \le \norm{\widehat{\mathbf{u}}_{T,\epsilon}}_{L^{2}}.
\end{equation}
By Plancherel, since \(\mathbf{u}^{T,\epsilon}(0)=\mathbf{u}_{T,\epsilon}\),
\begin{equation}
    \begin{aligned}
        \norm{\mathbf{u}^{T,\epsilon}(0)}_{H^{k}}^{2}
        &=\sum_{\abs{\bm{\alpha}}\le k}\norm{\bm{\xi}^{\bm{\alpha}}\widehat{\mathbf{u}}_{T,\epsilon}}_{L^{2}}^{2}
        \le \binom{k+3}{3}\norm{\mathbf{u}^{T,\epsilon}(0)}_{L^{2}}^{2}.
    \end{aligned}
\end{equation}
Setting \(C_{k}\coloneqq\sqrt{\binom{k+3}{3}}\), one has
\begin{equation}
    \norm{\mathbf{u}^{T,\epsilon}(0)}_{H^{k}}
    \le C_{k}\norm{\mathbf{u}^{T,\epsilon}(0)}_{L^{2}},
\end{equation}
with \(C_{k}\) depending only on \(k\) and independent of \(T\) and \(\epsilon\).
This is \autoref{thm:Linear_Instability}(ii).
The proof is now complete.

\section{The nonlinear instability}\label{sec_nonlinear_instability}

This section is devoted to proving \autoref{thm:nonlinear_instability_Hadamard}. The constants in this section are denoted by \( \mathsf{C}_{j} \) which depend on \( \nu \) and \( f \) unless otherwise specified. The symbol \( C \) will represent a generic positive constant depending on those parameters, which need not be labeled and may vary from line to line.
We use the energy notation
\begin{equation}
    \mathcal{E}(\mathbf{u}(t))\coloneqq\norm{\mathbf{u}(t)}_{H^{1}},\qquad
    \mathcal{E}_{\ini}^{\kappa}\coloneqq\norm{\mathbf{u}^{\kappa}(0)}_{H^{1}},
\end{equation}
as in Appendix~\ref{sec:nonlinear_energy_estimate_v2}.

By \autoref{thm:Linear_Instability} (i), for any given parameter set \( \kappa = ( T, \varepsilon_{\mathrm{lin}}) \) with \(  T  > 0 \) and \( \varepsilon_{\mathrm{lin}} > 0 \) one can construct the unstable solution \( (\mathbf{u}^{ \kappa},q^{ \kappa}) \) for the linearized system \eqref{linearized_perturbed_equation} satisfying 
\begin{equation}\label{growth_estimate}
    (e^{t \Lambda} - \varepsilon_{\mathrm{lin}}) \norm{\mathbf{u}^{ \kappa}(0)}_{L^2} \leq   \norm{\mathbf{u}^{ \kappa}(t)}_{L^2}\leq (e^{t \Lambda} + \varepsilon_{\mathrm{lin}}) \norm{\mathbf{u}^{ \kappa}(0)}_{L^2} 
\end{equation}
for any \( t\in [0, T ] \). We denote by \( \mathbf{u}_{\ini}^{\kappa} \) the initial data of \( \mathbf{u}^{ \kappa} \). In view of (ii) of \autoref{thm:Linear_Instability},  there exists a constant $L>0$ independent of \( \kappa \) such that 
\begin{equation} \label{def_of_L}
\frac{\norm{\mathbf{u}_{\ini}^{ \kappa}}_{L^2}}{\mathcal{E}_{\ini}^{\kappa}}  \geq  L .
\end{equation}  

Now, define the initial data with a rescaling factor \( \delta \in (0, \overline{\delta}_0) \) by
\[
    \mathbf{u} ^{\delta}_{\ini} \coloneqq \frac{\delta}{\mathcal{E}_{\ini}^{\kappa}} \mathbf{u} _{\ini}^{\kappa}.
\]  
The parameter \( \overline{\delta}_0 \) is small enough such that \autoref{prop:energy_estimate_for_u} holds. Clearly, the size of the initial value in \( H^1 \) satisfies
\begin{equation}
\begin{aligned}
     \mathcal{E}(\mathbf{u}^{\delta}_{\ini}) = \delta < 1 .
\end{aligned}
\end{equation} 
In addition, we denote the approximate solution \( \mathbf{u}^{\app} \) with the associated pressure by
\begin{equation}
    \label{Hadamard_estimate_3}
\begin{aligned}
    \mathbf{u}^{\app}(\mathbf{x},t) \coloneqq \frac{\delta}{\mathcal{E}_{\ini}^{\kappa}}  \mathbf{u}^{ \kappa}(\mathbf{x},t) ,\quad 
    q  ^{\app}(\mathbf{x},t)  \coloneqq\frac{\delta}{\mathcal{E}_{\ini}^{\kappa}} q^{ \kappa}(\mathbf{x},t),
\end{aligned}
\end{equation} 
Hence \(\mathbf{u}^{\app}\) satisfies the same relative \(L^{2}\) bounds as in \eqref{growth_estimate}, with \(\mathbf{u}^{\kappa}\) replaced by \(\mathbf{u}^{\app}\).

Here the parameters are $\kappa=(T,\varepsilon_{\mathrm{lin}})$ fixing the linear pseudo-mode and $\delta$ fixing the initial size. To make the bootstrap work we fix $\delta\in(0,\overline{\delta}_{0})$ and choose $\kappa=(T,\varepsilon_{\mathrm{lin}})$ satisfying
\begin{equation} \label{fixed_parameters}
    \begin{aligned}
         \delta^{-4}\lesssim T \leq T_{0},\quad \varepsilon_{\mathrm{lin}} < \tfrac12\min\qty{\varepsilon,1},\quad 
        \tfrac{\delta}{2}e^{\Lambda T/2} \geq \mathsf{C}_{0},
    \end{aligned}
\end{equation}
where $\varepsilon=\varepsilon(\overline{\delta}_{0},f,\nu)$ is the nonlinear threshold fixed later in \eqref{def_of_epsilon_0}, $\mathsf{C}_{0}=\mathsf{C}_{0}(\nu,f)>0$, and
\begin{equation}
    T_{0}=C_{0}\nu^{3}/\|\nabla\mathbf{u}_{\ini}^{\delta}\|_{L^{2}}^{4}
\end{equation}
is the maximal time of the energy estimate \eqref{maximal_time_of_energy_estimate}. The three conditions have a clear meaning:
\begin{itemize}
    \item[(i)] $T\gtrsim\delta^{-4}$ ensures the linear growth has enough time to drive the nonlinear solution to an $O(1)$ size before $T_{0}$ expires;
    \item[(ii)] $\varepsilon_{\mathrm{lin}}$ small ensures the pseudo-mode is $L^{2}$-close to $e^{t\Lambda}$ on $[0,T]$;
    \item[(iii)] $(\delta/2)e^{\Lambda T/2}\ge\mathsf{C}_{0}$ guarantees $\varepsilon$ can be chosen independently of $\delta$, otherwise the instability threshold would degenerate as $\delta\to0$.
\end{itemize} 

The compatibility of $T$, $T_{0}$ and the escape time \( T^{\delta} \) (defined in \eqref{time_of_instability}) is as follows. By \eqref{maximal_time_of_energy_estimate} there exists $K=K(\nu)>0$ such that
\begin{equation}
    T_{0}\ge K\delta^{-4}.
\end{equation}
Hence with the natural choice $T=K\delta^{-4}$, the third inequality in \eqref{fixed_parameters} reads
\begin{equation}
    (\delta/2)\exp(\Lambda K/2\delta^{4})\ge\mathsf{C}_{0},\quad\text{i.e.}\quad \Lambda\ge2\delta^{4}K^{-1}\ln(2\mathsf{C}_{0}/\delta).
\end{equation}
Since $\delta\le1$, one can pick $\mathsf{C}_{0}=\mathsf{C}_{0}(\nu,f)$ satisfying
\begin{equation}
    \Lambda=\sup_{\delta\in(0,1]}2\delta^{4}K^{-1}\ln(2\mathsf{C}_{0}/\delta),
\end{equation}
which is possible for $\Lambda>0$. Then
\begin{equation}
    T^{\delta}\le T/2<T\le T_{0},
\end{equation}
so the bootstrap interval $[0,T^{\delta}]$ lies strictly inside the domain of linear approximation and of the energy estimate. In particular the family $\kappa(\delta)=(T(\delta),\varepsilon_{\mathrm{lin}}(\delta))$ depends on $\delta$, whence $\mathbf{u}_{\ini}^{\delta}$ depends on $\delta$, as allowed by \autoref{def_linear_instability} and \autoref{thm:nonlinear_instability_Hadamard}.

Let \( \mathbf{u}^{\delta} \) be a local nonlinear strong solution to \eqref{nonlinear_perturbed_equation}, emanating from \( \mathbf{u} ^{\delta}_{\ini} \), with the associated pressure \( q ^{\delta} \).  Define the time of instability
\begin{equation}\label{time_of_instability}
\begin{aligned}
    T^{\delta} \coloneqq \frac{1}{\Lambda} \ln \frac{2\varepsilon}{\delta}, \qq{i.e.} \delta e^{\Lambda T^{\delta}} = 2 \varepsilon,
\end{aligned}
\end{equation}
and the escape times 
\begin{equation}\label{escape_times}
\begin{aligned}
    T^{*} & \coloneqq  \sup \qty{ t > 0\ \middle|\ \mathcal{E}(\mathbf{u}^{\delta}(t)) \leq \overline{\delta}_0 }, \\ 
    T^{* *} & \coloneqq \sup \qty{ t > 0\ \middle|\ \norm{\mathbf{u}^{\delta}(t)}_{L^2}   \leq 2 \delta e ^{\Lambda t}}.
\end{aligned}
\end{equation}
Then, \( \mathbf{u}^{\delta} \) satisfies \eqref{small_energy_assumption} for \( T^{*} \) and \( \overline{\delta}_0 \). Thus, recalling the estimate \eqref{estimate_with_small_energy_assumption_for_Hadamard_new}, one obtains that 
\begin{equation} \label{Hadamard_estimate_1}
    \begin{aligned} 
    &\norm{( \mathbf{u}^{\delta},\nabla \mathbf{u}^{\delta})(t)}_{L^2}^2 + \int_{0}^{t} \norm{(\nabla \mathbf{u}^{\delta}, \Delta \mathbf{u}^{\delta})(\tau)}_{L^2}^2 \dd{\tau}\\& \leq  C
    \qty(\mathcal{E}^2(\mathbf{u}^{\delta}_{\ini }) + \int_{0}^{t} \|\mathbf{u}^{\delta}(\tau)\|_{L^2}^2 \dd{\tau} ) \\ 
    &\leq    C \qty( \delta^2 + \frac{2}{\Lambda} \delta^2 e ^{2\Lambda t}) \leq \mathsf{C}_1 \delta^2 e^{2\Lambda t},
\end{aligned}
\end{equation} 
for any \( t \leq \min \qty{T^{\delta}, T^{*}, T^{* *} } \), where \( \mathsf{C}_1 \) does not depend on \( \delta \).

Subsequently, we denote the difference of nonlinear solution \( \mathbf{u}^{\delta} \) and the approximate solution by
\[
     \mathbf{u}^{\dif} =  \mathbf{u}^{\delta} -\mathbf{u}^{\app}  ,\quad 
      q ^{\dif}=   q ^{\delta} -  q ^{\app} .
\] 
The nonlinear solution has the integral expression
\begin{equation}
    \mathbf{u}^{\delta}(t) = e^{t\mathscr{L}} \mathbf{u}_{\ini}^{\delta} + \int_{0}^{t} e^{(t - s)\mathscr{L}} \mathscr{N}(\mathbf{u}^{\delta} (s)) \dd{s} ,
\end{equation}
where the nonlinear term \( \mathscr{N} \) is defined as
\begin{equation}
    \mathscr{N}(\mathbf{u}) \coloneqq \nabla \Delta ^{-1} \div \qty(\mathbf{u} \cdot \nabla \mathbf{u}) -\mathbf{u} \cdot \nabla \mathbf{u}.
\end{equation}
Since the approximate solution satisfies \(  \mathbf{u}^{\app} = e^{t\mathscr{L}} \mathbf{u}_{\ini}^{\delta} \), then the difference admits the integral representation
\begin{equation}
     \mathbf{u}^{\dif}(t) = \int_{0}^{t} e^{(t - s)\mathscr{L}} \mathscr{N}(\mathbf{u}^{\delta} (s)) \dd{s} .
\end{equation}
For the nonlinear term, we have the following estimate using Gagliardo-Nirenberg inequality 
\begin{equation}
    \norm{\mathscr{N}(\mathbf{u}) }_{L^2} \leq \norm{ \mathbf{u} \cdot \nabla \mathbf{u}}_{L^2}   
    \leq \norm{\mathbf{u}^{\delta} }_{L^{3}}  \norm{  \nabla  \mathbf{u}^{\delta} }_{L^6}
    \leq \norm{\mathbf{u}^{\delta} }_{L^{2}} ^{\frac{1}{2}} \norm{  \nabla ^2 \mathbf{u}^{\delta} }_{L^2}^{\frac{3}{2}}  .
\end{equation}
In light of this estimate, \eqref{eq:simp_semigroup_bnd} and \eqref{Hadamard_estimate_1}, fix \(\varepsilon_1\in(0,\frac{\Lambda}{2})\). We have 
\begin{equation}\label{Hadamard_estimate_4}
    \begin{aligned}
	\norm{ \mathbf{u}^{\dif}(t)}_{L^{2}} &  \leq C_{\varepsilon_1} \int_{0}^{t} e^{(t - s)(\Lambda+\varepsilon_1)} \norm{\mathbf{u}^{\delta}(s) }_{L^{2}} ^{\frac{1}{2}} \norm{  \nabla ^2 \mathbf{u}^{\delta}(s) }_{L^2}^{\frac{3}{2}}    \dd{s}
    \\&\leq C_{\varepsilon_1} \int_{0}^{t} e^{(t - s)(\Lambda+\varepsilon_1)} \delta^{\frac{1}{2}} e^{ \frac{1}{2} s \Lambda} \norm{  \nabla ^2 \mathbf{u}^{\delta} (s)}_{L^2}^{\frac{3}{2}}    \dd{s} 
    \\ 
    & \leq C_{\varepsilon_1} \delta^{\frac{1}{2}} e^{t (\Lambda+\varepsilon_1)}  \int_{0}^{t} e^{ s \qty( \frac{1}{2}  \Lambda - \frac{1}{4} (\Lambda+\varepsilon_1))} e^{ - \frac{3}{2} \cdot \frac{1}{2}(\Lambda+\varepsilon_1) s }  \norm{  \nabla ^2 \mathbf{u}^{\delta}(s) }_{L^2}^{\frac{3}{2}}    \dd{s} \\
 &   \leq C_{\varepsilon_1} \delta^{\frac{1}{2}} e^{t (\Lambda+\varepsilon_1)} \qty(\int_{0}^{t} e^{ s  (\Lambda - \varepsilon_1)} \dd{s})^{\frac{1}{4}} \qty( \int_{0}^{t}    e^{ - (\Lambda+\varepsilon_1) s }\norm{ \nabla ^2 \mathbf{u}^{\delta} (s) }_{L^2}^{2}    \dd{s} )^{\frac{3}{4}}. 
\end{aligned}   
\end{equation}
By using \autoref{lemma_ana_tool_1} and \eqref{Hadamard_estimate_1}, we have 
\begin{equation}
    \int_{0}^{t} e^{ - (\Lambda+\varepsilon_1) s } \norm{  \nabla ^2 \mathbf{u}^{\delta} (s) }_{L^2}^{2}    \dd{s} \leq C \delta^2 e^{(\Lambda - \varepsilon_1)t},
\end{equation}
 \eqref{Hadamard_estimate_4} then becomes
\begin{equation} 
    \begin{aligned}
	\norm{ \mathbf{u}^{\dif}(t)}_{L^{2}} &  \leq C \delta^{2} e^{t (\Lambda+\varepsilon_1)} \qty( e^{t  (\Lambda - \varepsilon_1)} )^{\frac{1}{4}} \qty( e^{t (\Lambda - \varepsilon_1)} )^{\frac{3}{4}} 
     \leq \mathsf{C}_2 \delta^{2}   e^{ 2 \Lambda t }  \qfor t \leq \min \qty{T^{\delta}, T^{*}, T^{* *} },
    \\  
\end{aligned}   
\end{equation}
since \(\varepsilon_1\in(0,\frac{\Lambda}{2})\) ensures \(\Lambda-\varepsilon_1>0\).

Next, we show that 
\begin{equation}\label{Hadamard_estimate_2}
\begin{aligned}
    T^{\delta} = \min  \qty{T^{\delta},T^{*}, T^{* *},  T },
\end{aligned}
\end{equation}
with \( \varepsilon \) defined as 
\begin{equation} \label{def_of_epsilon_0}
    \varepsilon \coloneqq \min \qty{
        \frac{\overline{\delta}_0}{4\sqrt{\mathsf{C}_1}} , \frac{1}{ 4 \mathsf{C}_2}, \mathsf{C}_0, \frac{L}{ 8 \mathsf{C}_2}
    }.
\end{equation}
If \( T^{*} < T^{\delta} \), then we have 
\begin{equation}
    \overline{\delta}_0 = \mathcal{E}(\mathbf{u}^{\delta}(T^{*})) \leq \sqrt{\mathsf{C}_1} \delta e^{\Lambda T^{*}} \leq \sqrt{\mathsf{C}_1} \delta e^{\Lambda T^{\delta}} = 2\sqrt{\mathsf{C}_1} \varepsilon \leq \frac{\overline{\delta}_0}{2}  < \overline{\delta}_0,
\end{equation}
which is a contradiction.
If \( T^{* *} < T^{\delta} \), then for \( \varepsilon_{\mathrm{lin}} <  \frac{1}{2} \), one gets
\begin{equation}
   \begin{aligned}
	 \norm{\mathbf{u}^{\delta}(T^{* *})}_{L^2} \leq {} &\norm{\mathbf{u}^{\app}(T^{* *})}_{L^2} + \norm{\mathbf{u}^{\dif}(T^{* *})}_{L^2}
      \\\leq  {} & \delta  (e ^{\Lambda T^{* *}} + \varepsilon_{\mathrm{lin}}) + \mathsf{C}_2 \delta^2  e^{ 2\Lambda  T^{* *}}   \\ 
      \leq {} & \delta e ^{\Lambda T^{* *}} (1 +  2\mathsf{C}_2 \varepsilon ) + \delta \varepsilon_{\mathrm{lin}}  \leq \frac{3}{2}\delta e^{\Lambda T^{* *}}+ \delta \varepsilon_{\mathrm{lin}} < 2\delta e^{\Lambda T^{* *}},
\end{aligned}
\end{equation}
which contradicts \eqref{escape_times}.  By \eqref{fixed_parameters} and \eqref{def_of_epsilon_0}, we have 
\begin{equation}
    T^{\delta} = \frac{1}{\Lambda} \ln \frac{2\varepsilon}{\delta} \leq \frac{ T }{2} <  T .
\end{equation}
Thus, \(  T  > T^{\delta} \).  \eqref{Hadamard_estimate_2} is now verified.

Finally, we show that the velocity is unstable in the $L^2$ norm. Recalling \eqref{Hadamard_estimate_3}, \eqref{fixed_parameters} and \eqref{def_of_epsilon_0}, we have 
\begin{equation}
\begin{aligned}
    \norm{\mathbf{u}^{\delta}(T^{\delta})}_{L^2} \geq {} & \norm{\mathbf{u}^{\app}(T^{\delta})}_{L^2} - \norm{\mathbf{u}^{\dif}(T^{\delta})}_{L^2} \\ \geq {} &\frac{\delta}{\mathcal{E}_{\ini}^{\kappa}} \norm{   \mathbf{u}^{ \kappa}(T^{\delta})}_{L^2} - \mathsf{C}_2 \delta^{2}   e^{ 2 \Lambda T^{\delta} } \\ 
    \geq {} & \delta (e^{T^{\delta} \Lambda } - \varepsilon_{\mathrm{lin}}) \frac{\norm{\mathbf{u}_{\ini}^{\kappa}}_{L^2}}{\mathcal{E}_{\ini}^{\kappa}} -\mathsf{C}_2 \delta^{2}   e^{ 2 \Lambda T^{\delta} }\\
    \geq {} & 2 \varepsilon   \qty(L   - 2 \mathsf{C}_2  \varepsilon) - L \varepsilon_{\mathrm{lin}} \delta
    \geq \frac{3}{2}L \varepsilon - L\varepsilon_{\mathrm{lin}} \delta   > L  \varepsilon, \\ 
\end{aligned}
\end{equation}
where we used the fact from \eqref{fixed_parameters}:
\begin{equation} 
    \begin{aligned}
        \varepsilon_{\mathrm{lin}} \delta  \leq  \varepsilon_{\mathrm{lin}}  <  \frac{1}{2} \varepsilon.
    \end{aligned}
\end{equation}  
By setting \( \varepsilon_0 = L\varepsilon \), one has \( T^{\delta}=\frac1\Lambda\ln\frac{2\varepsilon}{\delta}=\frac1\Lambda\ln\frac{2\varepsilon_0}{L\delta} \), and the proof of \autoref{thm:nonlinear_instability_Hadamard} is now complete.

\begin{appendices}
\section{Appendices}
\subsection{Nonlinear energy estimates}\label{sec:nonlinear_energy_estimate_v2}

In this section, we derive some nonlinear energy estimates for the perturbed problem, which are  used in the proof of the nonlinear instability, cf. \cite{Bedrossian2022}. To this end, let us assume that \( \psol \) is a sufficiently regular solution  in \( [0,T]\times \mathbb{R}^3 \) for some \( T > 0\) to the perturbed system \eqref{nonlinear_perturbed_equation}. 

Hereafter, we define
\begin{equation}
    \mathcal{E} (\mathbf{u}(t)) \coloneqq \norm{\mathbf{u}(t)}_{H^1}  ,\quad  \mathcal{E}_{\ini} \coloneqq  \norm{\mathbf{u}(0)}_{H^1}=  \norm{\mathbf{u}_{\ini}}_{H^1}.
\end{equation} 
and assume 
\begin{equation}\label{small_energy_assumption}
    \mathcal{E}(\mathbf{u}(t))\leq \overline{\delta}_0 \leq  1,\qfor t \in[0,T] .
\end{equation} 

One has 
\begin{equation}\label{energy_estimate_1_new}
	\begin{aligned}
		\frac{1}{2} \dv{}{t} \norm{\mathbf{u}}_{L^2}^2 + \nu \norm{\nabla \mathbf{u}}_{L^2}^2 =  -\int_{\mathbb{R}^3}^{} u_1 u_2 \dd{\mathbf{x}} \leq  \norm{u_1}_{L^2}\norm{u_2}_{L^2} \leq  \norm{\mathbf{u}}_{L^2}^2,
	\end{aligned}
\end{equation}
and
\begin{equation}
    \begin{aligned}
        &\frac{1}{2}  \dv{}{t}  \|\nabla \mathbf{u}\|_{L^2}^2 + \nu \|\Delta \mathbf{u}\|_{L^2}^2 \\ 
        &= \int_{\mathbb{R}^3} \nabla q \cdot \Delta \mathbf{u} \dd{\mathbf{x}} + \int_{\mathbb{R}^3} ( \mathbf{u} \cdot \nabla) \mathbf{u} \cdot \Delta \mathbf{u} \dd{\mathbf{x}}   \\&\quad+ \int_{\mathbb{R}^3}^{} u_2 \Delta u_1 \dd{\mathbf{x}} 
      + \int_{\mathbb{R}^3}^{}  y \partial_x \mathbf{u} \cdot \Delta   \mathbf{u} \dd{\mathbf{x}} 
         \\
        &= - \int_{\mathbb{R}^3} q \Delta(\nabla \cdot \mathbf{u}) \dd{\mathbf{x}} - \int_{\mathbb{R}^3} \nabla(( \mathbf{u} \cdot \nabla) \mathbf{u} ): \nabla \mathbf{u} \dd{\mathbf{x}} 
         \\&\quad - \int_{\mathbb{R}^3}^{} \nabla u_2 \cdot  \nabla u_1 \dd{\mathbf{x}} 
        - \int_{\mathbb{R}^3}^{} \partial_x \mathbf{u}\cdot  \partial_y \mathbf{u} \dd{\mathbf{x}} \\
        &= - \int_{\mathbb{R}^3}  \mathbf{u} \cdot \nabla \frac{|\nabla \mathbf{u}|^2}{2} \dd{\mathbf{x}} - \int_{\mathbb{R}^3} 
        (\nabla \mathbf{u} \cdot (\nabla \mathbf{u})^{\intercal} ): \nabla \mathbf{u}\dd{\mathbf{x}}   \\&\quad- \int_{\mathbb{R}^3}^{} \nabla u_2 \cdot  \nabla u_1 \dd{\mathbf{x}} 
        - \int_{\mathbb{R}^3}^{} \partial_x \mathbf{u}\cdot  \partial_y \mathbf{u} \dd{\mathbf{x}} \\
        &\leq \|\nabla \mathbf{u}\|_{L^3}^3 + 2\|\nabla \mathbf{u}\|_{L^2}^2
      \leq C \|\nabla \mathbf{u}\|_{L^2}^{3/2} \|\Delta \mathbf{u}\|_{L^2}^{3/2} + 2\|\nabla \mathbf{u}\|_{L^2}^2 \\ 
      &\leq \frac{\nu}{2} \|\Delta \mathbf{u}\|_{L^2}^2 + \frac{C}{\nu^3} \|\nabla \mathbf{u}\|_{L^2}^6 + 2\|\nabla \mathbf{u}\|_{L^2}^2.
        \end{aligned}
\end{equation}
Hence 
\begin{equation}\label{energy_estimate_3_new}
    \frac{1}{2}  \dv{}{t}   \|\nabla \mathbf{u}\|_{L^2}^2 + \frac{\nu}{2} \|\Delta \mathbf{u}\|_{L^2}^2 \leq  \frac{C}{\nu^3} \|\nabla \mathbf{u}\|_{L^2}^6 + 2\|\nabla \mathbf{u}\|_{L^2}^2 .
\end{equation}
Combining \eqref{energy_estimate_1_new} and \eqref{energy_estimate_3_new}, we have 
\begin{equation}
     \dv{}{t}  \norm{( \mathbf{u},\nabla \mathbf{u})}_{L^2}^2 +  \norm{(\nabla \mathbf{u}, \Delta \mathbf{u})}_{L^2}^2 \leq  C \qty( \norm{\nabla \mathbf{u}} _{L^2}^6 
    + \|\mathbf{u}\|_{L^2}^2).
\end{equation}
For small enough \( \overline{\delta}_0 \), the term \( \norm{\nabla \mathbf{u}} _{L^2}^6  \) on the right-hand side can be absorbed by the term \( \norm{\nabla \mathbf{u}} _{L^2}^2 \) of the left-hand side. One has
\begin{proposition}\label{prop:energy_estimate_for_u}
    There exists a constant \( \overline{\delta}_0 \leq 1 \), such that if a strong solution \( \mathbf{u}(t) \) of the system \eqref{nonlinear_perturbed_equation} satisfies the assumption \eqref{small_energy_assumption} for some \( T > 0 \), then the solution satisfies 
    \begin{equation} \label{estimate_with_small_energy_assumption_for_Hadamard_new}
        \begin{aligned}
            &  \norm{( \mathbf{u},\nabla \mathbf{u})(t)}_{L^2}^2 + \int_{0}^{t} \norm{(\nabla \mathbf{u}, \Delta \mathbf{u})(\tau)}_{L^2}^2 \dd{\tau} \leq  C
            \qty(\mathcal{E}^2_{\ini}+ \int_{0}^{t} \|\mathbf{u}(\tau)\|_{L^2}^2 \dd{\tau} ),
    \end{aligned}
    \end{equation}
    for all \( t\in [0,T] \), where the constant \( C \) depends on \( \nu \).
\end{proposition}

    \subsection{Existence of solutions}

 In view of the energy method in \cite{Bedrossian2017,Majda2002,Bedrossian2022} we have the existence of the solution to the system \eqref{nonlinear_perturbed_equation} in the divergence-free spaces introduced in Section~\ref{sec_Preliminary}. We then have the following theorem.

    \begin{theorem}[Local existence]\label{thm:local_existence}
        Fix $\nu>0$ and $\mathbf{u}_{\ini}\in H^1_\sigma(\mathbb{R}^3) := L^2_\sigma(\mathbb{R}^3)\cap H^1(\mathbb{R}^3)$. Then there exists a constant \( C_0 > 0 \) such that 
        \begin{equation}
            \label{maximal_time_of_energy_estimate}
            T_0  = C_0\frac{\nu^3}{\|\nabla \mathbf{u}_{\ini}\|_{L^2}^4}> 0
        \end{equation}
         and a unique solution $\mathbf{u}\in C([0,T_0);H^1_\sigma)\cap L^2((0,T_0);H^2_\sigma)$ of the Cauchy problem for the Navier-Stokes equation \eqref{nonlinear_perturbed_equation}, with initial value $\mathbf{u}_{\ini}$. 
        \end{theorem}

        \subsection{A technical lemma}
        
        \begin{lemma} \label{lemma_ana_tool_1}
Let \(a > b > 0\), \( T > 0 \) and \(K > 0\) be constants, and let \(F \in L^2([0,T])\) be a function such that for all \(0 \leq t \leq T\),  
\[
\int_0^t |F(s)|^2 \dd{s} \leq K e^{2at}.
\]  
Then there exists a constant \(C = \frac{aK}{a - b}\) such that for all \(0 \leq t \leq T\),  
\[
\int_0^t |F(s)|^2 e^{-2bs} \dd{s} \leq C e^{2(a - b)t}.
\]  
        \end{lemma}
        \begin{proof}
Define \(G(t) = \int_0^t |F(s)|^2 \dd{s}\), so \(G(t) \leq K e^{2at}\) for all \(0 \leq t \leq T\), and \(G'(s) = |F(s)|^2\) almost everywhere.

To estimate \(\int_0^t |F(s)|^2 e^{-2bs} \dd{s}\), we use integration by parts, which results in
\begin{equation}
    \begin{aligned}
        \int_0^t |F(s)|^2 e^{-2bs} \dd{s} & =  \left. \left[ G(s) e^{-2bs} \right] \right|_0^t + 2b \int_0^t G(s) e^{-2bs} \dd{s} \\ 
 & = G(t) e^{-2bt} + 2b \int_0^t G(s) e^{-2bs} \dd{s} \\
& \leq K e^{2(a -b)t} +  2b \int_0^t K e^{2as} \cdot e^{-2bs} \dd{s} \\
& \leq K e^{2(a -b)t} + \frac{bK}{a - b} \left( e^{2(a - b)t} - 1 \right).
    \end{aligned}
\end{equation}
 Using the condition \( a > b > 0 \), one has 
 \begin{equation}
    \begin{aligned}
        \int_0^t |F(s)|^2 e^{-2bs} \dd{s} & \leq   \frac{aK}{a - b} e^{2(a - b)t} .
    \end{aligned}
\end{equation} 
        \end{proof}

\end{appendices}

\section*{Acknowledgement} 
The work of Q.Wang is supported by the National Natural Science Foundation of Sichuan
Province (No.2025ZNSFSC0072) and the National Science Foundation of China (No.12301131
and No. 11901408).

\itemsep=0pt


\end{document}